\documentclass[oneside]{amsart}
\usepackage{amsmath,amssymb,color,amsthm,graphicx,cite}
\usepackage{color,bm}
\usepackage{enumerate} 
\usepackage[T1]{fontenc}
\usepackage[utf8]{inputenc} 
\usepackage[backref]{hyperref}
\RequirePackage{luatex85}
\usepackage[all]{xy}
\newtheorem{theorem}{Theorem}[section]
\newtheorem{proposition}[theorem]{Proposition}
\newtheorem{lemma}[theorem]{Lemma}
\newtheorem{corollary}[theorem]{Corollary}
\newtheorem{claim}{Claim}
\newtheorem{remark}{Remark}
\newtheorem{question}{Question}

\newtheorem{problem}{Problem}
\theoremstyle{definition}
\newtheorem{definition}{Definition}
\newtheorem{main}{Theorem}

\newtheorem{main_cor}[main]{Corollary}

\def\Z{\mathbb{Z} }
\def\N{\mathbb{N} }
\def\Q{\mathbb{Q} }
\def\R{\mathbb{R} }
\def\A{\mathbb{A} }
\def\T{\mathbb{T} }

\def\nbd{neighborhood }
\def\nbds{neighborhoods }

\def\Sv{\mathop{\mathrm{Sing}}(v)}

\def\Pv{\mathop{\mathrm{Per}}(v)}
\def\Cv{\mathop{\mathrm{Cl}}(v)}

\author{Tomoo Yokoyama}
\date{\today}
\address{Department of Mathematics, Faculty of Science, Saitama University, Shimo-Okubo 255, Sakura-ku, Saitama-shi, 338-8570 Japan\\}
\email{tyokoyama@rimath.saitama-u.ac.jp}
\thanks{The author was partially supported by JSPS Grant Number 20K03583 and 24K06733}
\subjclass[2010]{}

\makeatletter
\@namedef{subjclassname@2020}{%
  \textup{2020} Mathematics Subject Classification}
\makeatother

\title[Beyond Poincar{\'e} recurrence]{Beyond Poincar{\'e} recurrence: on a topological trichotomy of orbits and on the dependence of limit sets for flows on surfaces}
\date{\today}
\keywords{Flows on surfaces, recurrence, non-wandering sets, minimal flows, non-limit ciruits}

\makeatletter
\@namedef{subjclassname@2020}{%
  \textup{2020} Mathematics Subject Classification}
\makeatother

\subjclass[2020]{37E35,37B20,37G30,37C55}

\begin{document}

\maketitle

\begin{abstract}
This paper investigates the global structure of orbits for continuous flows on surfaces. One of our main results provides a topological trichotomy of orbits of flows with finitely many connected components of the singular point set on surfaces with finite genus and finite ends: each orbit is either recurrent, wandering, or corresponds to a circuit with wandering holonomy or a strictly limit circuit, two new concepts that we introduce in this paper. These concepts capture non-recurrent behaviors lying beyond the reach of classical tools such as the Poincar{\'e} recurrence theorem and the Poincar{\'e}-Bendixson theorem. Moreover, the non-recurrent non-wandering behaviors reveal the dependency between the existence of infinitely many non-degenerate $\omega$-limit sets and $\alpha$-limit sets, and provide a necessary and sufficient condition for a flow to be minimal on (possibly non-compact) surfaces.
\end{abstract}

\section{Introduction}

Birkhoff introduced the concepts of non-wandering points and recurrent points in \cite{birkhoff1927dynamical}, and using these concepts, we can describe and capture dynamical behaviors. 
Two classical results assure the existence of non-wandering orbits: the Poincar{\'e} recurrence theorem and the Poincar{\'e}-Bendixson theorem.
For instance, the existence of returns induced by the Poincar{\'e} recurrence theorem can be roughly described as the denseness of recurrent orbits in the non-wandering set for measure-preserving dynamical systems on finite measure spaces. 
On the other hand, the existence of non-wandering orbits outside the closure of the union of recurrent orbits is an obstruction to the denseness of recurrence as seen in the Poincar{\'e} recurrence theorem. 
However, no characterization of non-recurrent but non-wandering orbits was known even for flows on surfaces.

The Poincar{\'e}-Bendixson theorem captures the limit behaviors of orbits of flows and has been applied to various phenomena (e.g. \cite{bhatia1966application,du2021traveling,koropecki2019poincare,hajek1968dynamical,roussarie2020topological,roussarie2021some,roussarie2021some02,pokrovskii2009corollary}). 
Nevertheless, there exist non-wandering orbits which neither the Poincar{\'e}-Bendixson theorem nor recurrence descriptions capture even for flows with finitely many singular points on closed surfaces. 
This is because such non-wandering orbits are not contained in the union of $\omega$-limit sets, $\alpha$-limit sets, and the closure of the recurrent point set. 
Therefore, describing such non-wandering orbits serves as a complement to the Poincar{\'e}-Bendixson theorem and the Poincar{\'e} recurrence theorem. 
In this paper, we consider the following problem.

\begin{problem}\label{prob:nw_nonrec}
Let $S$ be a surface and $v$ be a continuous flow on $S$.
Characterize non-wandering orbits of $v$ that do not intersect the union of the $\alpha$-limit sets of all points, the $\omega$-limit sets of all points, and the closure of the set of recurrent points.
\end{problem}

While Problem~\ref{prob:nw_nonrec} focuses on topological constraints on specific invariant subsets, other topological constraints are known to govern the relationship between limit sets. Specifically, there is a known dependency between the types of $\alpha$-limit sets and $\omega$-limit sets \cite{yokoyama2023dependency}. 
However, the existence and nature of other such constraints remain unclear. 
For instance, one can take the perspective of going ``beyond Poincar{\'e} recurrence'' to explore dependencies that classical recurrence cannot capture. 
From this viewpoint, it is natural to ask whether there are constraints on their cardinality, such as the existence of infinitely many limit sets. 
Motivated by this, we investigate the following question.


\begin{question}\label{q:dependency_a_w}
Does the existence of infinitely many pairwise distinct non-degenerate $\omega$-limit sets of non-closed {\rm(i.e.} non-singular and non-periodic{\rm)} orbits imply the existence of infinitely many pairwise distinct non-degenerate $\alpha$-limit sets of non-closed orbits?
\end{question}

Here, non-degenerate $\omega$-limit (resp. $\alpha$-limit) sets are $\omega$-limit (resp. $\alpha$-limit) sets which are not singletons. 

\subsection{Main results}

In this paper, we provide an answer to both Problem~\ref{prob:nw_nonrec} and Question~\ref{q:dependency_a_w} for flows on compact surfaces with finitely many singular points.

More specifically, to describe non-recurrent non-wandering orbits, we introduce two new concepts, a circuit with wandering holonomy and a strictly limit circuit, and establish a topological trichotomy of orbits for flows with finitely many singular points on compact surfaces.
We also classify non-recurrent non-wandering orbits and obtain the equivalence of the existence of infinitely many non-degenerate $\omega$-limit sets and non-degenerate $\alpha$-limit sets. 
As applications of these results, we characterize the minimality of flows on compact surfaces and the denseness of closed orbits in the non-wandering set, and refine the Poincar{\'e} recurrence theorem for flows on surfaces.

\subsubsection{Dependency between the existence of infinitely many positive and negative limit behaviors}

To introduce the first of our main results, we recall some concepts for a flow $v$ on a surface. 
An orbit is periodic if it is homeomorphic to a circle, and an orbit $O$ is {\bf recurrent} if $O \subseteq \alpha(O) \cup \omega(O)$, where $\alpha(O)$ is the $\alpha$-limit set of $O$ and $\omega(O)$ is the $\omega$-limit set of $O$. 
An immersed image of an oriented circle by a continuous orientation-preserving mapping is a {\bf non-trivial circuit} if it is either a periodic orbit or a directed graph which is not a singleton and is the union of finitely many singular points and non-singular orbits such that the orientation is compatible with the dynamical orientation of the orbits.
A non-trivial circuit is a {\bf limit circuit} if it is either the $\omega$-limit or $\alpha$-limit set of a non-recurrent point. 
We are now ready to state Theorem~\ref{main:inifinie_limit_sets}, which demonstrates that the answer to Question~\ref{q:dependency_a_w} is affirmative for flows with finitely many singular points. 

\begin{main}\label{main:inifinie_limit_sets}
The following statements are equivalent for a flow $v$ with finitely many singular points on a compact surface:
\\
{\rm(1)} There are infinitely many pairwise distinct $\omega$-limit sets of non-closed orbits.  
\\
{\rm(2)} There are infinitely many pairwise distinct $\omega$-limit sets of non-closed orbits which are limit circuits. 
\\
{\rm(3)} There are infinitely many pairwise distinct $\omega$-limit sets of non-closed orbits which are periodic orbits.
\\
{\rm(4)} There are infinitely many pairwise distinct $\alpha$-limit sets of non-closed orbits. 
\\
{\rm(5)} There are infinitely many pairwise distinct $\alpha$-limit sets of non-closed orbits which are limit circuits. 
\\
{\rm(6)} There are infinitely many pairwise distinct $\alpha$-limit sets of non-closed orbits which are periodic orbits. 
\end{main}

Moreover, we characterize the existence of finitely many $\omega$-limit sets. 

\begin{main}\label{main:ch_finite_limit_sets}
The following statements are equivalent for a flow $v$ on a compact surface:
\\
{\rm(1)} There are at most finitely many pairwise distinct $\omega$-limit sets of points. 
\\
{\rm(2)} There are at most finitely many pairwise distinct $\alpha$-limit sets of points. 
\\
{\rm(3)} There are at most finitely many closed orbits.
\end{main}

\subsubsection{Mechanism of occurrence of the non-recurrent non-wandering behaviors}
To describe non-recurrent non-wandering behaviors more precisely, we recall some concepts. 
An orbit is singular if it consists of a point, and is closed if it is singular or periodic. 
For a flow $v$, denote by $\Sv$ (resp. $\Pv$, $\Cv$, $\mathrm{P}(v)$, $\mathrm{R}(v)$) the union of singular (resp. periodic, closed, non-recurrent, non-closed recurrent) orbits. 

Furthermore, by introducing the concept of a circuit with wandering holonomy (see Definition~\ref{def:cwh} in the next section), which returns near itself exactly once and never comes back again, and the concept of a strictly limit circuit (see Definition~\ref{def:sls} in the next section), which either attracts all nearby points or allows half of them to escape, we characterize the mechanism of occurrence of the non-recurrent non-wandering behaviors as follows. 

\begin{main}\label{main:rec_wandering_pt}
Let $v$ be a flow with finitely many singular points on a compact surface. 
For any orbit $O$ of $v$, exactly one of the three following statements holds:
\\
{\rm1.} The orbit $O$ is recurrent. 
\\
{\rm2.} The orbit $O$ is wandering and $\overline{O} - O = \alpha(O) \cup \omega(O)$. 
\\
{\rm3.} The orbit $O$ is non-recurrent and non-wandering. 
Moreover, it satisfies $\overline{O} - O = \alpha(O) \cup \omega(O) \subseteq \Sv$ and one of the conditions {\rm(3.1)}, {\rm(3.2)}, {\rm(3.3)}, and {\rm(3.4)} holds:
\\
{\rm(3.1)} There is a strictly limit circuit in the non-wandering set $\Omega(v)$ at $O \subset \mathop{\mathrm{int}} \mathrm{P}(v)$. 
\\
{\rm(3.2)} There is a circuit with wandering holonomy at $O \subset \mathop{\mathrm{int}} \mathrm{P}(v)$. 
\\
{\rm(3.3)} There is a circuit that is contained in $\overline{\Pv} - \Pv \subseteq \Omega(v)$ and contains $O$. 
\\
{\rm(3.4)} $O \subset \overline{\mathrm{R}(v)} - \mathrm{R}(v)$.
\end{main}

The previous theorem follows from a more general result for flows with finitely many connected components of the singular point set on (possibly non-compact) surfaces with finite genus and finite ends (see Theorem~\ref{main:rec_wandering}). 
Note that several techniques used to prove this theorem, such as the classification of orbits (proper, locally dense, or exceptional (see the definition in \S~\ref{sec:tpo})), holonomy, and the waterfall construction (cf. \cite[Lemma~3.3.7, p.86]{candel2000foliation}), originate from foliation theory (cf. \cite{candel2000foliation}).

Applying Theorem~\ref{main:rec_wandering_pt}, we obtain the following characterization, which answers Problem~\ref{prob:nw_nonrec} for flows with finitely many singular points on compact surfaces. 

\begin{main_cor}\label{main_cor:nw_nonrec}
The following statements are equivalent for any non-wandering orbit $O$ of a flow with finitely many singular points on a compact surface $S$: 
\\
{\rm(1)} The orbit $O$ is contained in the complement of the union of the $\alpha$-limit sets of all points, the $\omega$-limit sets of all points, and the closure of the set of recurrent points {\rm(i.e.} $O \subseteq S - (\bigcup_{x \in S} \alpha(x) \cup \omega(x) \cup \overline{\mathrm{R}(v)})${\rm)}. 
\\
{\rm(2)} There is a circuit with wandering holonomy at $O$.
\end{main_cor}

Notice that the previous characterization does not hold for flows with countably many singular points on compact surfaces (see Example~\ref{lem:counter_exam01}).


\subsection{Applications of Theorems~\ref{main:inifinie_limit_sets} and \ref{main:rec_wandering_pt}}

The main results Theorem~\ref{main:inifinie_limit_sets} and Theorem~\ref{main:rec_wandering_pt} imply the following applications. 

\subsubsection{Characterization of minimal flows}

Minimality has been studied from various perspectives in topological dynamics (see books  \cite{Auslander1988mini,Ellis1969lecture,Gottschalk1955top_dyn} for more details). 
For instance, examples of minimal flows on surfaces are constructed in \cite{aranson1973invariant,gardiner1985structure,Sacker1972existence}, and properties of minimal flows on surfaces are described in \cite{Athanassopoulos1995minimal,Athanassopoulos1997stable,Basener2006min,Gottschalk1963minimal,Smith1988transitive,ulcigrai2011absence}. 
On the other hand, the existence of dense orbits of flows on compact surfaces is topologically characterized in \cite{Marzougui2009dense,yokoyama2016topological}. 
Similarly, we would like to find a topological characterization of minimality for flows on {\rm(}possibly non-compact{\rm)} surfaces. 
Therefore, we pose the following problem. 

\begin{problem}\label{prob:minimality}
Find a necessary and sufficient condition for the minimality of flows. 
\end{problem}

By applying the two concepts of circuits, we establish the following necessary and sufficient condition for minimality, which answers Problem~\ref{prob:minimality} for flows on (possibly non-compact) surfaces. 

\begin{main}\label{main:minimal}
A flow on a connected surface with finite genus and finite ends is minimal if and only if it is non-wandering {\rm(i.e.} any point is non-wandering{\rm)} and the set difference $\overline{O} - O$ for any orbit $O$ contains non-singular points. 
\end{main}

Notice that both conditions (i.e. non-wandering property and the existence of non-singular points in the set difference) in Theorem~\ref{main:minimal} are necessary. 
The Denjoy flow and periodic flows demonstrate this.
Indeed, a Denjoy flow is not minimal, even though the set difference $\overline{O} - O$ contains non-singular points for any orbit $O$. 
Moreover, any periodic flow is non-wandering, but it fails the second condition since $\overline{O} - O = \emptyset$.

\subsubsection{Partial refinement of Poincar{\'e} recurrence theorem}

Applying Theorem~\ref{main:rec_wandering_pt}, essentially, we have the following characterization of the denseness of recurrence in the non-wandering set,  which is complementary to the Poincar{\'e} recurrence theorem for flows with finitely many singular points on compact surfaces.

\begin{main_cor}\label{main_cor:rec_th}
The following are equivalent for a flow with finitely many singular points on a compact surface: 
\\
{\rm(1)} The set of recurrent points is dense in the non-wandering set {\rm (i.e.} $\overline{\mathop{\mathrm{Cl}(v) \sqcup \mathrm{R}(v)}} = \Omega(v)${\rm)}.
\\
{\rm(2)} There are neither strictly limit circuits nor circuits with wandering holonomy. 
\end{main_cor}

Corollary~\ref{main_cor:rec_th} follows from a more general result for flows with finitely many connected components of the singular point set on surfaces with finite genus and finite ends (see Theorem~\ref{thm:rec_th}). 

The compactness (and finite genus) of the surface in Corollary~\ref{main_cor:rec_th} is necessary (see Lemma~\ref{ex:counter_exam01+} for more details). 
The existence of finitely many singular points in Corollary~\ref{main_cor:rec_th} is also necessary (see Lemma~\ref{lem:counter_exam01} for more details). 
Notice that the non-existence of two kinds of circuits is necessary. 
In fact, the non-existence of strictly limit circuits is also necessary even for spherical flows (see Lemma~\ref{ex:counter_exam_strict}). 
The necessity of the non-existence of circuits with wandering holonomy follows from Lemma~\ref{lem:counter_exam02}.

Applying Theorem~\ref{main:rec_wandering_pt} essentially, we obtain the following characterization of the denseness of closed orbits in the non-wandering set, which is a condition in the definition of Devaney's chaos. 

\begin{main}\label{thm:finite_dense}
The following conditions are equivalent for a flow $v$ with finitely many singular points on a compact surface:
\\
$(1)$ $\overline{\mathop{\mathrm{Cl}}(v)} = \Omega(v)$.
\\
$(2)$ There are neither non-closed recurrent orbits, nor strictly limit circuits, nor circuits with wandering holonomy.
\end{main}

Moreover, applying Corollary~\ref{main_cor:rec_th}, we obtain the following affirmative statement for \cite[Question~1]{boyd2015diffeomorphisms}, which asks: If the closed point set is closed, under what additional conditions on the dynamical system can we conclude that the closed point set coincides with the non-wandering set?
This generalizes the results for a continuous map of a closed interval \cite{xiong1981continuous,nitecki1982maps}. 

\begin{main_cor}\label{main_cor:a}
The following conditions are equivalent for a flow $v$ with finitely many singular points on a compact surface $S$:
\\
$(1)$ $\mathop{\mathrm{Cl}}(v) = \Omega(v)$.
\\
$(2)$ The union $\mathop{\mathrm{Cl}}(v)$ is closed, and there are neither non-closed recurrent orbits, nor non-periodic limit circuits, nor circuits with wandering holonomy.
\\
$(3)$ The union $\mathrm{P}(v)$ is open, and there are neither non-closed recurrent orbits, nor non-periodic limit circuits, nor circuits with wandering holonomy.

In any of the above cases, we have $S = \mathop{\mathrm{Cl}}(v) \sqcup \mathrm{P}(v) = \Omega(v) \sqcup \mathrm{P}(v)$.
\end{main_cor}

Corollary~\ref{main_cor:a} follows from a more general result for flows with finitely many connected components of the singular point set on surfaces with finite genus and finite ends (see Corollary~\ref{thm:main02}). 

The existence of finitely many connected components of the singular point set in Corollary~\ref{main_cor:a} is necessary (see Lemma~\ref{lem:counter_exam01} for details). 
We also remark that the non-existence of the three conditions in (2) is necessary. 
In fact, the non-existence of non-closed recurrent points is necessary due to the existence of irrational rotations on tori. 
The non-existence of non-periodic limit circuits is also necessary even for spherical flows (see Lemma~\ref{ex:counter_exam_strict}). 
The necessity of the non-existence of circuits with wandering holonomy follows from  Lemma~\ref{lem:counter_exam02}. 

\subsubsection{Finiteness of the non-wandering set}
To characterize the existence of finitely many non-wandering orbits, we introduce the concept of a non-wandering orbit with wandering holonomy (see Definition~\ref{def:orbit_wandering_hol} in the next section), which returns near itself exactly once and never comes back again, and obtain the following characterization from Theorem~\ref{main:rec_wandering_pt}.

\begin{main_cor}\label{main:c}
The following conditions are equivalent for a flow $v$ on a compact surface:
\\
{\rm(1)} The non-wandering set $\Omega(v)$ consists of finitely many orbits.
\\
{\rm(2)} Each recurrent orbit is a closed orbit, any strictly limit circuit consists of at most finitely many orbits, and there are at most finitely many closed orbits. 
\\
{\rm(3)} Each recurrent orbit is a closed orbit, any limit circuit consists of at most finitely many orbits,  and there are at most finitely many closed orbits. 
\\
{\rm(4)} Each recurrent orbit is a closed orbit, any limit circuit consists of at most finitely many orbits,  and there are at most finitely many pairwise distinct $\omega$-limit sets of points. 
\\
{\rm(5)} Each recurrent orbit is a closed orbit, any limit circuit consists of at most finitely many orbits,  and there are at most finitely many pairwise distinct $\alpha$-limit sets of points. 

In any of the above cases, the non-wandering set $\Omega(v)$ consists of closed orbits, strictly limit circuits, and non-wandering orbits with wandering holonomy. 
\end{main_cor}

\subsubsection{Contents of this paper}

The present paper consists of ten sections.
In the next section, we introduce fundamental concepts as preliminaries.
In \S~3, we discuss the properties of circuits with wandering holonomy and limit circuits.
In \S~4, we establish the existence of finitely many limit circuits for flows with finitely many singular points on compact surfaces (Theorem~\ref{main:fininite_existence_limit_circuits}). This result implies Theorem~\ref{main:inifinie_limit_sets} and Theorem~\ref{main:ch_finite_limit_sets}. 
Furthermore, we characterize the existence of infinitely many $\omega$-limit sets and examine the dependency between positive and negative limit behaviors.
In \S~5, we provide a classification of non-recurrent points assuming the non-existence of non-closed recurrent orbits (Lemma~\ref{lem:wandering_h}), which serves as a key lemma in this paper. 
Based on this, we prove Theorem~\ref{thm:finite_dense}.
In \S~6, we develop tools to extend the main results obtained in the introduction to a large class of flows, such as those satisfying the no-slip boundary condition. 
In particular, we prove a generalization of Theorem~\ref{thm:finite_dense}. 
In \S~7, we show a generalization of Theorem~\ref{main:rec_wandering_pt} for flows with finitely many connected components of the singular point set on surfaces with finite genus and finite ends. 
In addition, we characterize the denseness of recurrence in the non-wandering set, which refines the Poincar{\'e} recurrence theorem for flows with finitely many singular points on surfaces (Theorem~\ref{thm:rec_th}) implying Corollary~\ref{main_cor:rec_th}.
Furthermore, Theorem~\ref{main:minimal} is demonstrated. 
In \S~8, a characterization of the existence of finitely many non-wandering orbits is described, and the correspondence between the closed point set and the non-wandering set is characterized. 
Moreover, Corollary~\ref{main_cor:a} and Corollary~\ref{main:c} are shown. 
In \S~9, we construct examples of flows to demonstrate the necessity of the conditions in our main results and lemmas.
In the final section, we investigate the denseness of non-wandering flows on compact surfaces by applying our characterizations.

\section{Preliminaries}

In this section, we review fundamental notions from topology and dynamical systems, and we introduce the concepts of a strictly limit circuit, a circuit with wandering holonomy, and a non-wandering orbit with wandering holonomy. 

\subsection{Topological notion}
Denote by $\overline{A}$ the closure of a subset $A$ of a topological space, by $\mathop{\mathrm{int}}A$ the interior of $A$, and by $\partial A := \overline{A} - \mathop{\mathrm{int}}A$ the boundary of $A$, where $B - C$ is used instead of the set difference $B \setminus C$ when $B \subseteq C$.
A boundary component of a subset $A$ is a connected component of the boundary of $A$. 
A subset is {\bf locally dense} if its closure has a nonempty interior. 

\subsubsection{Surfaces}
By a {\bf surface}, we mean a two-dimensional paracompact smooth manifold, that need not be orientable.
A two-dimensional paracompact topological manifold is called a {\bf topological surface}.
By Rad{\'o}'s theorem~\cite{Rado1925smooth}, note that any topological surface is homeomorphic to a surface. 

\subsubsection{Curves and loops}
A curve is a continuous mapping $C: I \to X$ where $I$ is a non-degenerate connected subset of a circle $\mathbb{S}^1$. 
Here, a non-degenerate subset contains at least two points. 
A curve is simple if it is injective and embedded as a continuous mapping.
We also denote by $C$ the image of a curve $C$.
Denote by $\partial C := C(\partial I)$ the boundary of a curve $C$ if $C$ can be extended into a continuous map whose domain is $I \cup \partial I$, where $\partial I$ is the boundary of $I \subset \mathbb{S}^1$. 
Put $\mathop{\mathrm{int}} C := C \setminus \partial C$ if $\partial C$ is defined. 
A simple curve is a simple closed curve if its domain is $\mathbb{S}^1$ (i.e. $I = \mathbb{S}^1$).
A simple closed curve is also called a {\bf loop}. 
An {\bf arc} is a simple curve whose domain is an interval. 

\subsubsection{Intervals in arcs}
For an arc $I$ and any points $a, b \in I$, denote by $(a,b)_I$ the open subarc in $I$ from $a$ to $b$, by $[a,b]_I$ the closed subarc in $I$ from $a$ to $b$, by $(a,b]_I$ the half-open subarc in $I$ open at $a$ and closed at $b$, and by $[a,b)_I$ the half-open subarc in $I$ closed at $a$ and open at $b$.

\subsubsection{Concept ``parallel to the boundary''}
 
\begin{definition}
A subset which is either a loop or an annulus is {\bf parallel to the boundary} of an annulus $A$ if it is contained in $A$ and is homotopically non-trivial {\rm(i.e.} non-contractible {\rm)} in $A$. 
\end{definition}

\subsubsection{Ends and end completion}

We define the end completion of a topological space, which Freudenthal introduced \cite{Freudenthal1931end}, as follows: 

Consider the direct system $\{K_\lambda\}$ of compact subsets of a topological space $X$ and inclusions $K_\lambda \subseteq K_{\lambda +1}$ such that the interiors of $K_\lambda$ cover $X$.  
There is a corresponding inverse system $\{ \pi_0( X - K_\lambda ) \}$, where $\pi_0(Y)$ denotes the set of connected components of a space $Y$. 
Then the {\bf set of ends} of $X$ is defined to be the inverse limit of this inverse system, and any element of the set of ends is called an {\bf end} of $X$. 
Notice that $X$ has one end $x_{\mathcal{U}}$ for each sequence $\mathcal{U} := (U_i)_{i \in \mathbb{Z}_{>0}}$ with $U_i \supseteq U_{i+1}$ such that $U_i$ is a connected component of $X - K_{\lambda_i}$ for some $\lambda_i$. 
Consider the disjoint union $X_{\mathrm{end}}$ of $X$ and the set of ends as a set. 
A subset $V$ of the union $X_{\mathrm{end}}$ is an open neighborhood of an end $x_{\mathcal{U}}$ if there is some $i \in \mathbb{Z}_{>0}$ such that $U_i \subseteq V$, where $\mathcal{U} = (U_i)_{i \in \mathbb{Z}_{>0}}$. 
Then the resulting topological space $\bm{X_{\mathrm{end}}}$ is called {\bf the end completion} (or end compactification) of $X$. 
Note that the end completion is not compact in general.

\subsection{Notions of dynamical systems}
A {\bf flow} is a continuous $\R$-action on a topological space.
Let $v \colon \R \times X \to X$ be a flow on a topological space $X$.
For a point $x$ of $X$, we denote by $O(x)$ (or $O_v(x)$) the orbit of $x$, by $O^+(x)$ the positive orbit (i.e. $O^+(x) := \{ v(t,x) \mid t > 0 \}$), and by $O^-(x)$ the negative orbit (i.e. $O^-(x) := \{ v(t,x) \mid t < 0 \}$).

A point $x$ of $X$ is singular if $x = v(t,x)$ for any $t \in \R$, is non-singular if $x$ is not singular, and is periodic if there is a positive number $T > 0$ such that $x = v(T, x)$ and $x \neq v(t,x)$ for any $t \in (0, T)$.
A point is {\bf closed} if it is singular or periodic.
An orbit is singular (resp. periodic, closed) if it contains a singular (resp. periodic, closed) point. 

Notice that each closed orbit is a closed subset of a surface $S$, and conversely, that if the surface is compact, then each orbit that is a closed subset of $S$ is a closed orbit. 

Denote by $\mathop{\mathrm{Sing}}(v)$ the set of singular points of $v$, called the singular point set, and by $\mathop{\mathrm{Per}}(v)$  the union of periodic orbits of $v$, called the periodic point set.
The union $\bm{\mathop{\mathrm{Cl}}(v)} := \mathop{\mathrm{Sing}}(v) \sqcup \mathop{\mathrm{Per}}(v)$ is the union of closed orbits of $v$, called the closed point set, where $\sqcup$ denotes a disjoint union. 

\subsubsection{Recurrence and non-wandering property}
The $\omega$-limit (resp. $\alpha$-limit)  set of a point $x$ is $\omega(x) := \bigcap_{n\in \mathbb{R}}\overline{\{v(t,x) \mid t > n\}}$ (resp.  $\alpha(x) := \bigcap_{n\in \mathbb{R}}\overline{\{v(t,x) \mid t < n\}}$).
For an orbit $O$, define $\omega(O) := \omega(x)$ and $\alpha(O) := \alpha(x)$ for some point $x \in O$.
Note that an $\omega$-limit (resp. $\alpha$-limit) set of an orbit is independent of the choice of point in the orbit. 
An $\omega$-limit (resp. $\alpha$-limit) set of a point is  {\bf non-degenerate} if the $\omega$-limit (resp. $\alpha$-limit) set is not a singleton. 
A point $x$ of $X$ is {\bf recurrent} if $x \in \omega(x) \cup \alpha(x)$.
A {\bf Q-set} (or quasi-minimal set) is defined to be the closure of the orbit of a non-closed recurrent point.
Denote by $\mathrm{R}(v)$ the set of non-closed recurrent points.
Note that, since we would later like to form a partition of $X$ into orbits with distinct topological behaviors, we defined $\mathrm{R}(v)$ such that it is disjoint from $\mathop{\mathrm{Cl}}(v)$.

A point is {\bf wandering} if there are its neighborhood $U$ and a positive number $N$ such that $v(t,U) \cap U = \emptyset$ for any $t > N$. 
A point is {\bf non-wandering} if it is not wandering. 
Notice that a point is non-wandering if and only if for any of its neighborhoods $U$ and any positive number $N$, there is a number $t \in \mathbb{R}$ with $|t| > N$ such that $v(t,U) \cap U \neq \emptyset$.
Denote by $\bm{\Omega (v)}$ the set of non-wandering points, called the non-wandering set. 
A flow is {\bf non-wandering} if every point is non-wandering.  
A flow on a topological space $X$ is {\bf minimal} if every orbit is dense in $X$.  

\subsubsection{Invariant subsets}
A subset $A$ is positive (resp. negative) invariant if $v(t,A) \subseteq A$ for any $t \in \R_{>0}$ (resp. $t \in \R_{<0}$). 
A subset is {\bf invariant} (or saturated) if it is a union of orbits. 
The saturation of a subset is the union of orbits intersecting it.
Denote by $v(A)$ the saturation of a subset $A$.
Then $v(A) = \bigcup_{a \in A} O(a)$.
A nonempty closed invariant subset is minimal if it contains no proper nonempty closed invariant subsets. 
Notice that the restriction of a flow to a minimal subset is a minimal flow. 

\subsubsection{Separatrices and orbit arcs}
A non-singular orbit is a {\bf separatrix} if there is a singular point which is either its $\alpha$-limit set or its $\omega$-limit set. 
A separatrix is {\bf connecting} if each of its $\alpha$-limit and $\omega$-limit sets is a singular point. 
An {\bf orbit arc} is an arc contained in an orbit. 

\subsubsection{Topological equivalence}
A flow $v$ on a topological space $X$ is {\bf topologically equivalent} to a flow $w$ on a topological space $T$ if there is a homeomorphism $h \colon X \to T$ such that the image of any orbit of $v$ is an orbit of $w$ with preservation of the direction in time. 

To classify and analyze non-invariant subsets (e.g. small neighborhoods of sinks), we define the equivalence of restrictions of flows as follows.

\begin{definition}
The restriction of a flow $v$ to a subset $U$ of a topological space $X$ is {\bf topologically equivalent} to the restriction of a flow $w$ to a subset $V$ of a topological space $Y$ if there is a homeomorphism $h \colon U \to V$ satisfying the following conditions: 
\\
{\rm(1)} The images via $h$ of any orbit arcs of $v$ in $U$ are orbit arcs of $w$. 
\\
{\rm(2)} The inverse images via $h$ of any orbit arcs of $w$ in $V$ are orbit arcs of $v$. 
\\
{\rm(3)} Both $h$ and $h^{-1}$ preserve the directions of the orbit arcs of $v$ and $w$. 
\end{definition}

\subsubsection{Topological properties of orbits}\label{sec:tpo}

From now on, we will only consider flows on surfaces unless otherwise stated. 
%
An orbit is {\bf proper} if it is embedded as a submanifold and {\bf exceptional} if it is neither proper nor locally dense. 
A point is proper (resp. locally dense, exceptional) if its orbit is proper (resp. locally dense, exceptional).
Denote by $\mathrm{LD}(v)$ (resp. $\mathrm{E}(v)$, $\mathrm{P}(v)$) the union of locally dense orbits (resp. exceptional orbits, non-closed proper orbits). 
Note that, since we would later like to form a partition of a surface into orbits with distinct topological behaviors, we defined the invariant subset $\mathrm{P}(v)$ such that it is disjoint from $\mathop{\mathrm{Cl}}(v)$.
Moreover, an orbit on a paracompact manifold (e.g. a surface) is proper if and only if it has a neighborhood in which the orbit is closed \cite{yokoyama2019properness}.  
We have the following observation. 
\begin{lemma}[Lemma~2.1\cite{yokoyama2023topological}]\label{lem:decomp}
The following statements hold for a flow $v$ on a paracompact manifold $S$: 
\\
{\rm(1)} The union $\mathrm{P}(v)$ of non-closed proper orbits is the set of non-recurrent points. 
\\
{\rm(2)} The subset $\mathop{\mathrm{Cl}}(v) \sqcup \mathrm{LD}(v) \sqcup \mathrm{E}(v)$ is the set of recurrent points. 
\\
{\rm(3)} The subset $\mathrm{R}(v) = \mathrm{LD}(v) \sqcup \mathrm{E}(v)$ is the union of non-proper orbits. 
\\
{\rm(4)} $S = \mathop{\mathrm{Cl}}(v) \sqcup \mathrm{P}(v) \sqcup \mathrm{R}(v)$. 
\end{lemma}

\subsubsection{Flow box}
Let $v$ be a flow on a surface $S$ and $v_{\mathbf{e}}$ be the flow on $\R^2$ generated by a unit vector field $\mathbf{e} = (1,0)$. 

\begin{definition}
A subset $D \subset S$ is a {\bf trivial flow box} with respect to $v$ if the restriction $v\vert_D$ of $v$ to $D$ is topologically equivalent to the restriction $v_{\mathbf{e}}\vert_{\mathbb{D}}$ of the flow $v_{\mathbf{e}}$ to $\mathbb{D}$, where $\mathbb{D} := I \times J \subset \R^2$ is a rectangle for some intervals $I, J \in \{ (0,1), (0,1], [0,1), [0,1] \}$. 
\end{definition}

\begin{figure}
\begin{center}
\includegraphics[scale=0.35]{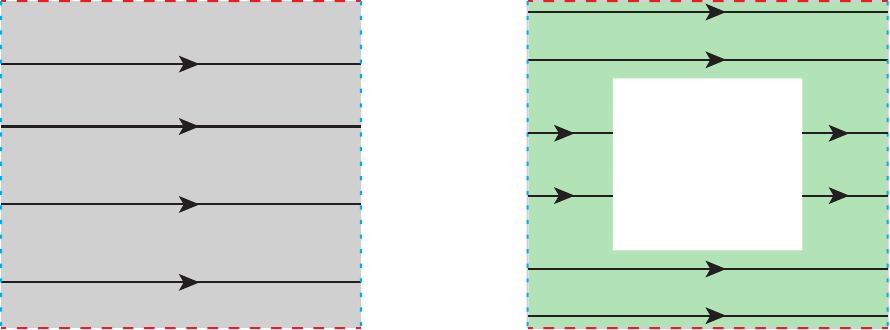}
\end{center}
\caption{Left, a trivial flow box; right, the intersection of a flow box and its small neighborhood of the boundary.}
\label{fig:flowbox}
\end{figure}

\begin{definition}
A subset $B \subset S$ is a {\bf flow box} with respect to $v$ if there is a closed subset $B'$ in the interior of $B$ and there is a homeomorphism $h \colon B - B' \to (I \times J) - [1/4,3/4]^2$ for some intervals $I, J \in \{ (0,1), (0,1], [0,1), [0,1] \}$ such that the restriction $v\vert_{B - B'}$ is topologically equivalent to the restriction  $v_{\mathbf{e}}\vert_{(I \times J) - [1/4,3/4]^2}$ of the flow $v_{\mathbf{e}}$ via $h$, as on the right in Figure~\ref{fig:flowbox}. 
\end{definition}

Note that what is typically referred to as a ``flow box'' in the standard literature corresponds to a ``trivial flow box'' in our terminology.

\subsubsection{One-sided loops and two-sided loops}

We also define the one-sidedness and two-sidedness for loops which is contained in either the boundary of a surface or the interior of the surface as follows. 

\begin{definition}
Let $\gamma$ be a loop in a surface $S$ with either $\gamma \subseteq \partial S$ or $\gamma \cap \partial S = \emptyset$. 
Then the loop $\gamma$ is {\bf one-sided} if either $\gamma \subseteq \partial S$ or it has a small neighborhood which is a M{\"o}bius band, and is {\bf two-sided} if it is not one-sided. 
\end{definition}

Notice that a loop $\gamma$ in the interior of $S$ is two-sided if and only if it has an open small annular neighborhood $U$ such that the complement $U - \gamma$ consists of two open annuli.

\subsubsection{Circuits}
By a {\bf cycle} or a periodic circuit, we mean a periodic orbit.

\begin{definition}\label{def:circuit}
An image of an oriented circle by a continuous orientation-preserving mapping is a {\bf circuit} if it is either a cycle or a directed graph which is the union of (possibly infinitely many) connecting separatrices and finitely many singular points such that the orientation of any edge of the directed graph coincides with the dynamical orientation of the corresponding connecting separatrix.
\end{definition}

A circuit is a {\bf trivial circuit} if it is a singular point.
A circuit is a {\bf non-trivial circuit} if it is not a singular point.

\subsubsection{Collars of non-trivial circuits, one-sided circuits, and two-sided circuits}

We define collars of non-trivial circuits and one-sided circuits as follows.  

\begin{definition}
An open annular subset $\mathbb{A}$ of a surface is a {\bf collar} of a closed invariant subset $\gamma$ if $\gamma$ is a boundary component of $\mathbb{A}$ and there is a neighborhood $U$ of $\gamma$ such that $\mathbb{A}$ is a connected component of the complement $U - \gamma$. 
\end{definition}

Note that not every circuit needs to admit a collar. 

\begin{definition}
A non-trivial circuit $\gamma$ is {\bf one-sided} at an non-singualr orbit $\mu \subset \gamma$ if for any small neighborhood $W$ of $\gamma$ there is a collar $V \subset W$ of $\gamma$ such that the union $V \sqcup \gamma$ is a neighborhood of $\mu$. 
\end{definition}

Notice that a circuit can be one-sided at one separatrix but not at another separatrix.

\begin{definition}
A non-trivial circuit $\gamma$ is {\bf two-sided} if it is not one-sided at any non-singular orbit in $\gamma$.  
\end{definition}

In other words, a non-trivial circuit $\gamma$ is two-sided if and only if there is a small neighborhood $U$ of $\gamma$ such that, for any collar $V \subset U$ of $\gamma$, the union $V \sqcup \gamma$ is not a neighborhood of any point in $\gamma \setminus \Sv$.
Moreover, we have the following observation. 

\begin{lemma}\label{lem:sided}
Let $v$ be a flow on a surface $S$ and $\mu$ a non-trivial circuit in $S$ which is a loop such that either $\mu \cap \partial S = \emptyset$ or $\mu \subseteq \partial S$. 
Then $\mu$ is one-sided {\rm(}resp. two-sided{\rm)} as a circuit if and only if it is one-sided {\rm(}resp.  two-sided  {\rm)} as a loop. 
\end{lemma}

\subsubsection{Limit circuits}

We define limit circuits as follows. 

\begin{definition}\label{def:limit_circuit}
A non-trivial circuit $\gamma$ is a {\bf limit circuit} if there is a point $x \notin \gamma$ with $\omega(x) = \gamma$ or $\alpha(x) = \gamma$. 
\end{definition}

\begin{figure}
\begin{center}
\includegraphics[scale=0.28]{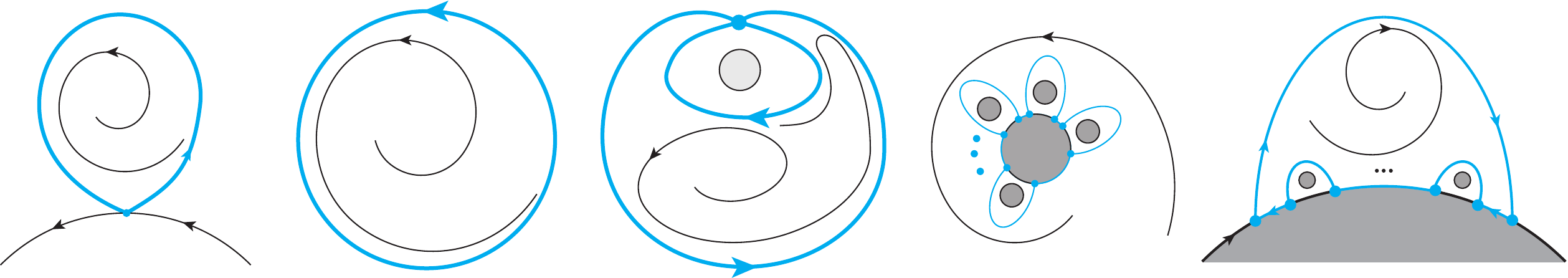}
\end{center}
\caption{Examples of limit circuits.}
\label{fig:limit_circuit}
\end{figure}

By the definition of a non-trivial circuit, notice that any limit circuit consists of closed orbits and non-recurrent orbits (see Figure~\ref{fig:limit_circuit}). 
We define semi-attracting and semi-repelling properties for circuits. 

\begin{definition}\label{def:semi-att}
A circuit $\gamma$ is a {\bf semi-attracting} circuit with respect to a small positive invariant collar $\A$ if $\omega(x)  = \gamma$ for any point $x \in \A$. 
\end{definition}

\begin{definition}\label{def:semi-rep}
A circuit $\gamma$ is a {\bf semi-repelling} circuit with respect to a small negative invariant collar $\A$ if $\alpha(x) = \gamma$ for any point $x \in \A$. 
\end{definition}

Then $\A$ is called a {\bf semi-attracting} (resp. a {\bf semi-repelling}) {\bf collar basin} of $\gamma$. 
%

\subsubsection{Limit cycle}

We define limit cycles as follows.

\begin{definition}\label{def:limit_cycle_attr}
A periodic orbit $\gamma$ is a {\bf semi-attracting limit cycle} with respect to a small positive collar $\A$ if $\gamma$ is a semi-attracting circuit with respect to $\A$. 
\end{definition}

\begin{definition}
A periodic orbit $\gamma$ is a {\bf semi-repelling limit cycle} with respect to a small negative collar $\A$ if $\gamma$ is a semi-repelling circuit with respect to$\A$. 
\end{definition}

\begin{definition}\label{def:limit_cycle}
A periodic orbit $\gamma$ is a {\bf limit cycle} with respect to a small collar $\A$ if $\gamma$ is either semi-attracting or semi-repelling limit cycle with respect to $\A$. 
\end{definition}

By \cite[Theorem~A]{yokoyama2021poincare}, limit cycles coincide with limit circuits in $\mathop{\mathrm{Per}}(v)$ as a set.
Thus, we establish the convention detailed in the following remark.

\begin{remark}\label{rem:counting}
We consider a limit cycle as two distinct semi-attracting {\rm(resp.} semi-repelling{\rm)} limit cycles if it is a semi-attracting {\rm(resp.} semi-repelling{\rm)} limit cycle with respect to two disjoint collars $\A$ and $\A'$, respectively, even if it is identical as a subset.
Thus, denote by $\widetilde{N}_{l,\omega} \in \Z_{\geq 0} \sqcup \{ \infty \}$ {\rm(resp.} $\widetilde{N}_{l,\alpha} \in \Z_{\geq 0} \sqcup \{ \infty \}${\rm)} {\rm(}used in {\rm Lemma~\ref{lem:counter_alpha_limit_cycles})} the number of pairwise distinct $\omega$-limit {\rm(resp.} $\alpha$-limit{\rm)} sets which are limit cycles counted according to this convention.
\end{remark}

\subsubsection{Poincar{\'e}-Bendixson theorem and Ma\u{\i}er-Markley theorem}

We recall a generalization of the Poincar{\'e}-Bendixson theorem for a flow with finitely many singular points, and the Ma\u{\i}er-Markley theorem, which are going to be frequently used in the next pages.

\begin{lemma}[cf. Corollary~5.8\cite{yokoyama2021poincare}]\label{lem:lc_transversal}
Let $v$ be a flow with finitely many singular points on a compact surface $S$.
Then the $\omega$-limit set of any non-closed orbit is one of the following exclusively:
\\
$(1)$ A singular point.
\\
$(2)$ A semi-attracting limit cycle.
\\
$(3)$ A semi-attracting limit non-periodic circuit.
\\
$(4)$ A locally dense Q-set.
\\
$(5)$ A transversely Cantor Q-set.
\end{lemma}

Although transversely Cantor Q-sets will not appear further in this paper, for the sake of completeness, the reader is referred to, for example, \cite{yokoyama2021poincare} for the definition of a transversely Cantor Q-set. 

\begin{lemma}[\cite{mayer1943trajectories,markley1970number} (cf. Remark 2 \cite{aranson1996maier})]\label{lem:Maier}
The following statements hold for a flow $v$ on a compact surface $S$:
\\
{\rm(}1{\rm)} The total number of Q-sets for $v$ is at most $g$ if $S$ is an orientable surface of genus $g$. 
\\
{\rm(}2{\rm)} The total number of Q-sets for $v$ is at most $\frac{p-1}{2}$ if $S$ is a non-orientable surface of non-orientable genus $p$.
\end{lemma}

The previous lemmas imply the following statement. 

\begin{lemma}\label{lem:annulus_nonrec}
Let $v$ be a flow on a compact surface $S$ and $A \subset S$ an open annulus. 
Then the following statements hold:
\begin{enumerate}[\rm (1)]
\item For any point $x \in A$ whose orbit is contained in $A$, we have that $x \notin \mathrm{R}(v)$. 
\item For any non-singular point $y \in S$ with $\omega(y) \subseteq A$, the $\omega$-limit set $\omega(y)$ is either a singular point or a limit circuit.
\end{enumerate}
\end{lemma}

\begin{proof}
By Gutierrez's smoothing theorem~\cite[Smoothing theorem]{gutierrez1986smoothing}, we may assume that $v$ is $C^1$. 
Denote by $X$ the $C^1$ vector field on $S$ which generates the $C^1$ flow $v$.
Using a $C^\infty$ bump function $\varphi \colon S \to [0,1]$ with $S - A = \varphi^{-1}(0)$, we obtain a $C^1$ vector field $\varphi X$ which generates a flow $v'$ with $S - A \subseteq \mathop{\mathrm{Sing}}(v')$. 
Then $O(x) = O_{v'}(x)$ for any $x \in A$ with $O(x) \subset A$. 
By the end completion of the open annulus $A$, which is also the two-point compactification, the open annulus $A$ can be considered as a subset of the resulting sphere $\mathbb{S}^2$. 
Moreover, by considering the two ends as singular points on the sphere, we obtain the resulting flow $v_{\mathbb{S}^2}$ on the sphere $\mathbb{S}^2$ from $v'\vert_A$. 
Then the subset $O(x)$ for any $x \in A$ with $O(x) \subset A$ can be considered as an orbit of $v_{\mathbb{S}^2}$. 
Lemma~\ref{lem:Maier} implies the non-existence of non-closed recurrent points of $v_{\mathbb{S}^2}$. 
Therefore, for any point $x \in A$ with $O(x) \subset A$, the orbit $O_{\mathbb{S}^2}(x) = O(x)$ is not non-closed recurrent with respect to $v_{\mathbb{S}^2}$ and so $v$, which implies $x \notin \mathrm{R}(v)$. 

Fix the $\omega$-limit set $\omega(y) \subseteq A$ of a non-singular point $y \in S$.
Since $\omega(y)$ is invariant, by the assertion (1), any orbit in $\omega(y)$ is not non-closed recurrent (i.e. $\omega(y) \cap \mathrm{R}(v) = \emptyset$).
Since any Q-set intersects $\mathrm{R}(v)$, Lemma~\ref{lem:lc_transversal} implies that the $\omega$-limit set $\omega(y)$ is either a singular point or a limit circuit.
\end{proof}

\subsubsection{Transversality}
We define transversality immediately as follows.  

\begin{definition}
A simple curve $C$ is {\bf transverse} to $v$ at a non-singular point $p \in C \cap \partial S$ if there are a small neighborhood $U$ of $p$ and a homeomorphism $h:U \to [-1,1] \times [0,1]$ with $h(p) = (0,0)$ such that $h^{-1}([-1,1] \times \{t \})$ for any $t \in [0, 1]$ is an orbit arc and $h^{-1}(\{0\} \times [0,1]) = C \cap U$. 
\end{definition}

\begin{definition}
A simple curve $C$ is {\bf transverse} to $v$ at a non-singular point $p \in C \setminus \partial S$ if there are a small neighborhood $U$ of $p$ and a homeomorphism $h:U \to [-1,1]^2$ with $h(p) = (0,0)$ such that $h^{-1}([-1,1] \times \{t \})$ for any $t \in [-1, 1]$ is an orbit arc and $h^{-1}(\{0\} \times [-1,1]) = C \cap U$. 
\end{definition}

A simple curve $C$ is transverse to $v$ if it is transverse to $v$ at every point of $C$.  
A simple curve $C$ which is transverse to $v$ is called a {\bf transverse arc}.
A simple closed curve is a {\bf closed transversal} if it is transverse to $v$.
By definition, a closed transversal $\mu$ is two-sided. 
We have the following observation on the existence of a collar, one of whose boundary components is a closed transversal.

\begin{lemma}\label{lem:existence_closed_transversal}
Let $v$ be a flow on a surface $S$ and $\gamma$ be a semi-attracting circuit with respect to a small positive invariant collar $\A$. 
Then there is a positive invariant collar $\A' \subset \A$ of $\gamma$ such that the boundary component of $\A'$ which is not $\gamma$ is a closed transversal. 
\end{lemma}

\begin{proof}
Since $\gamma$ is a semi-attracting circuit, there are a point $x \in S$ with $O^+(x) \subset \mathop{\mathrm{int}} \A$, a point $y \in \gamma = \omega(x)$, and an open transverse arc $I \subset \mathop{\mathrm{int}} \A$ whose boundary is $\{x,y\}$ such that $\vert I \cap O(x)\vert = \infty$.
From \cite[Lemma~3.2]{yokoyama2021poincare}, there are an orbit arc $C$ in $O^+(x) \subset \mathop{\mathrm{int}} \A$ and a transverse closed arc $J \subseteq I \subset \mathop{\mathrm{int}} \A$ such that the union $\mu := J \cup C$ is a loop in the annulus $\mathop{\mathrm{int}} \A$.
Moreover, for a small number $\varepsilon > 0$ such that the $\varepsilon$-neighborhood $B_{\varepsilon}(\mu)$ of $\mu$ with respect to a Riemannian distance of $S$ is contained in $\mathop{\mathrm{int}} \A$, there is a closed transversal $\nu \subset B_{\varepsilon}(\mu) \subset \mathop{\mathrm{int}} \A$. 
Therefore, we can choose a desired positive invariant collar $\A' \subset \mathop{\mathrm{int}} \A$ of $\gamma$ whose boundary is $\nu \sqcup \gamma$. 
\end{proof}

Note that, in the previous proof, we apply the waterfall construction (cf. \cite[Lemma~3.3.7, p.86]{candel2000foliation}) to $\mu$ in the construction of a closed transversal from $\mu$ using \cite[Lemma~3.2]{yokoyama2021poincare}.






\subsection{Strictly limit circuit and circuit with wandering holonomy}
To characterize the denseness of closed orbits in the non-wandering set, we introduce circuits with wandering holonomy and strictly limit circuits as follows. 

\begin{definition}\label{def:sls}
A non-periodic limit circuit $\gamma$ is a {\bf strictly limit circuit} at a separatrix $\mu \subseteq \gamma$ if either $\gamma$ is one-sided at $\mu$ or there is a transverse open arc $T \subset S - \mathrm{R}(v)$ intersecting $\mu$ such that $f_v|_{T_1 \cap \mathrm{dom}(f_v)}$ is either attracting or repelling, and that $f_v|_{T_2 \cap \mathrm{dom}(f_v)}$ has no periodic points,  where $T_1$ and $T_2$ are the two connected components of $T \setminus \mu$ and $f_v : \mathrm{dom}(f_v) \subseteq T \to T$ is the first return map.
\end{definition}

\begin{figure}
\begin{center}
\includegraphics[scale=0.35]{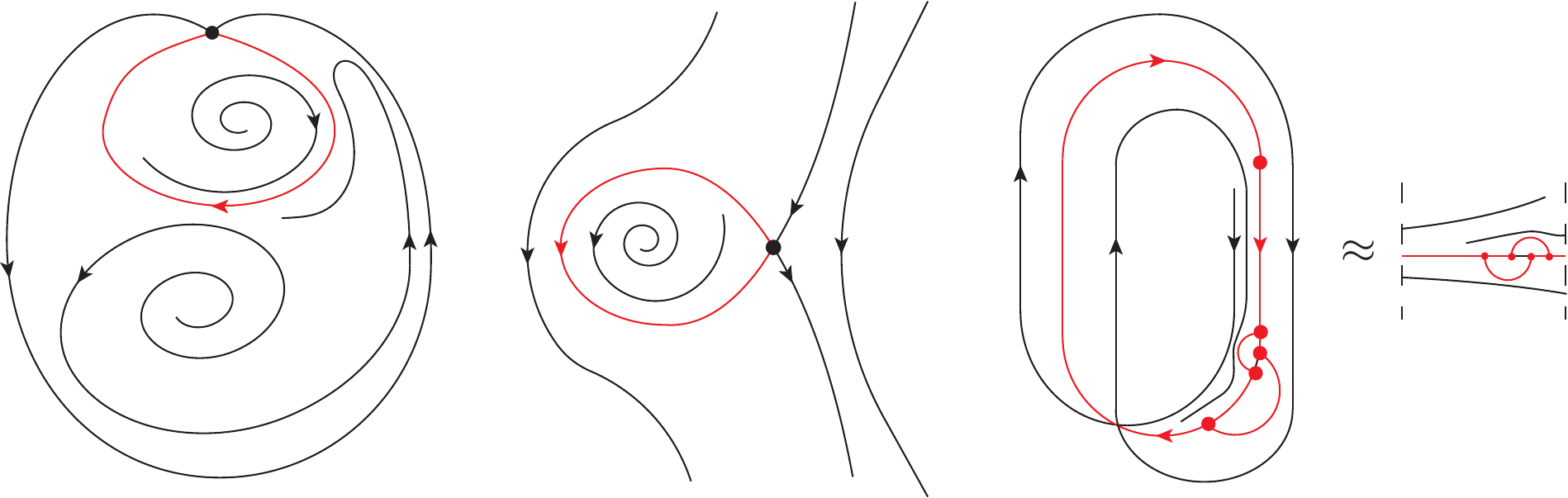}
\end{center}
\caption{Examples of strictly limit circuits.}
\label{fig:strict_limit_circuits}
\end{figure}

Note that a fixed point of a map is treated as a periodic point of the map in this paper.
Moreover, the term ``strictly limit circuit'' is named to emphasize the distinction from an ordinary limit circuit, which may allow the presence of invariant sets such as periodic orbits arbitrarily close to it. 
In contrast, a strictly limit circuit satisfies a more ``strict'' condition: only those orbits that wind around it are permitted to stay near it indefinitely.

\begin{figure}
\begin{center}
\includegraphics[scale=0.2]{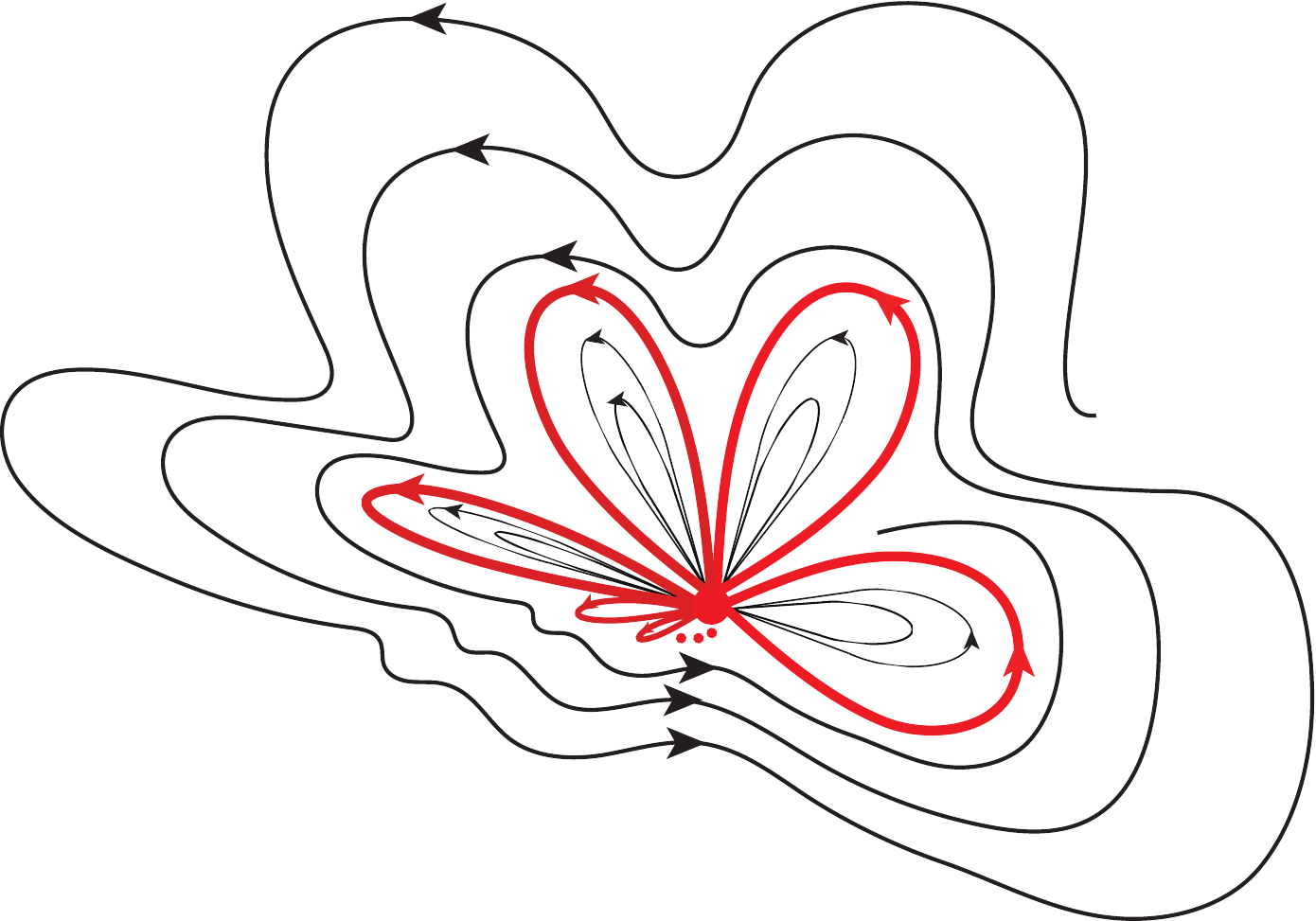}
\end{center}
\caption{A strictly limit circuit that consists of a degenerate singular point and infinitely many connecting separatrices, and its neighborhood which consists of a singular point and non-recurrent orbits.}
\label{NAC02}
\end{figure}

By Lemma~\ref{lem:lc_transversal}, any one-sided strictly limit circuit of a flow with finitely many singular points is either attracting or repelling, as in Figure~\ref{fig:strict_limit_circuits}. 
Note that there is a strictly limit circuit with infinitely many edges and such that any non-trivial circuit in it contains non-closed proper orbits as in Figure~\ref{NAC02}.

\begin{definition}\label{def:cwh}
A non-trivial non-periodic circuit $\gamma$ is a {\bf circuit with wandering holonomy} at a separatrix $O \subset \gamma$ if there are a point $x \in O \subset \gamma$ and an arbitrarily small open transverse arc $I$ containing $x$ such that the following conditions hold:
\\
{\rm(1)} There are two subarcs $I_1, I_2 \subset I$ with $I = I_1 \sqcup \{ x \} \sqcup I_2$. 
\\
{\rm(2)} The first return map $f_{v,I}$ on $I$ is defined from $I_1$ onto $I_2$ (i.e. $\mathrm{dom}(f_{v,I}) = I_1$ and $f_{v,I}(I_1) = I_2$) and reverses the orientation induced on $I_1$ and $I_2$ by $I$.
\\
{\rm(3)} The separatrix $\gamma$ is the connected component of $\partial D - (I - \{ x\})$ containing $x$, where $D$ is the union of open orbit arcs from points in $I_1$ to $I_2$. 
\end{definition}

This definition implies several topological and dynamical constraints on the flow.

\begin{figure}
\begin{center}
\includegraphics[scale=0.2]{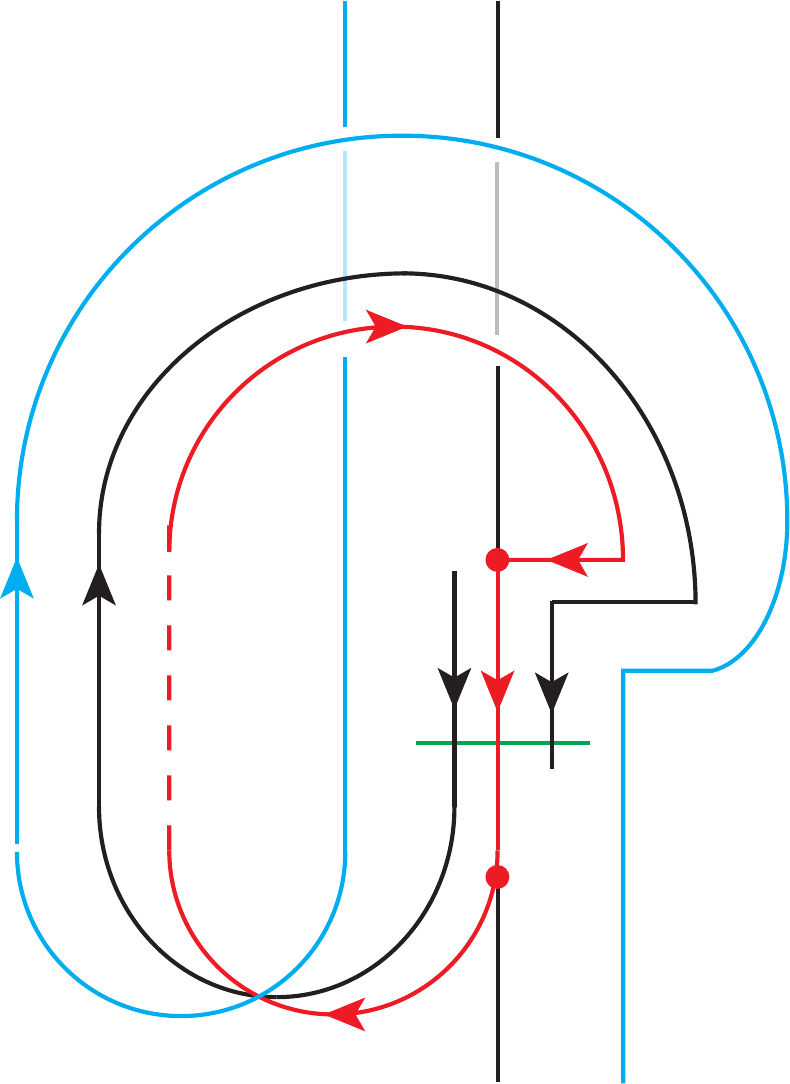}
\caption{An example of a circuit with wandering holonomy.}
\label{nondir}
\end{center}
\end{figure}

\begin{remark}
We make the following observations regarding the terminology and properties of circuits with wandering holonomy.
\begin{enumerate}
\item[{\rm(1)}] The term ``circuit with wandering holonomy'' is named to indicate that any point in a small neighborhood of a separatrix of the circuit that is not in the separatrix is wandering.
\item[{\rm(2)}] A circuit with wandering holonomy cannot exist on an orientable surface, and that, by definition, it cannot be a cycle.
\item[{\rm(3)}] There are flows with circuits with wandering holonomy as in Figure~\ref{nondir} and Figure~\ref{g2-ex}. 
\item[{\rm(4)}] No circuit with wandering holonomy is contained in the $\omega$-limit set or the $\alpha$-limit set of any point, or in the closure of the recurrent point set. 
\item[{\rm(5)}] No strictly limit circuit is contained in the closure of the recurrent point set. 
\end{enumerate}
\end{remark}

Similarly, we define non-wandering orbits with wandering holonomy as follows. 

\begin{definition}\label{def:orbit_wandering_hol}
A non-wandering orbit $O$ is a {\bf non-wandering orbit with wandering holonomy} if there are a non-singular point $x \in O$ and an arbitrarily small open transverse arc $I$ containing $x$ such that the first return map on $I$ is orientation reversing, the domain is nonempty, and the intersection of the domain and the codomain is empty. 
\end{definition}

Notice that any circuit with wandering holonomy contains non-wandering orbits with wandering holonomy. 

\subsection{Notation of a set}
In this paper, denote by $\mathbb{N}$ the set of positive integers (i.e. $\mathbb{N} = \Z_{> 0}$).

\section{Properties of circuits}

In this section, we discuss the properties of circuits with wandering holonomy and the properties of limit circuits, which are going to be utilized later in the proofs of Theorem~\ref{main:inifinie_limit_sets} and Theorem~\ref{thm:finite_dense} below. 

The following lemma ensures the existence of non-recurrent non-wandering points in a circuit with wandering holonomy. 

\begin{lemma}\label{lem:wandering_hol}
Let $v$ be a flow on a surface $S$. 
Then, each circuit with wandering holonomy is not periodic and contains points in $\Omega(v) - \overline{\Cv \sqcup \mathrm{R}(v)}$. 
\end{lemma}

\begin{proof}
Let $\gamma$ be a circuit with wandering holonomy at a separatrix $O$.  
Fix a non-singular point $x \in O$ and a small open transverse arc $I$ containing $x$ as in the definition of circuit with wandering holonomy. 
By the definition of a transverse arc, we have that $I \cap \Sv = \emptyset$. 
Denote by $f_{v,I}$ the first return map on the open transverse arc $I$. 

\begin{claim}\label{claim:001}
$I \cap \Pv = \emptyset$. 
\end{claim}
\begin{proof}[Proof of Claim~\ref{claim:001}]
Assume that $I \cap \Pv \neq \emptyset$. 
Then $\emptyset \neq I \cap \Pv \subseteq \mathrm{dom}(f_{v,I}) \cap \mathrm{codom}(f_{v,I}) = \emptyset$, which is a contradiction. 
\end{proof}

\begin{claim}\label{claim:002}
$I \cap \mathrm{R}(v) = \emptyset$. 
\end{claim}
\begin{proof}[Proof of Claim~\ref{claim:002}]
Assume that $I \cap \mathrm{R}(v) \neq \emptyset$. 
Fix a point $z \in I \cap \mathrm{R}(v)$. 
By time reversion if necessary, we may assume that $z \in \omega(z)$. 
The recurrence of $z$ implies that $z \in \mathrm{dom}(f_{v,I})$ and that there is a point $y \in  \mathrm{codom}(f_{v,I}) \cap O^+(z)$. 
Since the closed orbit arc between $z$ and $y$ is compact, the flow box theorem (cf. \cite[Theorem~4.2.6, p.95]{markley2023flows}, \cite[Theorem 1.1, p.45]{aranson1996introduction}) implies that there is an open interval $J \subseteq \mathrm{dom}(f_{v,I})$ containing $z$. 
The positive recurrence implies that $\emptyset \neq J \cap O^+(z) \subseteq J \cap \mathrm{codom}(f_{v,I}) \subseteq \mathrm{dom}(f_{v,I}) \cap \mathrm{codom}(f_{v,I}) = \emptyset$, which is a contradiction. 
\end{proof}

By Claims~\ref{claim:001} and \ref{claim:002}, the fact that $I \cap \Sv = \emptyset$, and the invariance of $\Cv \sqcup \mathrm{R}(v)$, we have that the saturation $v(I)$ is an open \nbd of $x$ such that $v(I) \cap (\Cv \sqcup \mathrm{R}(v)) = \emptyset$.
Therefore, we have that $x \in v(I) = v(I) \setminus \overline{\Cv \sqcup \mathrm{R}(v)}$. 
Thus $\gamma$ is not periodic, and there is a \nbd $V \subseteq v(I)$ of $x \in \gamma$ containing no closed points. 

\begin{claim}\label{claim:003}
The point $x$ is non-wandering. 
\end{claim}
\begin{proof}[Proof of Claim~\ref{claim:003}]
Fix any \nbd $U \subseteq V$ of $x$ and any $T > 0$.
Since $\gamma$ is a non-periodic non-trivial circuit, the circuit $\gamma$ is a directed graph that is the union of connecting separatrices and finitely many singular points.
Therefore, the $\omega$-limit set $\omega(x)$ is a singular point. 
The continuity of the flow $v$ implies that there is an open transverse arc $J \subset I \cap U$ containing $x$ such that the first return time of any point in $J$ with respect to $v$ is more than $T$ if the first return time exists.  
By the definition of a circuit with wandering holonomy, we have that $J \cap f_{v,I}(J \cap \mathrm{dom}(f_{v,I})) \neq \emptyset$. 
This implies that there is a positive number $t > T$ such that $\emptyset \neq  J \cap v(t,J) \subseteq U \cap v(t,U)$. 
Therefore, the point $x$ is non-wandering. 
\end{proof}

The previous claim implies that $x \in \Omega(v) - \overline{\Cv \sqcup \mathrm{R}(v)}$.
\end{proof}

\subsection{A sufficient condition to be a strictly limit circuit}

The following statement holds, which is used in the proof of Theorem~\ref{thm:finite_dense}. 

\begin{lemma}\label{lem:lim_cycle_strict}
Let $v$ be a flow on a compact surface $S$. 
Then every limit circuit containing a point $x \in S - \overline{\Cv \sqcup \mathrm{R}(v)}$ is a strictly limit circuit at $O(x) \subset \Omega(v)$. 
\end{lemma}

\begin{proof}
Let $\gamma$ be a limit circuit containing a point $x \in S - \overline{\Cv \sqcup \mathrm{R}(v)}$. 
By time reversion if necessary, we may assume that $\gamma$ is the $\omega$-limit set of a point. 
Then $x \in \Omega(v) \cap (\gamma \setminus \overline{\Cv \sqcup \mathrm{R}(v)}) = (\Omega(v) - \overline{\Cv \sqcup \mathrm{R}(v)}) \cap \gamma$. 
If $\gamma$ is one-sided at $O(x)$, then it is a strictly limit circuit at $O(x)$. 
Thus, we may assume that $\gamma$ is not one-sided at $O(x)$. 
Since $\gamma$ is an $\omega$-limit set, by Lemma~\ref{lem:lc_transversal}, there is a small semi-attracting collar basin $\A$ of the limit circuit $\gamma$.
Since $x \in \gamma \setminus \overline{\Cv \sqcup \mathrm{R}(v)}$, there is a small open transverse arc $I \subset S - \overline{\Cv \sqcup \mathrm{R}(v)}$ containing $x$. 
Denote by $I_-$ the connected component of $I - \{ x \}$ which does not intersect $\A$. 
By $I \cap \overline{\Pv} = \emptyset$, the first return map on $I_-$ has no periodic points, and so the circuit $\gamma$ is a strictly limit circuit at $O(x)$. 
\end{proof}

\subsection{Property of three parallel limit circuits} 

We have the following statements. 

\begin{lemma}\label{lem:union_of_annuli}
Let $(A_{1,n+1})_{n \in \N}$ be a sequence of open annuli on a surface $S$ of finite genus  such that $A_{1,n+1} \subseteq A_{1,n+2}$ and that $A_{1,n+1}$ is parallel to the boundary of $A_{1,n+2}$ for any $n \in \mathbb{N}$. 
Then the union $A := \bigcup_{n=1}^\infty A_{1,n+1}$ is an open annulus. 
\end{lemma}

\begin{proof}
By the definition of $A$, the open subset $A$ of $S$ is an orientable open surface. 
Since $S$ is of finite genus, so is $A$. 
From \cite[Theorem~3]{richards1963classification}, any orientable connected surface with  finite genus is homeomorphic to the resulting surface from an orientable closed surface by removing a closed totally disconnected subset. 
This implies that $A$ is an annulus if and only if both any loops which are homotopically non-trivial in $A$ are freely homotopic in $A$ and there is a loop which is homotopically non-trivial in $A$.

\begin{claim}\label{claim:ess_loop}
There is a loop which is homotopically non-trivial in $A$.    
\end{claim}

\begin{proof}[Proof of Claim~\ref{claim:ess_loop}]
Fix a loop $\mu \subset A_{1,2}$ which is parallel to the boundary of $A_{1,2}$.
By the definition of $A_{1,n+1}$, the loop $\mu$ is parallel to the boundary of $A_{1,n+1}$ for any $n \in \N$.
Assume that $\mu$ is homotopically trivial in $A$. 
Then there is a closed disk $D_\mu \subset A$. 
Since $D_\mu$ is compact and the family $\{ A_{1,n+1}\}_{n \in \N}$ is an open cover of $D_\mu$, there is a number $N \in \N$ such that $D \subset A_{1,n+1}$. 
This means that $\mu = \partial D$ is not parallel to the boundary of $A_{1,n+1}$, which contradicts that $\mu = \partial D$ is parallel to the boundary of $A_{1,n+1}$. 
Thus, the assertion holds. 
\end{proof}

\begin{claim}\label{claim:ess_loop02}
Any two loops which are homotopically non-trivial in $A$ are freely homotopic to each other in $A$.
\end{claim}

\begin{proof}[Proof of Claim~\ref{claim:ess_loop02}]
Fix any two loops $\gamma, \gamma'$ which are homotopically non-trivial in $A$. 
Since the loops are compact and the family $\{ A_{1,n+1}\}_{n \in \N}$ is an open cover of $\gamma \cup \gamma'$, there is an integer $N \in \mathbb{N}$ such that the open annulus $\bigcup_{n=1}^N A_{1,n+1} = A_{1,N+1}$ contains $\gamma \cup \gamma'$. 
Because the loops $\gamma$ and $\gamma'$ are homotopically non-trivial in $A$, so are they in the annulus $A_{1,n+1}$. 
This means that the loops $\gamma$ and $\gamma'$  are freely homotopic in $A_{1,n+1}$ and so in $A$.
\end{proof}
These claims imply the assertion.
\end{proof}

\begin{lemma}\label{lem:three_parallel_loops}
Let $v$ be a flow on a surface $S$, $\gamma_1, \gamma_2, \gamma_3$ closed transversals which are boundary components of semi-attracting collar basins $A_1, A_2, A_3$ of pairwise distinct semi-attracting limit circuits $\omega_1, \omega_2, \omega_3$ respectively, and $\mathbb{A}$ an open annulus whose boundary components are $\gamma_1$ and $\gamma_3$ and which is a \nbd of $\gamma_2$. 
Then $\omega_2 \subset \A$. 
Moreover, if the limit circuit $\omega_2$ is non-periodic, then $\mathbb{A}$ contains singular points. 
\end{lemma}

\begin{proof}
By the definition of a collar basin, the pairwise distinct property of $\omega_1, \omega_2, \omega_3$ implies that semi-attracting collar basins $A_1, A_2, A_3$ are pairwise disjoint and that $A_2 \cap \overline{A_1 \sqcup A_3} = \emptyset$. 
From $\partial \mathbb{A} = \gamma_1 \sqcup \gamma_3 \subset \overline{A_1 \sqcup A_3}$, we have $A_2 \cap \partial \mathbb{A} = \emptyset$. 
From $\gamma_2 \subset \mathbb{A} \cap \overline{A_2}$, the openness of $\mathbb{A}$ implies that $\mathbb{A} \cap A_2 \neq \emptyset$. 
By $A_2 \cap \partial \mathbb{A} = \emptyset$, the connectivity of the open annuli $\mathbb{A}$ and $A_2$ implies that $A_2 \subset \mathbb{A}$. 
From $\gamma_2 \subset \mathbb{A}$ and $\overline{A_2} = A_2 \sqcup \partial A_2$, we have $\emptyset = \omega_2 \cap (\gamma_1 \sqcup \gamma_3) =  \omega_2 \cap \partial \mathbb{A} = (\partial A_2 - \gamma_2) \cap \partial \mathbb{A}  = \partial A_2 \cap \partial \mathbb{A}  = \overline{A_2} \cap \partial \mathbb{A}$. 
Therefore, we obtain $\omega_2 \subset \overline{A_2} \subset \mathbb{A}$. 

Suppose that the limit circuit $\omega_2$ is non-periodic. 
Since the non-periodic limit circuit $\omega_2$ contains singular points, the annulus $\mathbb{A}$ contains singular points. 
\end{proof}

\subsection{Properties of infinitely many pairwise disjoint loops}

We state the properties of more general loops, including cycles, as follows.

\begin{lemma}\label{claim:003a}
Let $v$ be a flow with finitely many singular points on a compact surface $S$ and $(\mu_n)_{n \in \N}$ a sequence of pairwise disjoint loops. 
Suppose that the loop $\mu_n$ is either a periodic orbit or a closed transversal for any $n \in \N$. 
There is a strictly increasing sequence $(n_k)_{k \in \N}$ of non-negative integers and there is a sequence $(A_{n_k})_{k \in \N}$ of pairwise disjoint open annuli $A_{n_k}$ each of whose boundaries is $\mu_{n_k} \sqcup \mu_{n_{k+1}}$ satisfying the following properties:
\begin{enumerate}[\rm (1)]
\item The union $A := \bigcup_{k=1}^\infty (A_{n_k} \sqcup \mu_{n_{k+1}})$ is an open annulus one of whose boundary components is the loop $\mu_{n_1}$. 
\item For any $k \in \N$, the annulus $A_{n_k}$ and the loop $\mu_{n_{k+1}}$ are parallel to the boundary of $A$.
\item The set of connected components of $S - \bigcup_{k=1}^\infty \mu_{n_k}$ contained in $A$ coincides with the family $\{ A_{k}\}_{k \in \N}$ of annuli. 
\item For any $k<k' \in \N$, the union $A_{n_k} \sqcup \bigsqcup_{l = k+1}^{k'-1} (\mu_{n_l} \sqcup A_{n_l}) \subset A$ is an open annulus whose boundary is $\mu_{n_k} \sqcup \mu_{n_{k'}}$.
\end{enumerate}
\end{lemma}

\begin{proof}
From the disjoint property of $(\mu_n)_{n \in \N}$ and the finite genus of $S$, by taking a subsequence of $(\mu_n)_{n \in \N}$, we may assume that the loops $(\mu_n)_{n \in \N}$ are pairwise freely homotopic in $S$. 
By the definition of a one-sided loop, any one-sided loop in $(\mu_n)_{n \in \N}$ either is contained in the boundary $\partial S$ or has a small neighborhood which is a M{\"o}bius band.
Since 
$S$ has finite genus and finitely many boundary components, by taking a subsequence of $(\mu_n)_{n \in \N}$, we may assume that the loops $(\mu_n)_{n \in \N}$ are two-sided. 

Suppose that there are infinitely many loops $\mu_n$ each of which is homotopically non-trivial in $S$. 
By taking a subsequence of $(\mu_n)_{n \in \N}$, we may assume that any loop $\mu_n$ is homotopically non-trivial in $S$. 

\begin{claim}\label{claim:ex_annulus}
There are a strictly increasing sequence $(n_k)_{k \in \N}$ of non-negative integers and a sequence $(A_{n_1,n_k})_{k \in \N}$ of open annuli each of whose boundaries is $\mu_{n_1} \sqcup \mu_{n_k}$ satisfying the following properties:
\begin{enumerate}[\rm (1)]
\item $A_{n_1,n_k} \subsetneq A_{n_1,n_{k+1}}$ for any $k \in \N$. 
\item The union $A := \bigcup_{k=1}^\infty A_{n_1,n_k} $ is an open annulus one of whose boundary components is the loop $\mu_1$. 
\item For any $k \in \N_{>1}$, the loop $\mu_{n_k}$  is parallel to the boundary of $A$.
\end{enumerate}
\end{claim}

\begin{proof}[Proof of Claim~\ref{claim:ex_annulus}]
Since $\mu_1$ is two-sided, we refer to the transverse orientations of $\mu_1$ as the right side and the left side. 
Therefore, any collar of $\mu_1$ whose boundary consists of two distinct loops lies on either the right side or the left side exclusively. 
For any $n \in \mathbb{N}$, since $\mu_1$ and $\mu_{n+1}$ are freely homotopic in $S$, there is an open annulus $A_{1,n+1}$ whose boundary is $\mu_1 \sqcup \mu_{n+1}$. 
By taking a subsequence of $(\mu_n)_{n \in \N}$ and reversing the transverse orientations if necessary, we may assume that every $A_{1,n+1}$ lies on the left side. 
Therefore, for any $n \neq m \in \mathbb{N}$, either $A_{1,n+1} \subsetneq A_{1,m+1}$ or $A_{1,m+1} \subsetneq A_{1,n+1}$. 

Suppose that there is a natural number $N \in \mathbb{N}_{>1}$ such that $A_{1,N}$ contains infinitely many loops $\mu_{n_k}$. 
Since the loop $\mu_n$ for any $n \in \N$ is homotopically non-trivial, each $\mu_{n_k}$ is parallel to the boundary of the annulus $A_{1,N}$. 
Then we can define a total order $\sqsubset$ on $\{ \mu_{n_k} \}_{k \in \N}$ by $\mu_{n_{k}} \sqsubset \mu_{n_{k'}}$ if $A_{1,n_{k}} \subsetneq A_{1,n_{k'}}$. 
Since the closure $\overline{A_{1,N}}$ is a closed annulus and so is compact, there is a point $x \in \bigcap_{k=1}^\infty \overline{\bigsqcup_{l=k}^\infty \mu_{n_{l}}} \subset \overline{A_{1,N}}$. 
This implies that there is a strictly increasing sequence $(k_l)_{l \in \N}$ of positive integers such that $(\mu_{n_{k_l}})_{l \in \N}$ is monotonic with respect to the total order $\sqsubset$.
By taking a subsequence of $(n_k)_{k \in \N}$, we may assume that $(\mu_{n_{k}})_{k \in \N}$ is monotonic with respect to the total order $\sqsubset$.
Denote by $A_{n_1,n_k}$ the open annulus in $A_{1,N}$ whose boundary is $\mu_{n_1} \sqcup \mu_{n_{k+1}}$ for any $k \in \N$. 
Put $A := \bigcup_{k=1}^\infty A_{n_1,n_{k+1}}$. 
By Lemma~\ref{lem:union_of_annuli}, the open subset $A$ is an annulus.
Since the loop $\mu_n$ for any $n \in \N$ is homotopically non-trivial in $S$, the loop $\mu_{n_{k+1}} \subset A$ for any $k \in \mathbb{N}$ is parallel to the boundary of $A$. 
Therefore, the sequence $(A_{n_1,n_{k+1}})_{k \in \N}$ has the desired properties, which implies the claim.

Thus, we may assume that the annulus $A_{1,n}$ for any $n \in \mathbb{N}$ does not contain infinitely many loops $\mu_m$. 
By taking a subsequence of $(\mu_n)_{n \in \N}$, we may assume that $A_{1,n+1} \subsetneq A_{1,n+2}$ for any $n \in \mathbb{N}$.
Put $A := \bigcup_{n=1}^\infty A_{1,n+1}$. 
By Lemma~\ref{lem:union_of_annuli}, the open subset $A$ is an annulus.
Since the loop $\mu_n$ for any $n \in \N$ is homotopically non-trivial in $S$, the loop $\mu_{k+1} \subset A$ for any $k \in \mathbb{N}$ is parallel to the boundary of $A$. 
This means that $(A_{1,k+1})_{k \in \N}$ has the desired properties, which implies the claim.
\end{proof}

From Claim~\ref{claim:ex_annulus}, take a strictly increasing sequence $(n_k)_{k \in \N}$ as in Claim~\ref{claim:ex_annulus}. 
By taking a subsequence of $(\mu_n)_{n \in \N}$, we may assume that the sequence $(\mu_{n_k})_{k \in \N}$ is the sequence $(\mu_n)_{n \in \N}$. 
Since $A_{1,n} \sqcup \mu_{n+1} \subsetneq A_{1,n+1}$ for any $n \in \N_{>1}$, put $A_1 := A_{1,2}$ and $A_{n} := A_{1,n+1} - (A_{1,n} \sqcup \mu_{n+1})$ for any $n \in \mathbb{N}_{> 1}$. 
By construction, the sequence $(A_n)_{n \in \mathbb{N}}$ consists of pairwise disjoint open annuli $A_n$ with $\partial A_n = \mu_n \sqcup \mu_{n+1}$.  
Moreover, since $\mu_1$ is a boundary component of $A$, the set of connected components of $S - \bigcup_{n=1}^\infty \mu_{n}$ contained in $A$ coincides with the set of connected components of $A - \bigcup_{n=1}^\infty \mu_{n}$, which is $\{ A_{n}\}_{n \in \N}$. 
This completes the proof of the lemma in the case where there are infinitely many loops $\mu_n$ each of which is homotopically non-trivial.

Thus, we may assume that there are at most finitely many loops $\mu_n$ each of which is homotopically non-trivial.
By taking a subsequence of $(\mu_n)_{n \in \N}$, we may assume that any loop $\mu_n$ is homotopically trivial. 
For any $n \in \N$, denote by $D_n$ a closed disk whose boundary is $\mu_n$. 
By the Poincar{\'e}-Hopf theorem, any closed disk whose boundary is either a periodic orbit or a closed transversal contains singular points. 
By taking a subsequence of $(\omega_n)_{n \in \N}$, from the finiteness of singular points, we may assume that there is a finite subset $F \subseteq \Sv$ such that $F = \Sv \cap D_n = \Sv \cap D_m$ for any $n \neq m \in \N$.
For any $n \neq m \in \N$, since $\partial D_n \cap \partial D_m = \mu_n \cap \mu_m = \emptyset$, we have either $F \subset D_n \subset \mathop{\mathrm{int}}D_m$ or $F \subset D_m \subset \mathop{\mathrm{int}}D_n$. 
By taking a subsequence of $(\mu_n)_{n \in \N}$, we may assume that the sequence $(\mu_n)_{n \in \N}$ satisfies either $D_n \subset \mathop{\mathrm{int}}D_{n+1}$ for any $n \in \N$ or $D_{n+1} \subset \mathop{\mathrm{int}}D_{n}$ for any $n \in \N$. 
For any $n \in \N$, put $A_{n} := D_{n+1} - (D_{n} \sqcup \mu_{n+1})$ if the former case, and put $A_{n} := D_{n} - (D_{n+1} \sqcup \mu_{n})$ if the latter case.
Then $\partial A_n = \mu_n \sqcup \mu_{n+1}$ for any $n \in \N$.
By the construction of $(A_n)_{n \in \N}$, to show the assertion, it suffices to show that $A := \bigsqcup_{n=1}^\infty (A_{n} \sqcup \mu_{n+1})$ is an open annulus. 
Indeed, put $A_{1,2} := A_1$.
For any $n \in \mathbb{N}_{>1}$, the subset $A_{1,n+1} := A_1 \sqcup \bigsqcup_{i=2}^{n} \mu_i \sqcup A_i$ is an open annulus whose boundary is $\mu_1 \sqcup \mu_{n+1}$. 
By construction, we have that $A_{1,n+1} \subset A_{1,n+2}$ for any $n \in \mathbb{N}$. 
Applying Lemma~\ref{lem:union_of_annuli} to $(A_{1,n+1})_{n \in \N}$, the open subset $A$ is an annulus.
\end{proof}

\section{Finiteness of limit circuits for flows with finitely many singular points on compact surfaces}

In this section, we describe the existence of finitely many limit circuits for flows with finitely many singular points on compact surfaces (Theorem~\ref{main:fininite_existence_limit_circuits}), which implies Theorems~\ref{main:inifinie_limit_sets}, \ref{main:ch_finite_limit_sets}, and \ref{thm:finite_dense}, and Corollary~\ref{main:c}. 

\begin{theorem}\label{main:fininite_existence_limit_circuits}
Every flow with finitely many singular points on a compact surface contains at most finitely many pairwise distinct $\omega$-limit sets which are non-periodic limit circuits.
\end{theorem}

\begin{proof}
Assume there is a sequence $( \omega_n )_{n \in \N}$ of infinitely many pairwise distinct $\omega$-limit sets of points which are non-periodic limit circuits.
Fix a sequence $( A_n )_{n \in \N}$ of semi-attracting collar basins of $\omega_n$ respectively. 
From the definition of semi-attracting collar basin (i.e. Definition~\ref{def:semi-att}), we have $A_n \cap A_m = \emptyset$ for any $n \neq m \in \N$. 
By Lemma~\ref{lem:existence_closed_transversal}, by taking $A_n$ small if necessary, we may assume that the boundary $\partial A_n$ for any $n \in \N$ is a disjoint union of $\omega_n$ and a closed transversal.  
Denote by $\mu_n$ the closed transversal with $\partial A_n = \mu_n \sqcup \omega_n$ for any $n \in \N$. 
Then $(A_n \sqcup \mu_n) \cap (A_m \sqcup \mu_m) = \emptyset$ for any $n \neq m \in \N$. 
In particular, the sequence $( \mu_n )_{n \in \N}$ is pairwise disjoint. 
By the finiteness of the genus, by taking a subsequence of the sequence $(\omega_n)_{n \in \N}$, we may assume that closed transversals $\mu_n$ and $\mu_m$ are freely homotopic to each other in the surface for any $n \neq m \in \N$.

From Lemma~\ref{claim:003a}, by taking a subsequence of the sequence $(\mu_n)_{n \in \N}$, we may assume that, for any $n \in \N$, there is an open annulus $U_n$ whose boundary components are $\mu_n$ and $\mu_{n+1}$ such that $U_n \sqcup \mu_{n+1} \sqcup U_{n+1}$ is an open annulus and that the union $A := \bigsqcup_{i =1}^\infty U_i \sqcup \mu_{i+1}$ is an open annulus.
From the finiteness of singular points, by taking a subsequence of the sequence $(\mu_n)_{n \in \N}$, we may assume that $A$ contains no singular points.
Consider the union $W := U_1 \sqcup \mu_2 \sqcup U_2 \subset A$ which is an open annulus and whose boundary components are closed transversals $\mu_1$ and $\mu_3$ where $\mu_1$ (resp. $\mu_2$, $\mu_3$) is a boundary component of semi-attracting collar basin $A_1$ (resp. $A_2$, $A_3$) of the non-periodic semi-attracting limit circuit $\omega_1$ (resp. $\omega_2$, $\omega_3$). 
Since $\mu_1$, $\mu_2$, and $\mu_3$ are pairwise disjoint, Lemma~\ref{lem:three_parallel_loops} implies that $W \subset A$ contains singular points, which contradicts the non-existence of singular points in $A$. 
\end{proof}

Notice that the non-periodicity in the previous theorem is essential, because there is a flow with finitely many singular points on a compact surface that contains infinitely many pairwise distinct $\omega$-limit sets, which are limit cycles.
In fact, the flow $v_X$ generated by the vector field $X$ on a torus $\R^2/\Z^2$ induced by a vector field 
\[
\widetilde{X}(x,y) =
\begin{cases}
\left(x \sin \dfrac{\pi}{x}, 1 \right) & (x \in (0,1]) \\
(0,1) & (x = 0)
\end{cases}
\]
on $[0,1]^2$ is a non-singular flow which has infinitely many pairwise distinct $\omega$-limit sets which are limit cycles.
In particular, since the $y$-component of $\widetilde{X}$ is always $1$, the flow has no singular points. 
On the other hand, for each integer $n \geq 1$, the subset $\{1/n\} \times \R/\Z$ corresponds to a periodic orbit $O_n$ because the $x$-component vanishes at $x=1/n$. 
Moreover, these orbits are isolated except the periodic orbit $O_1$, and thus they are limit cycles that accumulate to the circle $O_1$.

\subsection{Dependency between $\alpha$-limit sets and $\omega$-limit sets}

We have the following lower bound of the number of pairwise distinct $\alpha$-limit sets which are limit cycles. 

\begin{lemma}\label{lem:counter_alpha_limit_cycles}
Let $v$ be a flow with finitely many singular points on an orientable {\rm(resp.} non-orientable{\rm)} connected compact surface $S$ of orientable {\rm(resp.} non-orient\-able{\rm)} genus $g$. 
Denote by $\widetilde{N}_{l,\omega}$ {\rm(resp.} $\widetilde{N}_{l,\alpha}${\rm)} the number of pairwise distinct $\omega$-limit {\rm(resp.} $\alpha$-limit{\rm)} sets which are limit cycles counted as in {\rm Remark~\ref{rem:counting}}, by $N_{l,\omega} \in \Z_{\geq 0} \sqcup \{ \infty \}$ {\rm(resp.} $N_{l,\alpha} \in \Z_{\geq 0} \sqcup \{ \infty \}${\rm)} the number of pairwise distinct $\omega$-limit {\rm(resp.} $\alpha$-limit{\rm)} sets which are limit cycles counted as a subset, by $N_\partial (S)$ the number of boundary components of $S$, and by $N_s$ the number of singular points.
Then the following statements hold:
\begin{enumerate}[\rm (1)]
\item If $g = N_s = 0$, then $S$ is a closed annulus and $\widetilde{N}_{l,\omega} = \widetilde{N}_{l,\alpha} \geq N_{l,\alpha} \geq N_{l,\omega}/2$. 
\item Suppose $g  + N_s > 0$. 
For any non-negative integer $k \in \Z_{\geq 0}$, if the inequality 
\[
N_{l,\omega} \geq  5N_s/2 + 5g/2 + (N_s - \chi(S) + 2)(N_s - \chi(S))/2 + N_\partial (S)/2 + k 
\]
holds, then $N_{l,\alpha} \geq k+1$.
\end{enumerate}
\end{lemma}

\begin{proof}
First, we show the assertion (1) as follows. 

\begin{claim}\label{claim:genus_zero}
If $g = N_s = 0$, then $S$ is a closed annulus, $\widetilde{N}_{l,\omega} = \widetilde{N}_{l,\alpha}$, and $N_{l,\alpha} \geq N_{l,\omega}/2$. 
\end{claim}

\begin{proof}[Proof of Claim~\ref{claim:genus_zero}]
Suppose that the genus is zero (i.e. $g=0$) and that there are no singular points (i.e. $N_s = 0$). 
By the Poincar{\'e}-Hopf theorem, we have that $\chi(S) = 0$. 
From $g=0$, the surface $S$ is a closed annulus. 
By the Poincar{\'e}-Hopf theorem, any closed disk whose boundary is a periodic orbit contains singular points. 
Therefore, the non-existence of singular points implies that any periodic orbit is parallel to the boundary $\partial S$, because $S$ is a closed annulus.
On the other hand, Lemma~\ref{lem:Maier} implies the non-existence of Q-sets. 
By Lemma~\ref{lem:lc_transversal}, the non-existence of singular points implies that the $\omega$-limit set of any non-closed point is a limit cycle, which is parallel to the boundary $\partial S$. 
Similarly, by the dual statement of Lemma~\ref{lem:lc_transversal}, the non-existence of singular points implies that the $\alpha$-limit set of any non-closed point is a limit cycle, which is parallel to the boundary $\partial S$. 
Therefore, each of the $\omega$-limit set and the $\alpha$-limit set of any non-closed point is a limit cycle parallel to $\partial S$. 
This implies that, for any limit cycle $\omega$ which is an $\omega$-limit set, there is a non-closed point $x \in S$ such that $\omega(x) = \omega$ and that the $\alpha$-limit set $\alpha(x)$ is also a limit cycle parallel to $\partial S$ of the closed annulus $S$. 
By symmetry, for any limit cycle $\alpha$ which is an $\alpha$-limit set, there is a non-closed point $y \in S$ such that $\alpha(y) = \alpha$. 
Since all limit cycles are pairwise parallel to the boundary $\partial S$ of the closed annulus $S$, we have that $\widetilde{N}_{l,\omega} = \widetilde{N}_{l,\alpha}$.
From $2 N_{l,\alpha} \geq \widetilde{N}_{l,\alpha} \geq N_{l,\alpha} $ and $\widetilde{N}_{l,\omega} \geq N_{l,\omega}$, we have that $\widetilde{N}_{l,\alpha} \geq N_{l,\alpha} \geq \widetilde{N}_{l,\alpha}/2 = \widetilde{N}_{l,\omega}/2 \geq N_{l,\omega}/2$. 
\end{proof}

Suppose that $g  + N_s > 0$. 
Let $n_+$ be the number of singular points whose indices are positive, $n_-$ the number of singular points whose indices are negative, and $n_{\mathrm{non}}$ the non-orientable genus. 
Then $N_s \geq n_+ + n_-$. 
If $S$ is orientable, then $n_{\mathrm{non}} = 0$. 
If $S$ is non-orientable, then $n_{\mathrm{non}} = g$. 
Define an integer $N_l$ by the minimal integer more than or equal to $5N_s/2 + 5g/2 + (N_s - \chi(S) + 2)(N_s - \chi(S))/2 + N_\partial (S)/2 + k$. 
Suppose that
\[
N_{l,\omega} \geq \frac{5N_s}{2} + \frac{5g}{2} + \frac{(N_s - \chi(S) + 2)(N_s - \chi(S)) + N_\partial (S)}{2} + k.
\]
By the definitions of $N_{l,\omega} \in \Z_{\geq 0 } \sqcup \{ \infty \}$ and $N_l \in \Z$, the following relation holds:
\begin{equation}\label{eq:bound_Nl}
N_{l,\omega} \geq N_l \geq \frac{5N_s}{2} + \frac{5g}{2} + \frac{(N_s - \chi(S) + 2)(N_s - \chi(S)) + N_\partial (S)}{2} + k.
\end{equation}

\begin{claim}\label{claim:no_01}
$N_l \geq 3/2 + g$. 
\end{claim}

\begin{proof}[Proof of Claim~\ref{claim:no_01}]
Suppose that  $\chi(S) = 2$. 
Then the Poincar{\'e}-Hopf theorem implies $N_s \geq 1$. 
Therefore, by the inequality~\eqref{eq:bound_Nl}, we have the following inequality: 
\[
\begin{split}
N_l - g \geq \dfrac{5N_s}{2} + \dfrac{3g}{2} + \dfrac{N_s(N_s - 2) + 0}{2} + 0  &= \dfrac{3N_s}{2} + \dfrac{3g}{2} + \dfrac{N_s^2}{2}  
\\
&\geq \dfrac{3}{2} + 0 + \dfrac{1}{2} = \dfrac{5}{2}
\end{split}
\]

Suppose that  $\chi(S) = 1$. 
Then the Poincar{\'e}-Hopf theorem implies $N_s \geq 1$. 
Therefore, by the inequality~\eqref{eq:bound_Nl}, we have the following inequality:   
\[
\begin{split}
N_l - g \geq \dfrac{5N_s}{2} + \dfrac{3g}{2} + \dfrac{(N_s + 1)(N_s - 1) + 0}{2} + 0  &= \dfrac{5N_s}{2} + \dfrac{3g}{2} + \dfrac{N_s^2}{2} - \dfrac{1}{2}  
\\
&\geq \dfrac{5}{2} + 0 + \dfrac{1}{2} - \dfrac{1}{2} = \dfrac{5}{2} 
\end{split}
\]

Suppose that  $\chi(S) \leq 0$. 
By the inequality~\eqref{eq:bound_Nl}, we have the following inequality:   
\[
\begin{split}
N_l - g \geq \dfrac{5N_s}{2} + \dfrac{3g}{2}  + \dfrac{(N_s + 2)N_s}{2} 
\end{split}
\]
If $N_s \geq 1$, then the following inequality holds:   
\[
\begin{split}
N_l - g \geq \dfrac{5}{2} + \dfrac{3g}{2}  + \dfrac{1 + 2}{2} \geq 4
\end{split}
\]
Since either $g > 0$ or $N_s >0$, we may assume that $g > 0$. 
By the inequality~\eqref{eq:bound_Nl}, we have the following inequality:   
\[
N_l - g \geq 0 + \dfrac{3}{2} + 0 \geq \dfrac{3}{2}
\]
\end{proof}

From $g \geq 0$, the previous claim implies $N_l \in \Z_{\geq 2}$. 
Fix semi-attracting collar basins $\{ A_n \}_{n=1}^{N_l}$ of semi-attracting limit cycles $\{\omega_n\}_{n=1}^{N_l}$ respectively. 
The semi-attracting properties imply that $A_n \cap A_m = \emptyset$ for any $n \neq m  \in \{1, \ldots , N_l\}$. 
By Lemma~\ref{lem:existence_closed_transversal}, we may assume that the boundary $\partial A_n$ is a disjoint union of $\omega_n$ and a closed transversal, denote by $\mu_n$, for any $n \in \{1, \ldots , N_l\}$. 
Taking $A_n$ small if necessary, we may assume that $(A_n \sqcup \mu_n) \cap (A_m \sqcup \mu_m) = \emptyset$ for any $n \neq m \in \Z_{\geq 0}$. 
In particular, the sequence $( \mu_n )_{n \in \{1, \ldots , N_l\}}$ is pairwise disjoint. 
Let $U_1, U_2, \ldots , U_{N}$ be the connected components of the complement $S - \bigcup_{i=1}^{N_l} \mu_i$ for some $N \in \mathbb{N}$. 
Then the number $N$ of connected components of the complement $S - \bigcup_{i=1}^{N_l} \mu_i$ is at least $N_l - g +1 \geq 5/2$, because of Claim~\ref{claim:no_01}. 
By $N > 1$, any connected component $U_i$ has a nonempty boundary in $S$, and so is neither a sphere nor a projective plane.
Moreover, from $N > 1$, any connected component of $S - \bigcup_{i=1}^{N_l} \mu_i$ which is an open annulus has disjoint boundary components. 

\begin{claim}\label{claim:no_02}
The number $\nu_{+}$ of connected components of the complement $S - \bigcup_{i=1}^{N_l} \mu_i$ which are open disks is at most $n_+$.
\end{claim}

\begin{proof}[Proof of Claim~\ref{claim:no_02}]
By the Poincar{\'e}-Hopf theorem, every closed disk whose boundary is a closed transversal contains at least one singular point whose index is positive. 
Since any boundary component of $U_i$ as a subset of $S$ for any $i \in \{1,2, \ldots ,N \}$ (i.e. any connected component of $\overline{U_i} - U_i$) is a closed transversal, the number $\nu_{+}$ of connected components of the complement $S - \bigcup_{i=1}^{N_l} \mu_i$ which are open disks is at most $n_+$.
\end{proof}

\begin{claim}\label{claim:no_03}
The number $\nu_{-}$ of connected components of the complement $S - \bigcup_{i=1}^{N_l} \mu_i$ whose Euler characteristics are negative is at most $n_-$.
\end{claim}

\begin{proof}[Proof of Claim~\ref{claim:no_03}]
By the Poincar{\'e}-Hopf theorem, every closed surface whose Euler characteristic is negative contains at least one singular point whose index is negative. 
For any $i \in \{1,2, \ldots ,N \}$ such that $\chi(U_i) < 0$, since the double of the compact surface $\overline{U_i}$ is a closed surface whose Euler characteristic is negative, because $\overline{U_i} - U_i$ consists of closed transversals, the connected component $U_i$ contains at least one singular point whose index is negative. 
Therefore, the number $\nu_{-}$ of connected components of the complement $S - \bigcup_{i=1}^{N_l} \mu_i$ whose Euler characteristics are negative is at most $n_-$.
\end{proof}

By renumbering the connected components $\{U_i \}_{i=1}^N$ of $S - \bigcup_{i=1}^{N_l} \mu_i$, we may assume that $U_1, \ldots , U_{\nu_+}$ are open disks, $U_{\nu_+ +1}, \ldots , U_{\nu_+ +\nu_-}$ are surfaces whose Euler characteristics are negative, $U_{\nu_+ + \nu_- +1}, \ldots , U_{\nu_+ + \nu_- + \nu_{\mathrm{non}}}$ are M{\"o}bius bands, and $U_{\nu_+ + \nu_- + \nu_{\mathrm{non}} +1}, \ldots , U_{N}$ are annuli, where $\nu_{\mathrm{non}} \leq n_{\mathrm{non}}$ is a non-negative integer. 
Denote by $n_t$ the number of annuli $U_{\nu_+ + \nu_- + \nu_{\mathrm{non}} +1}, \ldots , U_{N}$. 
Then $n_t = N - (\nu_+ + \nu_- + \nu_{\mathrm{non}}) \geq N - (n_+ + n_- + n_{\mathrm{non}})$ due to Claim~\ref{claim:no_02} and Claim~\ref{claim:no_03}. 
Moreover, the annuli $U_{\nu_+ + \nu_- + \nu_{\mathrm{non}} +1}, \ldots , U_{N}$ consist of at least $n_t - N_\partial (S)$ open annuli and at most $N_\partial (S)$ half-open annuli (i.e. annuli homeomorphic to $(0,1] \times \mathbb{S}^1$). 

\begin{claim}\label{claim:08}
The number $N_\partial (U_i)$ of boundary components of any connected component $U_i$ as a subset of $S$ {\rm(i.e.} number of ends of $U_i${\rm)} is at most $N_s - \chi(S) + 2$. 
\end{claim}
\begin{proof}[Proof of Claim~\ref{claim:08}]
Assume that there is a connected component $U_j$ whose number $N_\partial (U_j)$ of boundary components as a subset of $S$ is more than $N_s - \chi(S) +2$ (i.e. $N_\partial (U_j) > N_s - \chi(S) +2$). 
Since the Euler characteristic of any connected surface is at most 2 minus the number of its boundary components, we have that $\chi(U_j) \leq 2 - N_\partial (U_j) < -(N_s - \chi(S))$. 
By Claim~\ref{claim:no_02} and from $N_s \geq n_+ + n_-$, the Euler characteristic $\chi(S)$ satisfies the following inequality: 
\[
\chi(S) = \sum_{i=1}^N \chi(U_i) = \nu_+ + \sum_{i=\nu_+ + 1}^{\nu_+ + \nu_-} \chi(U_i) \leq n_+ + \chi(U_j)  < n_+ - (N_s - \chi(S)) \leq \chi(S)
\]
This is a contradiction. 
\end{proof}

\begin{claim}\label{claim:09}
$n_- \leq N_s - \chi(S)$. 
\end{claim}
\begin{proof}[Proof of Claim~\ref{claim:09}]
Assume that $n_- > N_s - \chi(S)$. 
By $n_+ + n_- \leq N_s$, from Claim~\ref{claim:no_02}, the Euler characteristic $\chi(S)$ satisfies the following inequality: 
\[
\chi(S) = \sum_{i=1}^N \chi(U_i) \leq \nu_+ + \sum_{i=\nu_+ + 1}^{\nu_+ +\nu_-} \chi(U_i) \leq n_+ - \nu_- \leq n_+ \leq N_s - n_-  <  \chi(S)
\]
This is a contradiction. 
\end{proof}

Denote by $N_{\partial, \mathrm{nonA}}$ the sum of the numbers of ends of the non-annular connected components $U_1, \ldots, U_{\nu_+ + \nu_- + \nu_{\mathrm{non}}}$. 
In other words, the sum $N_{\partial, \mathrm{nonA}}$ is the sum of the numbers of boundary components of the non-annular connected components $U_1, \ldots, U_{\nu_+ + \nu_- + \nu_{\mathrm{non}}}$ as a subset of $S$. 
 By Claims~\ref{claim:no_02}--\ref{claim:09}, from $\nu_{\mathrm{non}} \leq n_{\mathrm{non}}$, the sum $N_{\partial, \mathrm{nonA}}$ satisfies

\begin{equation}\label{eq:bound_N_partial}
N_{\partial, \mathrm{nonA}} \leq n_+ + (N_s - \chi(S) + 2)(N_s - \chi(S)) + n_{\mathrm{non}}
\end{equation}
because the sum of the numbers of boundary components of the disks $U_1, \ldots, U_{\nu_+}$ (resp. M{\"o}bius bands $U_{\nu_+ + \nu_- + \nu_{\mathrm{non}} +1}, \ldots , U_{N}$) as a subset of $S$ is at most $\nu_+ \leq n_+$ (resp. $n_{\mathrm{non}}$) and the sum of the numbers of boundary components of surfaces $U_{\nu_+ +1}, \ldots , U_{\nu_+ +\nu_-}$ whose Euler characteristics are negative is at most the multiplication of the number of such components $\nu_- \leq N_s - \chi(S)$ and the largest number $\max_{i \in \{ \nu_+ +1, \ldots , \nu_+ +\nu_- \}} N_\partial (U_i) \leq N_s - \chi(S) + 2$ of their boundary components. 

Since any open annulus has exactly two ends, the sum of the numbers of ends of open or half-open annuli $U_{\nu_+ + \nu_- + \nu_{\mathrm{non}} +1}, \ldots , U_{N}$ is $2 n_t - N_\partial (S)$. 
Because the integer $N_{\partial, \mathrm{nonA}}$ is the sum of the numbers of ends of the non-annular connected components $U_1, \ldots, U_{\nu_+ + \nu_- + \nu_{\mathrm{non}}}$, there are at least $n_t - N_\partial (S)/2 - N_{\partial, \mathrm{nonA}}/2 = (2 n_t - N_\partial (S) - N_{\partial, \mathrm{nonA}})/2$ closed transversals which are intersections of boundaries of annuli in $\{U_i\}_{i=N- n_t +1}^N$.
By the inequality \eqref{eq:bound_Nl}, the inequality holds:  
\begin{equation}\label{eq:inequality_var}
N_l - \dfrac{5g}{2} - N_s \geq \dfrac{3N_s}{2} + \dfrac{(N_s - \chi(S) + 2)(N_s - \chi(S))+ N_\partial (S)}{2} + k 
\end{equation}
We use the following inequalities:
\[
N \geq N_l - g +1, \quad n_{\mathrm{non}} \geq \nu_{\mathrm{non}}, \quad N_s \geq n_+ + n_-, \quad g \geq n_{\mathrm{non}}.
\]
Combining these with the inequalities \eqref{eq:bound_N_partial} and \eqref{eq:inequality_var}, we have that
\begin{align*}
& \hspace{14pt} n_t  
\\
&= N - (\nu_+ + \nu_- + \nu_{\mathrm{non}}) & \hspace{-190pt} (\text{by the definition of } n_t) 
\\
 &\geq N_l - g + 1 - (\nu_+ + \nu_- + n_{\mathrm{non}}) & \hspace{-290pt}  \text{(by $N \geq N_l - g+1$ and $n_{\mathrm{non}} \geq \nu_{\mathrm{non}}$)} 
 \\
 &\geq N_l - g + 1 - (n_+ + n_- + n_{\mathrm{non}})  
 & \hspace{-160pt} \text{(by $n_{\pm} \geq \nu_{\pm}$)} 
 \\
 &\geq N_l - g + 1 - N_s - n_{\mathrm{non}}  
 &  \hspace{-120pt}\text{(by $N_s \geq n_+ + n_-$)} 
 \\
 &\geq N_l - g + 1 - N_s + n_{\mathrm{non}}/2 - 3g/2  \hspace{-160pt} {}
 &  \text{(by $g \geq n_{\mathrm{non}}$)} 
 \\
 &= N_l - 5g/2 - N_s + n_{\mathrm{non}}/2 + 1 & \hspace{-160pt}  \\
 &\geq \dfrac{3N_s}{2} + \frac{(N_s - \chi(S) + 2)(N_s - \chi(S))+ N_\partial (S) }{2} + k + n_{\mathrm{non}}/2  + 1 
 &  \hspace{-200pt} \text{(by \eqref{eq:inequality_var})} 
 \\
 &\geq N_s + \frac{n_+ + (N_s - \chi(S) + 2)(N_s - \chi(S)) + n_{\mathrm{non}} - N_{\partial, \mathrm{nonA}}}{2} & \hspace{-200pt} 
 \\
 & \hspace{10pt} + N_\partial (S)/2 + N_{\partial, \mathrm{nonA}}/2 + k + 1 
 & \hspace{-170pt}\text{(by $N_s \geq n_+ + n_-$)} 
 \\
 &\geq N_s + k+1 + N_\partial (S)/2 + N_{\partial, \mathrm{nonA}}/2
& \hspace{-170pt} \text{(by \eqref{eq:bound_N_partial})}
\end{align*}
and so $n_t - N_\partial (S)/2 - N_{\partial, \mathrm{nonA}}/2 \geq N_s + k+1$. 
Therefore, there are $N_s + k+1$ triples $(\mu_{n_j}, U_{i_j}, U_{i'_j})_{j= 0}^{N_s+ k}$ such that $\mu_{n_j}$ is a closed transversal and the connected components $U_{i_j}$ and $U_{i'_j}$ are annuli  for any $j \in \{0, \ldots , N_s+ k\}$. 
For any $j \in \{0, \ldots , N_s+ k\}$, since the annulus $U_{i_j}$ has disjoint boundary components in $S$, the disjoint union $\mathbb{A}_j := U_{i_j} \sqcup \mu_{n_j} \sqcup  U_{i'_j}$ is an annulus, which is either closed, half-open, or open. 
Since $N_s = \vert \Sv \vert$, by removing $N_s$ triples from $(\mu_{n_j}, U_{i_j}, U_{i'_j})_{j= 0}^{N_s+ k}$, we may assume that $(\mu_{n_j}, U_{i_j}, U_{i'_j})_{j= 0}^{k}$ satisfies $\Sv \cap \bigsqcup_{j=0}^k \mathbb{A}_j = \emptyset$. 

\begin{claim}\label{claim:10}
For any $j \in \{0, \ldots , k\}$, the $\alpha$-limit set of any point in $\mu_{n_j}$ is contained in the annulus $\mathbb{A}_j$.
\end{claim}

\begin{proof}[Proof of Claim~\ref{claim:10}]
Suppose that $\A_j$ is a closed annulus. 
Since $S$ is connected, we have that $\A_j = S$ and so that $\A_j$ is invariant. 
Then the $\alpha$-limit set of a point $x \in \mu_{n_j}$ is contained in the annulus $\mathbb{A}_j$. 

Suppose that $\A_j$ is a half-open annulus. 
Then $\mu := \overline{\A_j} - \A_j \subset S - \partial S$ is a closed transversal. 
By transposing $U_{i_j}$ and $U_{i'_j}$ if necessary, we may assume that $\mu$ is a boundary component of $U_{i_j}$. 
Since $\mu$ and $\mu_{n_j}$ are boundary components of semi-attracting collar basins $A$ and $A_{n_j}$ of distinct semi-attracting limit cycles respectively, by the connectivity of the annulus $v(A_{n_j}) = v(\mu_{n_j})$, we obtain that $v(A_{n_j}) = v(\mu_{n_j}) \subset \A_j \setminus v(\mu) = \A_j \setminus v(A)$ and so that $\overline{v(A_{n_j})} \subseteq \A_j$.
Then the $\alpha$-limit set of a point $x \in \mu_{n_j}$ is contained in the annulus $\mathbb{A}_j$. 

Thus, we may assume that $\A_j$ is an open annulus. 
Then the boundary $\partial \A_j$ consists of disjoint two closed transversals $\mu_{i_j}$ and $\mu_{i'_j}$. 
The disjoint union $U_{\partial \mathbb{A}_j} := A_{i_j} \sqcup \mu_{i_j} \sqcup A_{i'_j} \sqcup \mu_{i'_j}$ is a \nbd of  $\partial \mathbb{A}_j = \mu_{i_j} \sqcup \mu_{i'_j}$ in the closed annulus $\overline{\mathbb{A}_j}$. 
Since $(A_n \sqcup \mu_n) \cap v(A_m) = \emptyset$ for any $n \neq m$, we have that $U_{\partial \mathbb{A}_j} \cap v(A_{n_j}) 
= \emptyset$. 
Since the open annulus $\mathbb{A}_j$ intersects the open annulus $A_{n_j}$, the connectivity of $v(A_{n_j})$ implies that $v(A_{n_j}) \subseteq \mathbb{A}_j - U_{\partial \mathbb{A}_j}$ and so that $\overline{v(A_{n_j})} \subset \mathbb{A}_j$. 
Then the $\alpha$-limit set of a point $x \in \mu_{n_j}$ is contained in the open annulus $\mathbb{A}_j$. 
\end{proof}

For any $j \in \{0, \ldots , k\}$, by the non-existence of singular points on $\A_j$,
the dual statement of Lemma~\ref{lem:annulus_nonrec} implies that the $\alpha$-limit set $\alpha(x) \subset \mathbb{A}_j$ for any point $x \in \mu_{n_j}$ is a limit cycle. 
This means that $N_{l,\alpha} \geq k+1$.
\end{proof}

Theorem~\ref{main:fininite_existence_limit_circuits} and the previous lemma imply the following dependency. 

\begin{proof}[Proof of Theorem~\ref{main:inifinie_limit_sets}]
From Lemma~\ref{lem:Maier}, by the finiteness of singular points, the total number of Q-sets for $v$ is finite and so assertions {\rm(1)} and {\rm(2)} (resp. {\rm(4)} and {\rm(5)}) are equivalent.  
By Theorem~\ref{main:fininite_existence_limit_circuits}, assertions {\rm(2)} and {\rm(3)} (resp. {\rm(5)} and {\rm(6)}) are equivalent. 
Lemma~\ref{lem:counter_alpha_limit_cycles} and the dual statement imply that assertions {\rm(3)} and {\rm(6)} are equivalent. 
\end{proof}

By considering the negation of assertions of Theorem~\ref{main:inifinie_limit_sets}, we obtain the following characterization of the existence of finitely many limit sets. 

\begin{theorem}\label{main:ifinite_limit_sets}
The following statements are equivalent for a flow $v$ with finitely many singular points on a compact surface:
\\
{\rm(1)} There are at most finitely many pairwise distinct $\omega$-limit sets of non-closed points which are limit cycles. 
\\
{\rm(2)} There are at most finitely many pairwise distinct $\omega$-limit sets of non-closed points which are limit circuits. 
\\
{\rm(3)} There are at most finitely many pairwise distinct $\omega$-limit sets of non-closed points. 
\\
{\rm(4)} There are at most finitely many pairwise distinct $\alpha$-limit sets of non-closed points which are limit cycles. 
\\
{\rm(5)} There are at most finitely many pairwise distinct $\alpha$-limit sets of non-closed points which are limit circuits. 
\\
{\rm(6)} There are at most finitely many pairwise distinct $\alpha$-limit sets of non-closed points. 
\end{theorem}

\subsubsection{Equivalence of existence of finitely many $\omega$-limit sets and $\alpha$-limit sets}

The previous result implies Theorem~\ref{main:ch_finite_limit_sets} as follows. 

\begin{proof}[Proof of Theorem~\ref{main:ch_finite_limit_sets}]
Theorem~\ref{main:ifinite_limit_sets} implies that assertions (1) and (2) are equivalent. 
If there are at most finitely many pairwise distinct $\omega$-limit sets of points, then there are at most finitely many closed orbits.

Conversely, suppose that there are at most finitely many closed orbits.
Then there are at most finitely many pairwise distinct $\omega$-limit sets of non-closed points which are limit cycles. 
Theorem~\ref{main:ifinite_limit_sets} implies that there are at most finitely many pairwise distinct $\omega$-limit sets and $\alpha$-limit sets of non-closed points. 
The finiteness of closed orbits implies that there are at most finitely many pairwise distinct $\omega$-limit sets and $\alpha$-limit sets of points. 
\end{proof}

\section{Topological characterization of the denseness of periodic orbits in the non-wandering set}

In this section, we provide a necessary and sufficient condition for the denseness of closed orbits in the non-wandering set for a flow with finitely many singular points on compact surfaces. 
This section aims to establish Lemma~\ref{lem:wandering_h} as an intermediate step toward Theorem~\ref{main:rec_wandering_pt}, and to derive 
Theorem~\ref{thm:finite_dense} as a corollary. 

We recall the following statement, which is used in the proofs of several results of this section.

\begin{lemma}[Lemma~2.3 \cite{yokoyama2016topological}]\label{lem:top23}
 Let $v$ be a flow on a compact surface $S$. 
 Then $\overline{\Cv \sqcup \mathrm{LD}(v)} \cap \mathrm{E}(v) = \emptyset$ and $\overline{\Cv \sqcup \mathrm{E}(v)} \cap \mathrm{LD}(v) = \emptyset$. 
Moreover, we have $\mathrm{E}(v) \subseteq \mathop{\mathrm{int}} \overline{\mathrm{P}(v)}$.  
\end{lemma}

Note that the compactness in the previous lemma is necessary (see Example ~\ref{ex:counter_exam01+} for details). 

\subsection{Necessary conditions}
We present the following necessary conditions for the denseness of closed orbits with respect to a flow on a surface, which serve as an intermediate step toward Theorem~\ref{thm:finite_dense} and are in fact a special case of it.

\begin{proposition}\label{prop:necessary_dense}
Let $v$ be a flow with finitely many singular points on a compact surface. 
If $\overline{\mathop{\mathrm{Cl}}(v)} = \Omega(v)$, then there are neither non-closed recurrent points, nor strictly limit circuits, nor circuits with wandering holonomy.
\end{proposition}

To show the previous proposition, we show the following lemma. 

\begin{lemma}\label{thm41}
Let $v$ be a flow on a compact surface $S$.
If $\overline{\mathop{\mathrm{Cl}}(v)} = \Omega(v)$, then $\mathrm{R}(v) = \emptyset$ $( \mathrm{i.e.}$ $S = \mathop{\mathrm{Cl}}(v) \sqcup \mathrm{P}(v)$ $)$. 
\end{lemma}

\begin{proof}
Suppose that $\overline{\mathop{\mathrm{Cl}}(v)} = \Omega(v)$. 
Since recurrent orbits are non-wandering, we have $\mathrm{R}(v) \subseteq  \Omega(v) = \overline{\mathop{\mathrm{Cl}}(v)}$.
From the compactness of $S$, Lemma~\ref{lem:top23} implies that $\overline{\mathop{\mathrm{Cl}}(v)} \cap \mathrm{R}(v) = \emptyset$, and hence $\mathrm{R}(v) = \emptyset$.
By Lemma~\ref{lem:decomp}, we have $S = \mathop{\mathrm{Cl}}(v) \sqcup \mathrm{P}(v)$. 
\end{proof}

The converse of Lemma~\ref{thm41} is not true in general (see Lemma~\ref{lem:counter_exam00} for details). 
We have the following property.

\begin{lemma}\label{lem:wandering_hol_one-sieded}
Let $v$ be a flow with finitely many singular points on a compact surface $S$ and $x \in S$ a non-singular point.  
If there is a strictly limit circuit which is one-sided at $O(x)$, then $x \in \Omega(v) - \overline{\Cv \sqcup \mathrm{R}(v)}$.
\end{lemma}

\begin{proof}
The finiteness of singular points implies that $x \in S - \Sv = S - \overline{\Sv}$. 
Let $\gamma$ be a strictly limit circuit which is one-sided at $O(x)$. 
By time reversion if necessary, we may assume that there is a point $y \in S$ with $\gamma = \omega(y)$. 
Then $x \in \gamma = \omega(y) \subseteq \Omega(v)$. 
Moreover, the non-periodic non-trivial limit circuit $\gamma$ consists of connecting separatrices and finitely many singular points and $\emptyset \neq \gamma \setminus \Sv \subseteq \mathrm{P}(v)$. 
By Lemma~\ref{lem:lc_transversal}, there is a semi-attracting collar basin $A_0$ of the limit circuit $\gamma$ such that $\gamma$ is a boundary component of $A_0$. 
Since $\gamma$ is one-sided at $O(x)$, there is a small collar $A \subset A_0$ of $\gamma$ such that the union $A \sqcup \gamma$ is a neighborhood of the point $x \in \mathrm{P}(v) \cap \gamma$. 
%
The smallness of $A$ implies that we may assume that $\partial A$ consists of two connected components. 
Since $\gamma$ is a boundary component of $A$, let $\mu$ be another connected component of $\partial A$. 
Then $\gamma \sqcup \mu = \partial A$.
By Lemma~\ref{lem:existence_closed_transversal}, we may assume that $\mu$ is a closed transversal. 
Since $\gamma$ consists of finitely many orbits in $\Sv \sqcup \mathrm{P}(v)$, the union $A \sqcup O(x)$ is also a neighborhood of $x$. 

\begin{claim}\label{claim:13}
$A \cap (\Cv \sqcup \mathrm{R}(v)) = \emptyset$. 
\end{claim}

\begin{proof}[Proof of Claim~\ref{claim:13}]
Fix any point $z \in A$.
Since $A$ is a semi-attracting collar basin of $\gamma$, the $\omega$-limit set  $\omega(z)$ is the limit circuit $\gamma \subseteq \Sv \sqcup \mathrm{P}(v)$. 
This implies that $z \notin \Cv$. 
Therefore, we obtain $A \cap \Cv = \emptyset$. 
By Lemma~\ref{lem:annulus_nonrec}, if $O(z)$ is contained in the open annulus $A$, then $x \notin \mathrm{R}(v)$. 
Thus, we may assume that $O(z)$ is not contained in the open annulus $A$. 
Then $O(z) \cap \partial A \neq \emptyset$. 
Since $\partial A$ consists of the invariant subset $\gamma$ and the closed transversal $\mu$, there is a point $z' \in O(z) \cap \mu$ such that $O^-(z') \cap (A \sqcup \mu) = \emptyset$ and $O^+(z') \subset A$. 
Because $\mu$ is the closed transversal and the set difference $\partial A - \mu = \gamma$ is invariant, we have that $(\alpha(z') \sqcup \omega(z')) \cap \mu = (\alpha(z') \sqcup \gamma) \cap \mu = \alpha(z') \cap \mu = \emptyset$. 
This implies that $z' \in \mu$ is non-recurrent and so is $z \in O(z')$. 
Therefore, we obtain $A \cap \mathrm{R}(v) = \emptyset$. 
\end{proof}

By the previous claim, the union $A \sqcup O(x)$ is a neighborhood of $x$ with $(A \sqcup O(x)) \cap (\Cv \sqcup \mathrm{R}(v)) = \emptyset$. 
This implies that $x \notin \overline{\Cv \sqcup \mathrm{R}(v)}$. 
Since $x \in \Omega(v)$, we have that $x \in \Omega(v) - \overline{\Cv \sqcup \mathrm{R}(v)}$. 
\end{proof}

\begin{lemma}\label{lem:slc_01}
Let $v$ be a flow with finitely many singular points on a compact surface $S$.
Every strictly limit circuit $\gamma$ contains a point $x$ in $\Omega(v) - \overline{\mathop{\mathrm{Cl}}(v) \sqcup \mathrm{R}(v)}$ such that $\gamma$ is a strictly limit circuit at $O(x)$. 
\end{lemma}

\begin{proof}
Let $\gamma$ be a strictly limit circuit of the flow $v$. 
Lemma~\ref{lem:wandering_hol_one-sieded} implies that we may assume that $\gamma$ is not one-sided at $O(x)$. 
Then there are a separatrix $\mu \subseteq \gamma$ and a transverse open arc $T \subset S - \mathrm{R}(v)$ intersecting $\mu$ uniquely at the point $x$ such that $f_v|_{T_1}$ is either attracting or repelling, and that $f_v|_{T_2}$ has no periodic points, where $f_v$ is the first return map on $T$, the subsets $T_1$ and $T_2$ are the two connected components of $T \setminus \mu$.
By time reversion if necessary, we may assume that $f_v|_{T_1}$ is attracting. 
Then there is a point $y \in T_1$ with $\omega(y) = \gamma$. 
Since $\gamma$ is a non-periodic non-trivial circuit, the circuit $\gamma$ is a directed graph which is the union of connecting separatrices and finitely many singular points. 
Then $\gamma \subseteq (\mathop{\mathrm{Sing}}(v) \sqcup \mathrm{P}(v)) \cap \Omega(v)$ and so $x \in \mu \cap T \subset \gamma \cap \mathrm{P}(v) \cap \Omega(v)$. 
Since $T \subset S - \mathrm{R}(v)$ is a transverse open arc, we have $x \notin \overline{\mathrm{R}(v)}$. 

\begin{claim}\label{claim:0a}
$x \notin \overline{\mathop{\mathrm{Cl}}(v)}$.
\end{claim}
\begin{proof}[Proof of Claim~\ref{claim:0a}]
By the flow box theorem, there is an open \nbd $U$ of $x$ containing no singular points.
Assume that $x \in \overline{\mathop{\mathrm{Cl}}(v)}$.
Since $U$ has no singular points, the closedness of the singular point set implies $x \in \overline{\mathop{\mathrm{Per}}(v)}$. 
Because $f_v|_{T_1}$ is attracting, we have $x \in \overline{\mathop{\mathrm{Per}}(v) \cap T_2}$, which contradicts that $f_v|_{T_2 \cap \mathrm{dom}(f_v)}$ has no periodic points. 
\end{proof}

By the previous claim, we obtain that $x \in \gamma \cap (\Omega(v) - \overline{\mathop{\mathrm{Cl}}(v) \sqcup \mathrm{R}(v)})$. 
\end{proof}

By Lemmas~\ref{lem:wandering_hol} and \ref{lem:slc_01}, the following statement holds.

\begin{proposition}\label{prop:056}
Let $v$ be a flow with finitely many singular points on a compact surface. 
If $\Omega(v) = \overline{\mathop{\mathrm{Cl}}(v) \sqcup \mathrm{R}(v)}$, then there are neither strictly limit circuits nor circuits with wandering holonomy.
\end{proposition}

Therefore, Proposition~\ref{prop:necessary_dense} follows from the previous proposition together with Lemma~\ref{thm41}.
Note that the converse of Proposition~\ref{prop:necessary_dense} is not true without the assumption of finitely many singular points (see Lemma~\ref{lem:counter_exam01} for details). 
Moreover, we have the following statement. 

\begin{corollary}\label{prop:necessary_dense02}
Let $v$ be a flow with finitely many singular points on a compact surface. 
If $\mathop{\mathrm{Cl}}(v) = \Omega(v)$, then $\mathop{\mathrm{Cl}}(v)$ is closed, and there are neither non-closed recurrent points, nor non-periodic limit circuits, nor circuits with wandering holonomy.
\end{corollary}

\begin{proof}
Suppose that $\mathop{\mathrm{Cl}}(v) = \Omega(v)$. 
Since $\Omega(v)$ is closed, the closed point set $\mathop{\mathrm{Cl}}(v)$ is closed with $\overline{\mathop{\mathrm{Cl}}(v)} = \mathop{\mathrm{Cl}}(v) = \Omega(v)$. 
Proposition~\ref{prop:056} implies that there are neither non-closed recurrent points nor circuits with wandering holonomy.
Since any point in $\omega$-limit and $\alpha$-limit sets of any point is non-wandering and so closed, there are no non-periodic limit circuits. 
\end{proof}

\subsection{Sufficient conditions}

In this subsection, under the hypothesis of the finiteness of singular points, we show that any non-singular point in the boundary of the periodic point set is contained in a circuit.
To show this, we first make the following observation.

\begin{lemma}\label{lem:circuit-}
Let $v$ be a flow with finitely many singular points on a compact surface $S$ and $I$ an open transverse arc whose boundary contains a non-singular point $x_\infty \in S$. 
Suppose that there is a sequence $(x_i)_{i \in \mathbb{N}}$ of fixed points of the first return map $f_I$ on $I$ generated by $v$ that satisfies the following three conditions: 
\\
{\rm(1)} The sequence $(x_i)_{i \in \mathbb{N}}$ is monotone on $I$ and converges to $x_\infty$ such that the orbits $O(x_i)$ and $O(x_j)$ are disjoint for any $i \neq j$.
\\
{\rm(2)} For any $i,j \in \mathbb{N}$, loops $O(x_i)$ and $O(x_j)$ are freely homotopic to each other.
\\
{\rm(3)}  For any $i \in \mathbb{N}$, there is an invariant open annulus $\A_i$ containing $(x_i,x_{i+1})_I$ and intersecting no singular points of $v$ such that $\A_i$ is a connected component of $S - \bigsqcup_{j \in \mathbb{N}} O(x_j)$ and that $\partial \A_i= O(x_i) \sqcup O(x_{i+1})$. 

Then the boundary component of $\A := \bigsqcup_{i \in \mathbb{N} } \A_i \sqcup O(x_{i+1})$ containing $x_\infty$ is a circuit in $\overline{\Pv} \subseteq S - \mathrm{R}(v)$. 
\end{lemma}

\begin{proof}
By Gutierrez's smoothing theorem, 
we may assume that $v$ is $C^1$. 
By the construction, the union $\A = \bigsqcup_{i \in \mathbb{N} } \A_i \sqcup O(x_{i+1})$ is an invariant open annulus intersecting no singular points of $v$. 
Take a subinterval $J := \bigsqcup_{i \in \mathbb{N}} (x_i,x_{i+1})_I \sqcup \{ x_{i} \} \subseteq I$ and the first return map $f_J$ on $J$ generated by $v$. 
Then $J = (x_1, x_\infty)_I$ and the point $x_i$ for any $i \in \mathbb{N}$ is a fixed point with respect to $f_J$. 
Let $\mu$ be the boundary component of $\A$ containing $x_\infty$. 
Then $\mu \subset \overline{\bigcup_{j \in \mathbb{N}} O(x_j)} \subseteq \overline{\Pv}$. 
Lemma~\ref{lem:decomp} and Lemma~\ref{lem:top23} imply that $\mathrm{R}(v) \cap \overline{\mathop{\mathrm{Per}}(v)} = \emptyset$ and so that $\mathrm{R}(v) \cap \mu = \emptyset$.

Thus, it suffices to show that the boundary component $\mu$ is a circuit with respect to $v$. 
Indeed, 
since $v$ is $C^1$, there is an integrable continuous vector field $X$ generating $v$. 
Because $J$ is transverse to $v$, by the flow box theorem (cf. \cite[Theorem~4.2.6, p.95]{markley2023flows}, \cite[Theorem 1.1, p.45]{aranson1996introduction}), there is a trivial flow box $U$ with $J \subset U$ and $x_\infty \in \partial U$. 
Then we can identify the restriction $X \vert_U$ of the vector field $X$ with a unit vector field $(1,0) = \partial/\partial x_1$ on a closed disk $[0,1]^2$ as on the left in Figure~\ref{fig:flowbox_perturbation_att}. 
Since the union $\A$ is an invariant open annulus, for any $i \in \N$, the restriction $f_J$ near $x_i$ is orientation-preserving. 
For any $i \in \N$, the closure $\overline{\A_i}$ contains no singular points and is an invariant closed annulus whose boundary $\partial \A_i = O(x_i) \sqcup O(x_{i+1})$ is a disjoint union of two periodic orbits and which contains the closed interval $[x_i,x_{i+1}]_I$ connecting the boundary components.

\begin{claim}\label{claim:0b}
For any $i \in \N$, the restriction $f_J\vert_{[x_i,x_{i+1}]_I}$ is a homeomorphism on the closed interval $[x_i,x_{i+1}]_I$ such that  $f_J([x_i,x_{i+1}]_I) = [x_i,x_{i+1}]_I$. 
\end{claim}

\begin{proof}[Proof of Claim~\ref{claim:0b}]
Fix any $i \in \N$.
Since the Poincar{\'e}-Hopf theorem implies the fact that any closed disk whose boundary is a periodic orbit contains singular points, applying Lemma~\ref{lem:annulus_nonrec} to $\A$ and the dual statment, the non-existence of singular points on $\A_i$ implies that each of the $\omega$-limit set and the $\alpha$-limit set of any point in $\A_i$ is a limit cycle parallel to the boundary components of $\overline{\A_i}$. 
This implies that the restriction $f_J\vert_{[x_i,x_{i+1}]_I}$ and its inverse map $(f_J\vert_{[x_i,x_{i+1}]_I})^{-1}$ are well-defined continuous surjections. 
Since the inverse map $(f_J\vert_{[x_i,x_{i+1}]_I})^{-1}$ is well-defined, the restriction $f_J\vert_{[x_i,x_{i+1}]_I}$ is injective.
Therefore, the restriction $f_J\vert_{[x_i,x_{i+1}]_I}$ is a homeomorphism on the closed interval $[x_i,x_{i+1}]_I$. 
\end{proof}

\begin{figure}
\begin{center}
\includegraphics[scale=0.7]{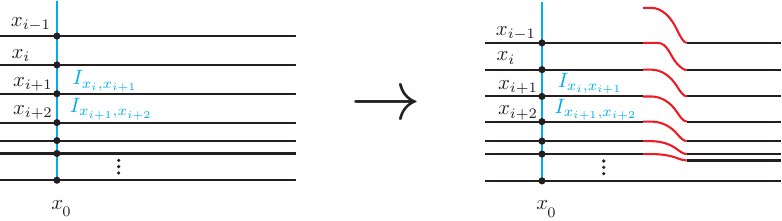}
\end{center}
\caption{Perturbation of a unit vector field $(1,0)$ in the flow box $U$ using a bump function, and orbits generated by the resulting vector field as in the right.}
\label{fig:flowbox_perturbation_att}
\end{figure} 
By using a bump function $\varphi \colon U \to \R_{\geq 0}$ on the trivial flow box $U = [0,1]^2$ and replacing the restriction $X \vert_U$ of the vector field $X$ with $X \vert_U - \varphi \partial/\partial x_2$ which generates orbits as on the right in Figure~\ref{fig:flowbox_perturbation_att}, we obtain the flow $w$ generated by the resulting vector field satisfying the following conditions: 
\\
{\rm(1)} The flow $w$ has finitely many singular points. 
\\
{\rm(2)} The first return map $f_{w,J}$ on $J = (x_1, x_\infty)_I$ with respect to $w$ is semi-attracting toward $x_\infty$.
\\
{\rm(3)} The restriction $v|_{S - \A}$ is topologically equivalent to the restriction $w|_{S - \A}$ via the identity map on $S - \A$. 
\\
{\rm(4)} $\A \cap \mathop{\mathrm{Sing}}(w) = \emptyset$.
\\ 
The invariance of $\A$ with respect to $v$ and $w$ implies that $\mathrm{R}(w) = \mathrm{R}(v) \subset S - \A$.
Moreover, the $\omega$-limit set $\omega_w(x)$ of a point $x \in J$ with respect to $w$ contains the non-singular point $x_\infty$ such that $\omega_w(x) \subseteq \mu \subset S - \mathrm{R}(w) = S - \mathrm{R}(v)$.
By Lemma~\ref{lem:lc_transversal}, the $\omega$-limit set $\omega_w(x)$ is a limit circuit with respect to $w$. 

\begin{claim}\label{claim:bdy}
$\mu = \omega_w(x)$. 
\end{claim}

\begin{proof}[Proof of Lemma~\ref{claim:bdy}]
By Claim~\ref{claim:0b}, the set difference $\A - J$ is an open trivial flow box with respect to $w$.
Since the surface $S$ is Fr{\'e}chet-Urysohn, the attraction of $f_{w,J}$ implies that any point in $\mu$ is contained in $\overline{O^+(x)} - (O^+(x) \sqcup \{ x \}) \subseteq \omega_w(x)$. 
This means that $\mu = \omega_w(x)$. 
\end{proof}

Claim~\ref{claim:bdy} implies that $\mu = \omega_w(x) \subset S - \A$ is a circuit with respect to $w$. 
Since the restriction $v|_{S - \A}$ corresponds to the restriction $w|_{S - \A}$, the boundary component $\mu$ of $\A$ is the circuit with respect to $v$. 
\end{proof}

To show a sufficient condition for being contained in a circuit, we have the following observation.

\begin{figure}[t]
\begin{center}
\includegraphics[scale=0.5]{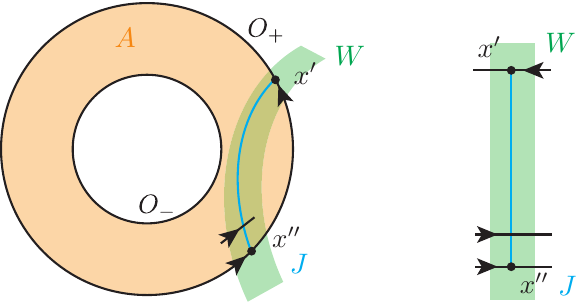}
\end{center}
\caption{The annulus $A$ and the open arc $J \subset A$ with $\partial J \subset O_+$.}
\label{fig:trans_orientation}
\end{figure} 

\begin{lemma}\label{lem:inh_ori}
Let $v$ be a flow on a surface $S$ and $A \subseteq S$ an invariant closed annulus whose boundary consists of disjoint two periodic orbits $O_-$ and $O_+$. 
Then there is no closed transverse arc $J \subset A$ with either $\partial J \subset O_-$ or $\partial J \subset O_+$. 
\end{lemma}
 
\begin{proof}
Assume that there is a closed transverse arc $J \subset A$ with $\partial J \subset O_+$ as on the left in Figure~\ref{fig:trans_orientation}.
By the flow box theorem, there is an open flow box $W$ containing $J$ as in  Figure~\ref{fig:trans_orientation}. 
Denote by $x'_i$ and $x''_i$ the points in $\partial J$. 
The flow $v$ defines inherent directions on $O_+$ at $x''$.
However, the trivial flow box $W$ containing $J$ implies a direction at $x'$ that is opposite to the direction induced by following the orbit arc from $x''$ to $x'$, which is a contradiction.
Therefore, there is no open arc $J \subset A$ with $\partial J \subset O_+$. 
By the same argument for $O_-$, the assertion holds. 
\end{proof}

The following lemma provides a sufficient condition for being contained in a circuit.

\begin{lemma}\label{lem:circuit}
Let $v$ be a flow with finitely many singular points on a compact surface $S$. 
Each non-singular point in $\overline{\Pv}$ is contained in a circuit in $\overline{\Pv} \subseteq S - \mathrm{R}(v)$. 
\end{lemma}

\begin{proof}
By Lemma~\ref{lem:top23}, we have that $\overline{\Pv} \subseteq S - \mathrm{R}(v)$. 
Fix a non-singular point $x_\infty \in \overline{\Pv}$. 
If $x_\infty \in \Pv$, then the orbit of $x_\infty$ is a periodic circuit, which implies the assertion. 
Thus, we may assume that $x_\infty$ is not periodic.  
By taking the lift of the flow $v$ to the double of $S$ if necessary, we may assume that the point $x_\infty$ is not contained in the boundary $\partial S$. 
From the flow box theorem, there is an open flow box $D$ centered at $x_\infty$. 
Then there is a small closed transverse arc $I \subset D$ whose boundary contains $x_\infty$, and there is a sequence $(x_i)_{i \in \N}$ of periodic points of the first return map $f_I$ in $\mathop{\mathrm{int}} I$ (i.e. $x_i \in \mathop{\mathrm{Per}}(f_I) \cap \mathop{\mathrm{int}} I$) which monotonically converges to $x_\infty$ on $I$.
We define a natural total order on the interval $I$ such that $x_\infty$ is the minimal element. 
With respect to this total order, the sequence $(x_i)_{i \in \N}$ converges monotonically to $x_\infty$.
Since any periodic orbit intersects the closed arc $I$ at most finitely many times, by replacing $x_i $ with a point in $O(x_i) \cap I$, by taking a subsequence of $(x_i)_{i \in \N}$, we may assume that any $x_i$ is minimal in the points in $O(x_i) \cap I$. 
Then $O(x_i) \cap [x_{i},x_{\infty}]_I = \{ x_i\}$ for any $i \in \N$. 
Since $S$ is of finite genus, by taking a subsequence of $(x_i)_{i \in \N}$, we may assume that, for any $i,j \in \N$, loops $O(x_i)$ and $O(x_j)$ are disjoint and freely homotopic to each other.

To apply Lemma~\ref{lem:circuit-} to a subsequence of $(x_i)_{i \in \N}$, it suffices to show that there is a subsequence $(x_{i_k})_{k \in \N}$ of $(x_i)_{i \in \N}$ which consists of fixed points of the first return map $f_{(x_{i_1},x_\infty)_I}$ on $(x_{i_1},x_\infty)_I$ generated by $v$ and satisfies condition (3) in Lemma~\ref{lem:circuit-}, where  $(x_{i_0},x_\infty)_I$ and $(x_{i_k})_{k \in \N}$ serve as $I$ and $(x_i)_{i \in \N}$ satisfying the conditions of Lemma~\ref{lem:circuit-}.
From Lemma~\ref{claim:003a}, by taking a subsequence of $(x_i)_{i \in \N}$, we may assume that there is a sequence $(\A_i)_{i \in \N}$ of pairwise disjoint open annuli $\A_{i}$ whose boundary is $O(x_i) \sqcup O(x_{i+1})$ satisfying the following properties: 
\begin{enumerate}[\rm (P1)]
\item The union $A := \bigcup_{k=1}^\infty (\A_{k} \sqcup O(x_{k+1}))$ is an open annulus one of whose boundary components is the loop $O(x_1)$. 
\item For any $k \in \N$, the loop $O(x_{k+1})$ and the annulus $\A_{k}$ are parallel to the boundary of $A$.
\item The set of connected components of $S - \bigcup_{k=1}^\infty O(x_k)$ contained in $A$ coincides with the family $\{ \A_{k}\}_{k \in \N}$ of annuli.
\item For any $k<k' \in \N$, the union $A_{k,k'} := \A_{k} \sqcup \bigsqcup_{l = k+1}^{k'-1} (O(x_l) \sqcup \A_{l}) \subset A$ is an open annulus whose boundary is $O(x_k) \sqcup O(x_{k'})$.
\end{enumerate}

\begin{claim}\label{claim:014}
For any positive numbers $i < j \in \N$, there is no open annulus whose boundary is $O(x_i) \sqcup O(x_j)$ containing $x_\infty$. 
\end{claim}
\begin{proof}[Proof of Claim~\ref{claim:014}]
Assume that there are numbers $M < N \in \N$ and an open annulus $A_{M,N}$ whose boundary is $O(x_M) \sqcup O(x_N)$ such that $x_\infty \in A_{M,N}$. 
Since $(x_i)_{i \in \mathbb{N}}$ converges to $x_\infty$, by taking a subsequence of $(x_i)_{i \in \mathbb{N}}$, we may assume $x_i \in A_{M,N}$ for any $i > N \in \mathbb{N}$. 
By the invariance of $A_{M,N}$, we have $O(x_i) \subset A_{M,N}$ for any $i > N \in \mathbb{N}$.
For any $i > N \in \mathbb{N}$, since $O(x_i) \subset \partial \A_{i}$, the intersection $\emptyset \neq A_{M,N} \cap \A_{i}$ is a nonempty open subset and so there is an integer $i'$ with $M \leq i' \leq N$ such that $\A_{i'} \cap \A_i \neq \emptyset$, which contradicts the pairwise disjoint property of $(\A_i)_{i \in \N}$. 
\end{proof}

Applying Lemma~\ref{lem:inh_ori} to the closed annulus $\overline{A_{i,j}}$, for any positive numbers $i < j \in \N$, there is no open arc $J \subset I \cap A_{i,j}$ with either $\partial J \subset O(x_i)$ or $\partial J \subset O(x_j)$. 

\begin{claim}\label{claim:016-}
The following statements hold for any $i< j \in \N$: 
\begin{enumerate}[\rm (1)]
\item $(x_{i},x_{j})_I = (x_1,x_{\infty}) \cap A_{i,j}$. 
\item $\{ x_j \} = (x_1,x_{\infty}) \cap O(x_j)$. 
\end{enumerate}
\end{claim}

\begin{proof}[Proof of Claim~\ref{claim:016-}]
Fix any $i< j \in \N$.
Claim~\ref{claim:014} implies that $x_\infty \notin A_{i,j}$. 
By $O(x_k) \cap [x_{k},x_{\infty}]_I = \{ x_k \}$ for any $k \in \N$, we obtain that $[x_{i},x_\infty]_I \cap O(x_{i}) = \{ x_{i} \}$ and $[x_{j},x_\infty]_I \cap O(x_{j}) = \{ x_{j} \}$. 
From $(x_{j},x_\infty)_I \subset (x_{i},x_\infty)_I$ and $\partial A_{i,j} = O(x_i) \sqcup O(x_{j})$, we have that  $[x_{j},x_\infty]_I \cap O(x_{i}) = \emptyset$ and so that $[x_{j},x_\infty]_I \cap \partial A_{i,j} = [x_{j},x_\infty]_I \cap (O(x_i) \sqcup O(x_{j})) = \{ x_{j} \} $. 
If $A_{i,j} \cap (x_{j},x_{\infty})_I \neq \emptyset$, then $[x_{j},x_{\infty}]_I - \{ x_j \} \subset A_{i,j}$ as shown on the right in Figure~\ref{fig:para_arc}
and so $x_\infty \in [x_{j},x_{\infty}]_I - \{ x_j \} \subset A_{i,j}$, which contradicts $x_\infty \notin A_{i,j}$. 
\begin{figure}
\begin{center}
\includegraphics[scale=1.15]{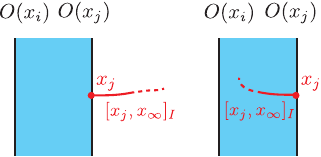}
\end{center}
\caption{Left: A possible pattern; Right: Impossible patterns.}
\label{fig:para_arc}
\end{figure} 
Thus, we obtain that $A_{i,j} \cap (x_{j},x_{\infty})_I = \emptyset$. 
By the transversality of $(x_{j},x_{\infty})_I$, there is a point $y \in (x_{i},x_{j})_I$ with $(y,x_{j})_I \subset A_{i,j}$. 
In particular, we have $A_{i,j} \cap (x_{i},x_{j})_I \neq \emptyset$. 
By $x_i \in O(x_i) \cap [x_{i},x_{j}]_I \subseteq O(x_i) \cap [x_{i},x_{\infty}]_I = \{ x_i \}$, we have that $O(x_i) \cap [x_{i},x_{j}]_I = \{ x_i \}$. 

Assume that $\partial A_{i,j} \cap (x_{i},x_{j})_I \neq \emptyset$. 
From $\partial A_{i,j} = O(x_i) \sqcup O(x_j)$ and $O(x_i) \cap (x_{i},x_{j})_I = \emptyset$, we have $O(x_j) \cap (x_{i},x_{j})_I = \partial A_{i,j} \cap (x_{i},x_{j})_I \neq \emptyset$.  
By the transversality of $(x_{i},x_{j})_I$, there is a point $z \in \partial A_{i,j} \cap (x_{i},x_{j})_I = O(x_j) \cap (x_{i},x_{j})_I$ as shown on the right in Figure~\ref{fig:para_arc_001}
\begin{figure}
\begin{center}
\includegraphics[scale=1.05]{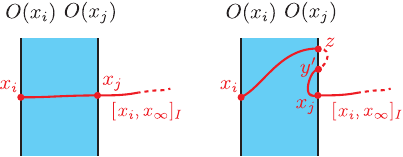}
\end{center}
\caption{Left: A possible pattern; Right: Impossible pattern.}
\label{fig:para_arc_001}
\end{figure} 
and there is a point $z' \in (z,y)_I$ with $(z,z')_I \cap A_{i,j} = \emptyset$. 
Then $x_i < z < z' < y < x_{j}$. 
Since $(y,x_{j})_I \subset A_{i,j}$ and $(z,z')_I \cap A_{i,j} = \emptyset$, there is a point $y' \in [z',y]_I$ such that $(y',x_{j})_I \subset A_{i,j}$ and $\partial (y',x_{j})_I \subset O(x_{j})$, which contradicts the non-existence of such an open subarc $J := (y',x_{j})_I$ of $I$ due to Lemma~\ref{lem:inh_ori}. 

Thus, we have that $\partial A_{i,j} \cap (x_{i},x_{j})_I = \emptyset$. 
Since $\partial A_{i,j} = O(x_i) \sqcup O(x_j)$, we have $[x_{i},x_j)_I \cap O(x_{j}) = \emptyset$. 
By $[x_{i},x_\infty]_I \cap O(x_{i}) = \{ x_{i} \}$ and $[x_{i},x_\infty]_I \cap O(x_{j}) = \{ x_{j} \}$, we obtain the following equality:  
\[
\begin{split}
[x_{i},x_{\infty}]_I \cap \partial A_{i,j} &= [x_{i},x_{\infty}]_I \cap (O(x_{i}) \sqcup O(x_{j}))
= \{ x_i, x_{j} \}
\end{split}
\]
Since $(x_{i},x_{j})_I \cap A_{i,j} \neq \emptyset$, we have that $(x_{i},x_{j})_I = (x_{i},x_{\infty})_I \cap A_{i,j}$. 

For any $i' < j' \in \N$, since $(x_{i'},x_{j'})_I = (x_{i'},x_{\infty})_I \cap A_{i',j'}$, by the pairwise disjoint property of $(A_k)_{k \in \N}$, we have $\emptyset = A_{i',j'} \cap (A_{1,i'} \sqcup O(x_{i'})) \supseteq A_{i',j'} \cap (x_1,x_{i'}]_I $. 
For any $i' < j' \in \N$, we obtain that $(x_{i'},x_{j'})_I = (x_{i'},x_{\infty})_I \cap A_{i',j'} = ((x_1,x_{i'}]_I \sqcup (x_{i'},x_{\infty})_I) \cap A_{i',j'} = (x_1,x_{\infty})_I \cap A_{i',j'}$, which implies assertion (1). 

Fix any $j' \in \N$. 
From the pairwise disjoint property of $(\A_{i'})_{i' \in \N}$ of open invariant subsets, since $O(x_{j'}) \subset \partial \A_{j'}$, assertion (1) implies the following inclusion: 
\[
\begin{split}
O(x_{j'}) \cap ((x_1,x_\infty)_I - \{ x_{j'}\})  &= O(x_{j'}) \cap (\bigsqcup_{i' \in \N} (x_{i'},x_{i'+1})_I \sqcup \bigsqcup_{i' \neq j'} \{ x_{i'}\}) 
\\
& \subseteq O(x_{j'}) \cap (\bigsqcup_{i' \in \N} \A_{i'} \sqcup \bigsqcup_{i' \neq j'} O(x_{i'})) = \emptyset
\end{split}
\]
Therefore, we have $\{ x_{j'} \} = (x_1,x_{\infty}) \cap O(x_{j'})$. 
\end{proof}

Claim~\ref{claim:016-} implies the following equality:  
\[
\begin{split}
(x_1,x_{\infty})_I &= (x_1,x_2)_I \sqcup \bigsqcup_{l=2}^\infty \left( \{ x_{l} \} \sqcup (x_{l},x_{l+1})_I \right)
\\
&= ((x_1,x_\infty)_I \cap A_{1,2}) \sqcup \bigsqcup_{l=2}^\infty \left( \{ x_{l} \} \sqcup ((x_1,x_{\infty}) \cap A_{l,l+1}) \right)
\\
&= (x_1,x_\infty)_I \cap  \left( A_{1,2} \sqcup \bigsqcup_{l=2}^\infty \left( O(x_{l}) \sqcup A_{l,l+1}\right) \right)
= (x_1,x_\infty)_I \cap A
\end{split}
\]
Therefore, we have $(x_1,x_\infty)_I \subset A$. 

\begin{claim}\label{claim:016}
For any $j \in \N$, the connected component of $S - \bigcup_{k \in \mathbb{N}} O(x_k)$ which contains $(x_j,x_{j+1})_I$ is the open annulus $\A_j = A_{j,j+1}$. 
\end{claim}

\begin{proof}[Proof of Claim~\ref{claim:016}]
Since the sequence $(x_i)_{i \in \mathbb{Z}_{\geq 0}}$ is monotonically decreasing in $I$, by $(x_1,x_\infty)_I \subset A$, 
Claim~\ref{claim:016-} and the property (P3) imply the assertion. 
\end{proof}

From the finiteness of singular points, by taking a subsequence of $(x_i)_{i \in \N}$, we may assume that, for any $k \in \N$, the annulus $\A_k = A_{k,k+1}$ contains no singular points. 
By replacing $I$ with $(x_1,x_\infty)_I$ and taking a subsequence of $(x_i)_{i \in \N}$ if necessary, we may assume that the sequence $(x_i)_{i \in \mathbb{Z}_{\geq 0}}$ consists of fixed points of the first return map $f_I$ on $I$ generated by $v$. 
From Claim~\ref{claim:016}, applying Lemma~\ref{lem:circuit-} to $(x_i)_{i \in \mathbb{N}}$, the point $x_\infty$ is contained in a circuit in $\overline{\Pv}$. 
\end{proof}

\subsection{Equivalence for the finite case}

In the remainder of this section, we prove Theorem~\ref{thm:finite_dense}. 

\begin{lemma}\label{lem:bdry_circuit}
Let $v$ be a flow with finitely many singular points on a compact surface $S$, $x \in S$ a non-recurrent point,  $I$ an open transverse arc containing $x$, and $f_I$ the first return map on $I$. 
Suppose that $f_I$ is orientation-reversing and $f_I(I_-) = I_+$, where $I_-, I_+$ are the connected components of $I - \{ x \}$. 
Denote by $D$ the union of open orbit arcs from points in $I_-$ to the image of $f_I$. 
The connected component of $\partial D - (I_- \sqcup I_+)$ containing $x$ is a circuit. 
\end{lemma}

\begin{proof}
By definition, the union $D$ is an open flow box as on the left in Figure~\ref{nondir02}.
\begin{figure}
\begin{center}
\includegraphics[scale=0.15]{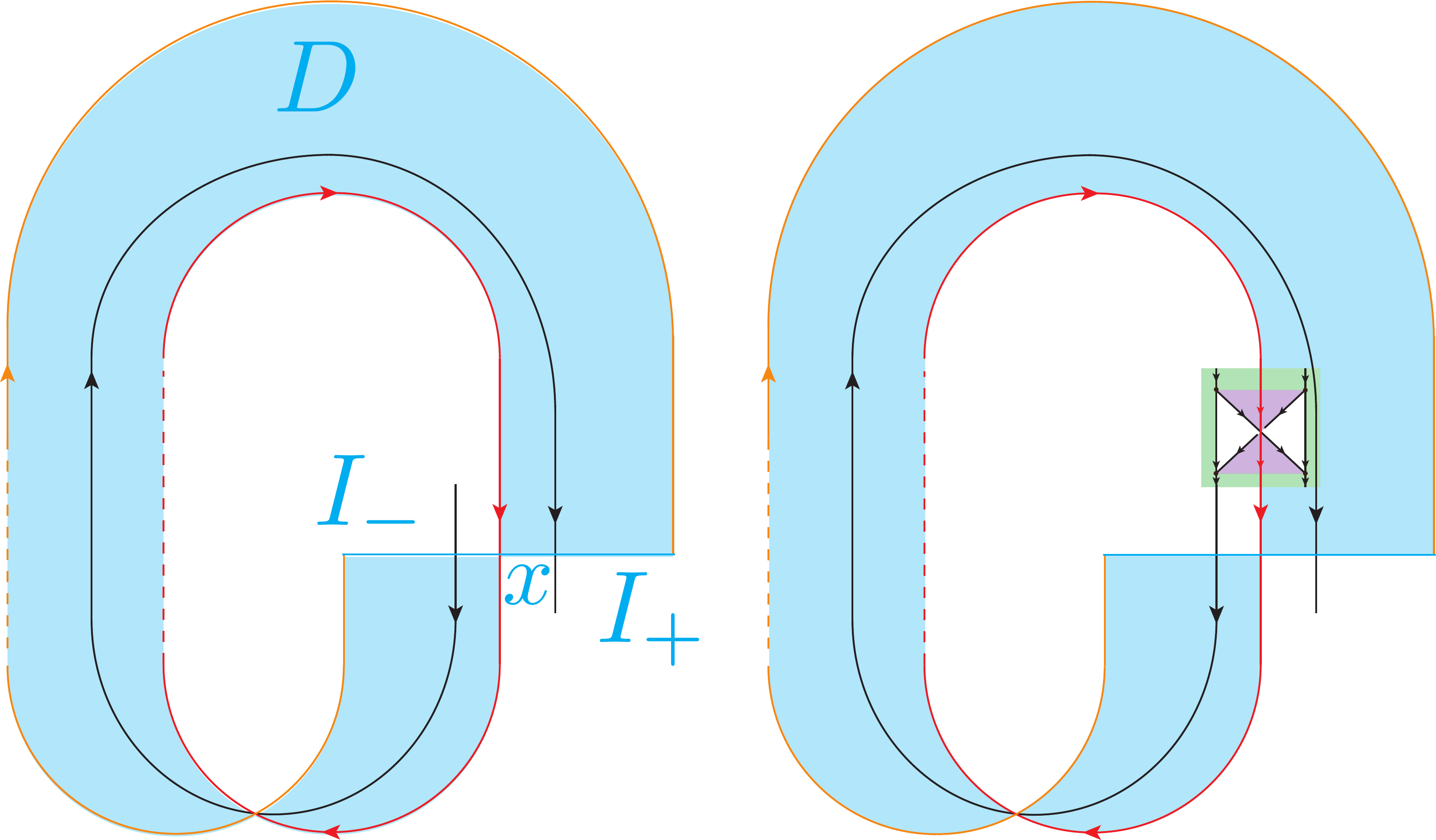}
\end{center}
\caption{Left, a trivial open flow box; right, the modified flow by replacing a trivial open flow box with a one-punctured M{\"o}bius band as in Figure~\ref{twist01}.}
\label{nondir02}
\end{figure} 
Let $\gamma$ be the connected component of $\partial D - (I_- \sqcup I_+)$ containing $x$. 
Replace a small trivial open flow box containing a point in $O^-(x)$ with a one-punctured M{\"o}bius band as in Figure~\ref{nondir02} and Figure~\ref{twist01}. 
We have that the resulting surface is compact and that the resulting flow $w$ has finitely many singular points. 
\begin{figure}
\begin{center}
\includegraphics[scale=0.2]{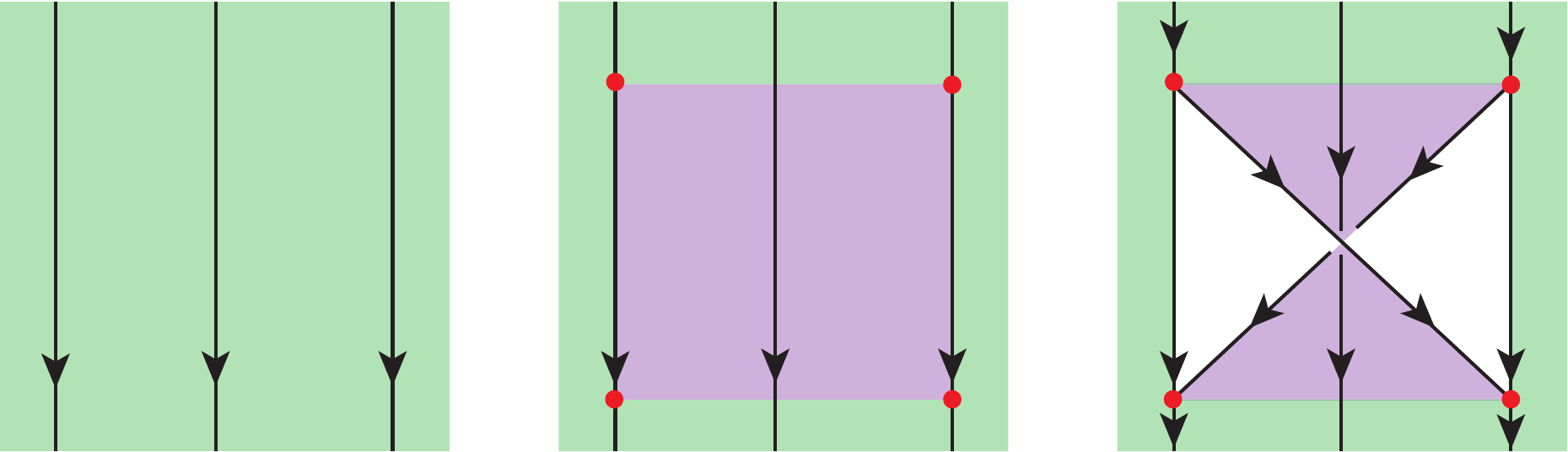}
\end{center}
\caption{Left, a trivial open flow box $D$; middle, an open flow box with exactly four singular points; right, a one-punctured M{\"o}bius band obtained from the open flow box with exactly four singular points by cutting an open disk, twisting the disk, and pasting two disjoint closed intervals.}
\label{twist01}
\end{figure} 
By deforming the M{\"o}bius band isotopically in a direction transverse to the flow $w$ near $I_+$ to slightly modify the flow inside the band, we may assume that the first return map $f_{w,I}$ on $I_-$ by $w$ is a continuous injective mapping $f_{w,I} \colon I_- \to I_-$ which is attracting. 
Then $\gamma$ is the $\omega$-limit set of a point in $I_-$ near $x$ with respect to $w$ and is a boundary component of an open annulus contained in $w(I_-)$. 
By Lemma~\ref{lem:lc_transversal}, the $\omega$-limit set $\gamma$ with respect to $w$ is a semi-attracting limit circuit with respect to $w$. 
Since the operation described in Figures~\ref{nondir02} and Figure~\ref{twist01} preserves $\gamma$, the subset $\gamma$ is also a circuit with respect to $v$.
\end{proof}

To show the existence of either a strictly limit circuit or a circuit with wandering holonomy, we have the following observation by the same argument of the proof of Lemma~\ref{lem:inh_ori}. 

\begin{figure}[t]
\begin{center}
\includegraphics[scale=0.5]{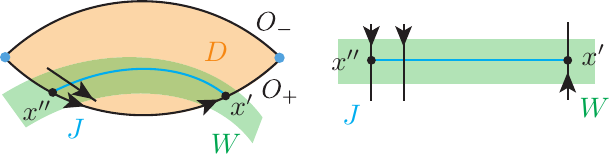}
\end{center}
\caption{The disk $A$ and the open arc $J \subset D$ with $\partial J \subset O_+$.}
\label{fig:trans_orientation02}
\end{figure}

\begin{lemma}\label{lem:inh_ori_02}
Let $v$ be a flow on a surface $S$ and $D \subseteq S$ an invariant closed disk whose boundary consists of two singular points $\alpha$ and $\omega$ and non-singular orbits $O_-$ and $O_+$ with $\alpha = \alpha(O_-) = \alpha(O_+)$ and $\omega = \omega(O_-) = \omega(O_+)$. 
Then there is no closed transverse arc $J \subset D$ with either $\partial J \subset \overline{O_-}$ or $\partial J \subset \overline{O_+}$. 
\end{lemma}
 
\begin{proof}
Assume that there is a closed transverse arc $J \subset D$ with $\partial J \subset \overline{O_+}$ as on the left in Figure~\ref{fig:trans_orientation02}. 
From the transversality of $J$, we have that $J \cap \{ \alpha, \omega \} = \emptyset$ and so that $\partial J \subset O_+$. 
By the flow box theorem, there is an open flow box $W$ containing $J$ as in Figure~\ref{fig:trans_orientation02}. 
Denote by $x'_i$ and $x''_i$ the points in $\partial J$. 
The flow $v$ defines inherent directions on $O_+$ at $x''$.
However, the trivial flow box $W$ containing $J$ implies a direction at $x'$ that is opposite to the direction induced by following the orbit arc from $x''$ to $x'$, which is a contradiction.
Therefore, there is no open arc $J \subset A$ with $\partial J \subset O_+$. 
By the same argument for $O_-$, the assertion holds. 
\end{proof}

We have the following statement.

\begin{lemma}\label{lem:finite_ends}
Let $v$ be a flow on a closed surface $S$. 
Any connected component $S_0$ of $S - \overline{\mathrm{R}(v)}$ is an open surface which is homeomorphic to the resulting surface from a closed surface $T$ by removing finitely many points in $T$. 
\end{lemma}
 
\begin{proof}
Fix any connected component $S_0$ of $S - \overline{\mathrm{R}(v)}$. 
Since the complement $S - \overline{\mathrm{R}(v)}$ is an open surface, the connected component $S_0$ is also an open surface. 
From \cite[Theorem~3]{richards1963classification}, the surface $S_0$ is homeomorphic to the resulting surface from a closed surface $T$ by removing a closed totally disconnected subset $T'$. 
Therefore, we identify $S - \overline{\mathrm{R}(v)}$ with $T - T'$.  
Fix any Riemannian metric on $T$. 

Assume that $T'$ is infinite. 
There is an accumulation point $x \in \overline{T' \setminus \{ x \}}^T$, where $\overline{A}^{T}$ is the closure in $T$ of a subset $A \subseteq T$. 
Then there is a sequence $(x_n)_{n \in \mathbb{N}}$ of points $x_n \in T' \setminus \{ x \}$ converging to $x$. 
Therefore, there are a sequence $(D_n)_{n \in \mathbb{N}}$ of pairwise disjoint closed disks $D_n \subset T$ each of whose boundaries $\gamma_n := \partial D_n$ is a loop contained in the complement $T - T' = S - \overline{\mathrm{R}(v)}$ and such that $x_n \in D_n \cap T'$. 
Since $D_n \cap D_m = \emptyset$ for any $n \neq m \in \mathbb{N}$, we have that $T' \setminus D_k \neq \emptyset$ for any $k \in \mathbb{N}$. 

\begin{claim}\label{claim:gamma_noncont}
For any $k \in \mathbb{N}$, the loop $\gamma_k$ is homotopically non-trivial in $S$. 
\end{claim}

\begin{proof}[Proof of Claim~\ref{claim:gamma_noncont}]
Assume that there is a loop $\gamma_k$ which is homotopically trivial in $S$. 
Let $D$ be a closed disk in $S$ whose boundary is $\gamma_k$. 
By $\overline{\mathrm{R}(v)} \cap (\bigsqcup_{n \in \mathbb{N}} \gamma_n) = \emptyset$, since any closed disk contains no $Q$-sets, the disk $D$ does not intersect $\overline{\mathrm{R}(v)}$. 
Then $D \subset S - \overline{\mathrm{R}(v)} = T - T'$. 
Since the disk $D_k$ intersects $T'$, we have that $D \neq D_k$ and so that $D \cap D_k = \gamma_k$. 
Then $T = D \cup D_k$ is a sphere such that $T' \subset T - D \subset D_k$, which contradicts $T' \not\subseteq D_k$. 
\end{proof}

By the finiteness of the genus of the closed surface $S$, by taking a subsequence of $(D_n)_{n \in \mathbb{N}}$, we may assume that the loops $\gamma_n$ and $\gamma_m$ for any $n \neq m \in \mathbb{N}$ are freely homotopic in $S$. 
By Claim~\ref{claim:gamma_noncont}, the loop $\gamma_n$ for any $n \in \mathbb{N}$ is homotopically non-trivial. 
Therefore, for any $n \neq m \in \mathbb{N}$, there is a closed annulus whose boundary is $\gamma_n \sqcup \gamma_m$. 
By renumbering $(D_n)_{n \in \mathbb{N}}$, we can take a closed annulus $A$ whose boundary is $\gamma_1 \sqcup \gamma_3$ such that $\gamma_2 \subset \mathop{\mathrm{int}}A$. 
From $\overline{\mathrm{R}(v)} \cap (\bigsqcup_{n \in \mathbb{N}} \gamma_n) = \emptyset$, we have that $\overline{\mathrm{R}(v)} \cap \partial A = \emptyset$.
Lemma~\ref{lem:annulus_nonrec} implies that $\mathrm{R}(v) \cap \mathop{\mathrm{int}}A = \emptyset$ and so that $\overline{\mathrm{R}(v)} \cap \mathop{\mathrm{int}}A = \emptyset$. 
Then $\overline{\mathrm{R}(v)} \cap A = \emptyset$. 
This means that $A \subset S - \overline{\mathrm{R}(v)} = T - T'$. 
On the other hand, since $\partial D_2 = \gamma_2 \subset A$ and $D_2 \cap \partial A= D_2 \cap (\gamma_1 \sqcup \gamma_3) = \emptyset$, we obtain that $D_2 \subset A \subset T - T'$, which contradicts $D_2 \cap T' \neq \emptyset$. 
\end{proof}

We state the existence of either a strictly limit circuit or a circuit with wandering holonomy. 

\begin{lemma}\label{lem:wandering_h}
Let $v$ be a flow with finitely many singular points without non-closed recurrent points on a compact surface $S$. 
For any point $x$ in $\Omega(v) - \overline{\Cv}$, there is either a strictly limit circuit at $O(x)$ or a circuit with wandering holonomy at $O(x)$.
\end{lemma}

\begin{proof}
Since each recurrent orbit is closed, we have that $\mathrm{R}(v) = \emptyset$. 
Let $d$ be a Riemannian distance on $S$. 
Fix any point $x$ in $\Omega(v) - \overline{\Cv}$. 

\begin{claim}\label{claim:lc_only}
For any non-closed point $z$ with $x \in \omega(z)$ {\rm(resp.} $x \in \alpha(z)${\rm)}, the limit set $\omega(z)$ {\rm(resp.} $\alpha(z)${\rm)} is a limit circuit containing $x$. 
\end{claim}

\begin{proof}[Proof of Claim~\ref{claim:lc_only}]
By Lemma~\ref{lem:lc_transversal}, from the non-existence of Q-sets, the $\omega$-limit or $\alpha$-limit set of a non-closed point is either a singular point, a limit cycle, or a non-periodic limit circuit. 
Since $x$ is non-singular, the $\omega$-limit or $\alpha$-limit set of any non-closed point which contains $x$ is a limit circuit. 
\end{proof}

If there is a limit circuit containing $x \in \Omega(v) - \overline{\Cv}$, then Lemma~\ref{lem:lim_cycle_strict} implies the circuit is a strictly limit circuit at $O(x)$. 
Thus, we may assume that $x$ is not contained in a limit circuit. 
Since $O(x)$ is not a closed orbit, by Claim~\ref{claim:lc_only}, the point $x$ is not contained in the $\omega$-limit or $\alpha$-limit set of any point. 
Since any limit circuit of $v$ is a closed subset which does not contain $x$, by $x \notin \overline{\Cv}$, the finiteness of non-periodic limit circuits due to Theorem~\ref{main:fininite_existence_limit_circuits} implies that the union of limit circuits and $\overline{\Cv}$ does not contain $x$. 
%

\begin{claim}\label{claim:closed}
It suffices to prove the lemma in the case where $S$ is a closed surface.
\end{claim}

\begin{proof}[Proof of Claim~\ref{claim:closed}]
By taking the lift $v_{\mathrm{db}}$ of the flow $v$ to the double $S_{\mathrm{db}}$ of $S$, the subset $S \subset S_{\mathrm{db}}$ is an invariant subset of $v_{\mathrm{db}}$ and so the property that $x$ is contained in a circuit with wandering holonomy at $O_v(x) = O_{v_{\mathrm{db}}}(x)$ is invariant under taking $v_{\mathrm{db}}$. 
Thus, by taking $v_{\mathrm{db}}$ if necessary, we may assume that the surface $S$ is a closed surface. 
\end{proof}

By the previous claim, there is an open trivial flow box $B$ containing $x$ which does not intersect the union of limit circuits and $\overline{\Cv}$. 
Since $x$ is non-recurrent, by taking $B$ small if necessary, we may assume that the intersection $O(x) \cap B$ is an open interval. 

Let $T \subset B$ be a small closed transverse arc whose interior intersects $x$, and $f_{T}$ the first return map on $T$. 
Then $T$ does not intersect the union of limit circuits and $\overline{\Cv}$. 
Because the intersection $O(x) \cap B$ is an open interval, we have that $x \notin v(T - \{ x\})$. 
Since $x \in \Omega(v) - \overline{\Cv}$, there are a connected component $T_+$ of $T - \{ x \}$, a sequence $(a_n)_{n \in \mathbb{N}}$ of non-closed points in $T_+ \cap \mathrm{dom}(f_T)$, and a sequence $(N_n)_{n \in \mathbb{N}}$ of positive integers such that the sequence $(a_n)_{n \in \mathbb{N}}$ is totally ordered on $T$ and monotonically converges to $x$ on $T$, $d(x, a_n) < 1/n$, and $d(x, f_{T}^{N_{n}}(a_n)) < 1/n$ for any $n \in \mathbb{N}$. 
Then $(x,a_{n+1})_T \subset (x,a_n)_T$ for any $n \in \mathbb{N}$.
By shortening $T_+$ if necessary, we may assume that $a_1$ is a boundary component of $T$. 
We define a natural total order on the interval $T$ such that $a_1$ is the maximal element. 
With respect to this total order, the sequence $(a_n)_{n \in \mathbb{N}}$ converges monotonically to $x$ in $T$.
Since $x \notin v(T - \{ x\})$, we obtain that $x \notin \bigcup_{n \in \mathbb{N}} O(a_n)$. 
Denote by $a_-$ the minimal element of $T$. 
Then $T_+ = (x,a_1]_T$ and $T = [a_-,a_1]_T$. 

\begin{claim}\label{claim:00}
It suffices to prove the lemma in the case where $O(a_n) \neq O(a_m)$ for any $n \neq m \in \mathbb{N}$.
\end{claim}

\begin{proof}[Proof of Claim~\ref{claim:00}]
Assume that there is an orbit $O$ containing infinitely many points in $(a_n)_{n \in \mathbb{N}}$. 
Then there is a subsequence $(a_{k_n})_{n \in \mathbb{N}}$ of $(a_n)_{n \in \mathbb{N}}$ such that either $a_{k_{n+1}} \in O^+(a_{k_n})$ for any $n \in \mathbb{N}$ or $a_{k_{n+1}} \in O^-(a_{k_n})$ for any $n \in \mathbb{N}$. 
Therefore, we have that $x \in \alpha(O) \cup \omega(O)$, which contradicts the fact that $x$ is not contained in the $\omega$-limit or $\alpha$-limit set of any orbit. 

Thus, any orbit contains at most finitely many points in $(a_n)_{n \in \mathbb{N}}$.
By taking a subsequence of $(a_n)_{n \in \mathbb{N}}$ if necessary, we may assume that $O(a_n) \neq O(a_m)$ for any $n \neq m \in  \mathbb{N}$.
\end{proof}

\begin{claim}\label{claim:01}
There are at most finitely many limit cycles each of which is the $\alpha$-limit or the $\omega$-limit set of some point in $(a_n)_{n \in \mathbb{N}}$. 
\end{claim}

\begin{proof}[Proof of Claim~\ref{claim:01}]
Assume that there are infinitely many pairwise disjoint limit cycles $C_{n_k}$ each of which is the $\omega$-limit set of $a_{n_k}$ in some subsequence $(a_{n_k})_{k \in \N}$ of $(a_{n})_{n \in \N}$. 
By taking a subsequence of $(a_{n})_{n \in \N}$, we may assume that $n_k = k$ for any $k \in \N$. 
Since $x \notin \overline{\Cv}$ and $\bigsqcup_n C_n \subseteq \Cv$, we have $x \notin \overline{\bigsqcup_n C_n}$. 
From the compactness of $S$, by taking a subsequence of $(C_n)_{n \in \N}$ if necessary, we may assume that the periodic orbits in $(C_n)_{n \in \N}$ are pairwise freely homotopic in $S$. 
By taking a subsequence of $(C_n)_{n \in \N}$, from Lemma~\ref{claim:003a}, we may assume that there is a sequence $(\A_{n})_{n \in \N}$ of pairwise disjoint open annuli $\A_n$ with $\partial \A_n = C_n \sqcup C_{n+1}$ satisfying the following properties:
\begin{enumerate}[\rm (P1)]
\item The union $\A := \bigcup_{k=1}^\infty (\A_{n} \sqcup C_{n+1})$ is an open annulus one of whose boundary components is the loop $C_{1}$. 
\item For any $n \in \N$, the annulus $\A_{n}$ and the loop $C_{n+1}$ are parallel to the boundary of $\A$.
\item The set of connected components of $S - \bigcup_{n=1}^\infty C_{n}$ contained in $\A$ coincides with the family $\{ \A_{n}\}_{n \in \N}$ of annuli. 
\item For any $n<n' \in \N$, the union $\A_{n} \sqcup \bigsqcup_{l = n+1}^{n'-1} (C_l \sqcup \A_{l}) \subset \A$ is an open annulus whose boundary is $C_n \sqcup C_{n'}$.
\end{enumerate}

We will show that $x \notin \A$. 
Indeed, assume that $x \in \A$. 
Since $x \notin \overline{\mathop{\mathrm{Cl}}(v)}$, we have that $x \in \A - \bigcup_{k=1}^\infty  C_{n+1} = \bigcup_{k=1}^\infty \A_{n}$ and so that there is an integer $N \in \N$ such that $x \in \A_N \subset \overline{\A_N} = \A_N \sqcup C_N \sqcup C_{N+1}$. 
By $\lim_{n \to \infty} d(x, a_n) = 0$, the \nbd $\A_N$ of $x$ contains all $(a_n)_{n \in \N_{> L}}$ for some $L >N$. 
Fix any $L' \in \N_{> L}$. 
Then $a_{L'} \in \A_N$ and so $C_{L'} = \omega(a_{L'}) \subset \overline{\A_N} = \A_N \sqcup C_N \sqcup C_{N+1}$. 
Since $(C_n)_{n \in \N}$ is pairwise disjoint, we obtain that $C_{L'} \subset \A_N \subset  S - \bigcup_{n=0}^\infty C_{n} \subset S - C_{L'}$, which is a contradiction.  

Let $\A_x$ be the connected component of the complement $S - \overline{\bigsqcup_{n \in \N} C_n}$ containing $x$. 
By $\lim_{n \to \infty} d(x, a_n) = 0$, the \nbd $\A_x$ of $x$ contains all $(a_n)_{n \in \N_{> L}}$ for some $L >0$. 
On the other hand,  by (P4), the interior $\A_{L+1} \sqcup C_{L+2} \sqcup \A_{L+2}$ of the closed invariant annulus $C_{L+1} \sqcup \A_{L+1} \sqcup C_{L+2} \sqcup \A_{L+2} \sqcup C_{L+3}$ contains the limit cycle $C_{L+2} = \omega(a_{L+2})$ and so contains $a_{L+2} \in \A_x \cap (\A_{L+1} \sqcup \A_{L+2})$. 
Since $\A_x$ is a connected component of $S - \overline{\bigsqcup_{m>0} C_m} \subset S - \bigsqcup_{m \geq L} C_m$, we obtain either $x \in \A_x \subset \A_{L+1}$ or $x \in \A_x \subset \A_{L+2}$, which contradicts $x \notin \A_{L+1} \sqcup \A_{L+2}$. 

By symmetry of the time reversion, there are at most finitely many limit cycles each of which is the $\alpha$-limit set of some point in $(a_n)_{n \in \mathbb{N}}$. 
\end{proof}

\begin{claim}\label{claim:02}
The set of the $\alpha$-limit sets and the $\omega$-limit sets of points in $(a_n)_{n \in \mathbb{N}}$ is finite. 
\end{claim}
\begin{proof}[Proof of Claim~\ref{claim:02}]
Assume that there are infinitely many pairwise distinct $\omega$-limit sets of points $a_n$. 
From Claim~\ref{claim:01}, there are at most finitely many limit cycles. 
From the non-existence of Q-sets, by taking a subsequence of $(a_n)$ if necessary, we may assume that each of the $\omega$-limit and $\alpha$-limit sets of any point $a_n$ is neither a limit cycle nor a Q-set. 
By Lemma~\ref{lem:lc_transversal}, the $\omega$-limit and $\alpha$-limit sets of any point $a_n$ are either singular points or non-periodic limit circuits. 
By the finiteness of singular points, there are infinitely many pairwise distinct non-periodic limit circuits which are the $\omega$-limit sets of some points in $(a_n)_{n \in \mathbb{N}}$, which contradicts the existence of finitely many non-periodic limit circuits due to Theorem~\ref{main:fininite_existence_limit_circuits}. 
By symmetry of the time reversion, the claim holds. 
\end{proof}

By taking a subsequence of $(a_n)_{n \in \mathbb{N}}$ if necessary, we may assume that there are closed invariant subsets $\alpha$ and $\omega$ in $S - T$ such that $\alpha = \alpha(a_n)$ and $\omega = \omega(a_n)$ for any $n  \in \mathbb{N}$. 
From the non-existence of non-closed recurrent points of $v$, each of $\alpha$ and $\omega$ is either a singular point or a limit circuit.

\begin{claim}\label{claim:reduction_sing}
It suffices to prove the lemma in the case where each of $\alpha$ and $\omega$ is a singular point.
\end{claim}

\begin{proof}[Proof of Claim~\ref{claim:reduction_sing}]
Since $x$ is not contained in any limit circuit, cutting limit circuits and collapsing the new boundary components into singletons if necessary, we may assume that each of $\alpha$ and $\omega$ is not a limit circuit, because the property that a circuit has wandering holonomy at $O(x)$ is invariant under the cutting and collapsing operations. 
\end{proof}

%
%

By taking a subsequence of $(a_n)_{n \in \mathbb{N}}$ if necessary, we may assume that the closed intervals $\overline{O(a_n)}$ and $\overline{O(a_m)}$ are homotopic relative to $\{\alpha, \omega \}$ for any $n,m \in \mathbb{N}$. 
Then $x \notin \overline{O(a_n)} = O(a_n) \sqcup \{\alpha, \omega \}$ for any $n \in \mathbb{N}$. 
From Claim~\ref{claim:00}, we have that $\overline{O(a_n)} \neq \overline{O(a_m)}$ for any $n,m \in \mathbb{N}$. 
Therefore, for any $n,m \in \mathbb{N}$, there is an open disk whose boundary is a loop $O(a_n) \sqcup O(a_{n+1}) \sqcup \{\alpha, \omega \}$. 
Since the sequence $(a_n)_{n \in \Z_{\geq 0}}$ monotonically converges to $x$ on $T$ and any closed interval $\overline{O(a_n)}$ intersects the closed transverse arc $T$ at most finitely many times, by replacing $a_n$ with a point in $O(a_n) \cap T$, 
and 
taking a subsequence of $(a_n)_{n \in \mathbb{N}}$, 
we may assume that any $a_n$ is minimal in the points in $O(a_n) \cap T$ with respect to the total order on $T$. 
Then $[a_n,x]_T \cap \overline{O(a_n)} = \{ a_n \}$. 
Since $S$ has finite genus, by taking a subsequence of $(a_n)_{n \in \mathbb{N}}$, we may assume that any connected component of $S - \left(\{\alpha, \omega \} \sqcup \bigsqcup_{j \in \mathbb{N}} O(a_j)\right)$ which intersects $T_+ = (x,a_1]_T$ and whose boundary is $O(a_i) \sqcup O(a_{i'}) \sqcup \{\alpha, \omega \}$ for some $i < i' \in \mathbb{N}$ is an open disk. 

\begin{claim}\label{claim:021}
It suffices to prove the lemma in the case where, for any $n \in \mathbb{N}$, there is an invariant open disk $U_n$ whose boundary is $O(a_n) \sqcup O(a_{n+1}) \sqcup \{\alpha, \omega \}$ such that $U_n$ is the connected component of $S - \left(\{\alpha, \omega \} \sqcup \bigsqcup_{j \in \N} O(a_j)\right)$ which contains $(a_{n+1},a_n)_T$ and that $x \notin U_n$. 
\end{claim}
\begin{proof}[Proof of Claim~\ref{claim:021}]
Suppose that there are positive numbers $M < N \in \mathbb{N}$ and an open disk $U_{M,N}$ whose boundary is a loop $O(a_M) \sqcup O(a_N) \sqcup \{\alpha, \omega \}$ such that $x \in U_{M,N}$. 
Since $(a_n)_{n \in \mathbb{N}}$ converges to $x$, by taking a subsequence of $(a_n)_{n \in \mathbb{N}}$, we may assume $a_n \in U_{M,N}$ for any $n> N \in \N$. 
By the invariance of $U_{M,N}$, we have $O(a_n) \subset U_{M,N}$ for any $n > N \in \N$.
Since the closed disk $\overline{U_{M,N}} = U_{M,N} \sqcup O(a_M) \sqcup O(a_N) \sqcup \{\alpha, \omega \}$ is homotopically trivial, the closed arc $O(a_n) \sqcup \{\alpha, \omega \}$ is freely homotopic to the closed arc $O(a_N) \sqcup \{\alpha, \omega \}$ in the closed disk $\overline{U_{M,N}}$ for any $n > N \in \N$.  
By taking a subsequence of $(a_n)_{n \in \mathbb{N}}$, we may assume that either $U_{M,n} \subset U_{M,n+1}$ for any $n \in \N_{> N}$, or $U_{N,n} \subset U_{N,n+1}$ for any $n \in \N_{> N}$, where $U_{M,m}$ (resp. $U_{N,m}$) is the open disk in $U_{M,N}$ whose boundary is $O(a_M) \sqcup O(a_m) \sqcup \{\alpha, \omega \}$ (resp. $O(a_N) \sqcup O(a_m) \sqcup \{\alpha, \omega \}$) for any $m \in \N_{> N}$.
Therefore, there is the connected component $U_n$ of $S - \left(\{\alpha, \omega \} \sqcup \bigsqcup_{j \in \Z_{>N}} O(a_j)\right)$ which contains $(a_{n+1},a_n)_T$ and is an invariant open disk. 
By taking a subsequence of $(a_n)_{n \in \mathbb{N}}$, we may assume that the claim holds. 

Thus, we may assume that, for any positive numbers $i < j \in \mathbb{N}$, there is no open disk which contains $x$ and whose boundary is a loop $O(a_i) \sqcup O(a_j) \sqcup \{\alpha, \omega \}$. 
For any $i < j \in \mathbb{N}$, let $U_{i,j}$ be the invariant open disk whose boundary is a loop $O(a_i) \sqcup O(a_j) \sqcup \{\alpha, \omega \}$ with $x \notin U_{i,j}$. 
By Lemma~\ref{lem:inh_ori_02}, there is no open arc $J \subset T \cap U_{i,j}$ with either $\partial J \subset O(a_j)$ or $\partial J \subset O(a_i)$. 
Since $[a_{i},x]_T \cap \overline{O(a_{i})} = \{ a_{i} \}$ and $x \notin U_{i,j}$ for any $i<j \in \mathbb{N}$, this non-existence of such open arcs implies that $(a_{i},x)_T \cap \overline{O(a_j)} = \{ a_j \}$ and that the arc $(a_{i},x)_T$ intersects $U_{i,j}$ as shown on the left in Figure~\ref{fig:rot_arc}. 
\begin{figure}
\begin{center}
\includegraphics[scale=0.7]{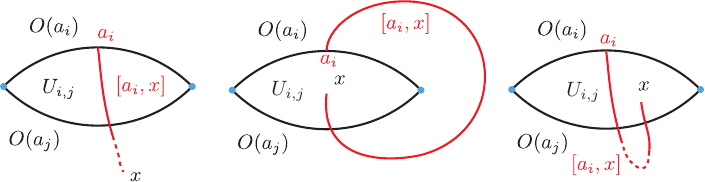}
\end{center}
\caption{Left: A possible pattern. Middle and Right: Impossible patterns.}
\label{fig:rot_arc}
\end{figure} 
In particular, any small arc from $a_i$ contained in $(a_{i},x)_T$ intersects $U_{i,j}$. 
Then $[a_{i},x]_T \cap \partial U_{i,j} = [a_{i},x]_T \cap (\overline{O(a_{i})
} \cup \overline{O(a_{j})}) = \{ a_{i}, a_j \}$ and so $(a_i,a_{j})_T \subset U_{i,j}$. 
Therefore, the connected component $U_i := U_{i,i+1}$ of $U - \left(\{\alpha, \omega \} \sqcup \bigsqcup_{j \in \mathbb{N}} O(a_j)\right)$ has the desired properties. 
\end{proof}


The previous claim implies that $x \notin \overline{U_n}$ for any $n \in \mathbb{N}$. 
For any $n \in \mathbb{N}$, denote by $U_{n-}$ (resp. $U_{n+}$) the connected component of $\overline{U_n} - [a_{n+1},a_n]_T$ containing $\alpha$ (resp. $\omega$).  

\begin{claim}\label{claim:022}
It suffices to prove the lemma in the case where, for any $n \in \mathbb{N}$, we have that $\{ \alpha \} = \bigcup_{a \in [a_n,a_{n+1}]_T} \alpha(a)$ and $\{ \omega \} = \bigcup_{a \in [a_n,a_{n+1}]_T} \omega(a)$. 
\end{claim}

\begin{proof}[Proof of Claim~\ref{claim:022}]
By the finiteness of singular points, by taking a subsequence of $(a_n)_{n \in \mathbb{N}}$, we may assume that, for any $n \in \mathbb{N}$, the disk $U_n$ contains no singular points. 
Since the Poincar{\'e}-Hopf theorem implies that any periodic orbit in an open disk bounds a singular point, by taking a subsequence of $(a_n)_{n \in \mathbb{N}}$, we may assume that, for any $n \in \mathbb{N}$, the disk $U_n$ contains no periodic points. 
Fix any $n \in \mathbb{N}$. 
Then $\overline{U_{n}} \cap \Cv = \{ \alpha, \omega \}$. 
Moreover, the loop $\partial \overline{U_{n}} = \{ \alpha, \omega \} \sqcup O_n \sqcup O_{n+1}$ is a directed graph but is not a non-trivial circuit because the orientation is not compatible with the dynamical orientation of the orbits.

Assume that there is a limit circuit $\gamma \in \overline{U_{n}}$. 
By time reversion if necessary, we may assume that $\gamma$ is the $\omega$-limit set of a point $p \in S$. 
Lemma~\ref{lem:lc_transversal} implies that $\gamma$ is semi-attracting. 
Since the loop $\partial \overline{U_{n}}$ is not a limit circuit, there is a non-singular orbit $\nu \subset \gamma$. 
Because $\gamma$ is a semi-attracting limit circuit, there is an open transverse arc $I \subset U_n$ whose boundary contains a point $\nu$, and there is a point in $I$ whose orbit intersects $I$ infinitely many times. 
The invariance of $U_n$ implies that $v(I) \subset U_n$. 
By \cite[Lemma~3.2]{yokoyama2021poincare}, there is a closed transversal in the open disk $U_n$. 
Then there is a disk $B_n \subset U_n$ whose boundary is the closed transversal. 
The Poincar{\'e}-Hopf theorem implies $B_n \subset U_n$ contains a singular point, which contradicts the non-existence of singular points in $U_n$.  
Thus, the $\alpha$-limit or the $\omega$-limit set of any point $a \in U_n$ is a singular point in $\overline{U_{n}} \cap \Cv = \{ \alpha, \omega \}$. 

Since the set difference $\overline{U_{n}} - [a_n,a_{n+1}]_T$ consists of the positive invariant disk $U_{n+}$ and the negative invariant disk $U_{n-}$ with $U_{n+} \cap \Cv = \{ \omega \}$ and $U_{n-} \cap \Cv = \{ \alpha \}$, 
This implies that $\alpha = \alpha(a)$ and $\omega = \omega(a)$ for any $a \in [a_n,a_{n+1}]_T$. 
\end{proof}

\begin{claim}\label{claim:filling}
It suffices to prove the lemma in the case where, for any $n \in \mathbb{N}$, we have that 
$U_n \sqcup O(a_n) \sqcup O(a_{n+1}) = \bigcup_{t \in \R} v(t,[a_{n+1},a_n]_T)$. 
\end{claim}

\begin{proof}[Proof of Claim~\ref{claim:filling}]
Fix any $n \in \mathbb{N}$. 
Since the closed disk $\overline{U_n}$ does not contain $x$, it suffices to show that, by collapsing $\overline{U_{n-}} \setminus v([a_{n+1},a_n]_T)$ (resp. $\overline{U_{n+}} \setminus v([a_{n+1},a_n]_T)$) into a singleton, we may assume that the conditions $O(a_n) \sqcup O(a_{n+1}) \subset \overline{U_n}$ are preserved by the collapse and we can identify the resulting singular points with $\alpha$ and $\omega$.
Indeed, consider the lift $w$ of the restriction $v\vert_{\overline{U_n}}$ to the double $S_n$ of the closed disk $\overline{U_n}$. 
Then $S_n$ is a sphere and $w$ consists of two singular points and non-recurrent points. 
Then the lift $T'$ of $[a_n,a_{n+1}]_T$ is a closed transversal of $w$ such that $w(T')$ is an open annulus with $\alpha = \alpha_w (a)$ and $\omega = \omega_w (a)$ for any $a \in T'$. 
Denote by $D_+$ (resp. $D_-$) the connected component $S_n - T'$ containing $\omega$ (resp. $\alpha$). 
Collapsing $D_+ \setminus w(T')$ (resp. $D_- \setminus w(T')$) into a singleton $\omega'$ (resp. $\alpha'$), the resulting flow $w'$ on the resulting sphere $S'_n$ satisfies that $\omega' = \omega_{w'} (a)$ and $\alpha' = \alpha_{w'}(a)$ for any $a \in S'_n - \{ \alpha', \omega' \} = w'(T')$. 
Then the collapse preserves $O(a_n) \sqcup O(a_{n+1}) \subset S_n$.
By construction, through the collapse of $\overline{U_{n-}} \setminus v([a_{n+1},a_n]_T)$ (resp. $\overline{U_{n+}} \setminus v([a_{n+1},a_n]_T)$) into a singleton, the orbits $O(a_n) \sqcup O(a_{n+1}) \subset \overline{U_n}$ are preserved and the resulting singular points can be identified with $\alpha$ and $\omega$.
\end{proof}

\begin{claim}\label{claim:saturation}
$U_n \sqcup O(a_n) \sqcup O(a_{n+1}) = \bigsqcup_{t \in \R} v(t,[a_{n+1},a_n]_T) = \bigsqcup_{a \in [a_{n+1},a_n]_T} O(a)$ for any $n \in \mathbb{N}$. 
\end{claim}

\begin{proof}[Proof of Claim~\ref{claim:saturation}]
Fix any $n \in \mathbb{N}$. 
For any $a \neq a' \in [a_{n+1},a_{n}]_T$, by the negative invariance of $U_{n-}$ and the positive invariance of $U_{n+}$, we have that $a' \notin O(a)$ and so that $O(a) \cap O(a') = \emptyset$. 
This implies that $U_n \sqcup O(a_n) \sqcup O(a_{n+1}) = \bigsqcup_{t \in \R} v(t,[a_{n+1},a_n]_T) = \bigsqcup_{a \in [a_{n+1},a_n]_T} O(a)$. 
\end{proof}

Put $U_\infty := \bigsqcup_{n \in \mathbb{N}} U_n \sqcup O(a_{n+1})$. 
By $O(a_n) = \bigsqcup_{t \in \R} \{ v(t,a_n) \}$ for any $n \in \mathbb{N}$, Claim~\ref{claim:saturation} implies that the subset $U_\infty = \bigsqcup_{t \in \R}v(t,(x,a_1)_T) = \bigsqcup_{t \in \R}v(t,T_+ - \{ a_1 \} ) = \bigsqcup_{a \in T_+ - \{ a_1 \}} O(a)$ is an invariant open disk and an invariant trivial flow box with $\{ \alpha \} = \bigcup_{a \in U_\infty} \alpha(a)$ and $\{ \omega \} = \bigcup_{a \in U_\infty} \omega(a)$ such that the sequence $(a_n)_{n \in \mathbb{N}}$ in $T_+ \subset U_\infty \sqcup O(a_1)$ converges to $x \in \partial T_+ \cap \partial U_\infty$. 
For any $a \in T_+$, since $O(a) \cap O(a') = \emptyset$ for any $a' \neq a \in T_+$, we have $T_+ \cap O(a) = \{ a \}$. 
Put $T_- := T - (T_+ \sqcup \{ x \}) = [a_-, x)_T$. 
%
Then $T = T_- \sqcup [x,a_1]_T$. 
Since $T_+ \cap O(a) = \{ a \}$ for any $a \in T_+$, we have that $ T_+ \cap O(a_n) = \{ a_n \}$ for any $n \in \N$. 
For any $n \in \N$, by $f_{T}^{N_{n}}(a_n) \in T$ and $N_n \in \N$, we obtain that $f_{T}^{N_{n}}(a_n) \in T_- \cap O^+(a_n)$.

\begin{claim}\label{claim:hol_non_ori}
It suffices to prove the lemma in the case where, the intersection $T \cap (U_\infty \sqcup O(a_1))$ is the disjoint union $(x,a_1]_T \sqcup [f_T(a_1),x)_T  = (x,a_1]_T \sqcup  f_T((x,a_1]_T)$ such that $\mathop{\mathrm{dom}}f_T = (x,a_1]_T$ and $\mathop{\mathrm{codom}}f_T = [f_T(a_1),x)_T$. 
\end{claim}

\begin{proof}[Proof of Claim~\ref{claim:hol_non_ori}]
For any $j \in \N_{>1}$, since $T_+ \cap \overline{U_{1,j}} = [a_j,a_1]_T$, by the transversality of the closed transverse arc $T$ and the compactness of $T \cap \overline{U_{1,j}}$ and $\overline{U_{1,j}} = U_{1,j} \sqcup \{\alpha, \omega \} \sqcup O(a_1) \sqcup O(a_j)$, the intersection $T \cap \overline{U_{1,j}}$ is a finite disjoint union of $[a_j,a_1]_T$ and $k_j$ closed subarcs in $T_-$ as shown in Figure~\ref{fig:rot_arc_02} for some $k_j \in \N$. 
\begin{figure}
\begin{center}
\includegraphics[scale=1.1]{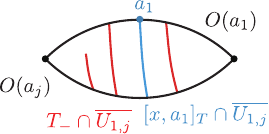}
\end{center}
\caption{The intersection $T \cap \overline{U_{1,j}}$.}
\label{fig:rot_arc_02}
\end{figure} 
For any $j \in \N_{>1}$, if a boundary component of such a closed subarc $\gamma \subseteq T_-$ is contained in $\overline{U_{1,j}}$, then the boundary component of $\gamma$ is $a_- \in T_-$. 
Therefore, such $k_j$ closed subarcs for any $j \in \N_{>1}$ contain at most one closed subarc whose boundary intersects $U_{1,j}$. 
Thus, the transversality of the closed transverse arc $T$ implies that $\vert k_j - k_{j'} \vert \leq 1$ for any $j,j' \in \N_{>1}$. 
Moreover, the intersection $T \cap U_\infty$ is a finite disjoint union of $(x,a_1)_T$ and $k_\infty$ open subarcs $T'_1, \ldots , T'_{k_\infty}$ for some $k_\infty \in \{ k_2, k_2 +1 \}$. 
Denote by $T'_- \in \{ T'_1, \ldots , T'_{k_\infty} \}$ the subarc closest to $x$ in $T_- \sqcup \{ x \} = [a_-,x]_T$. 
Since $\lim_{n \to \infty} f_{T}^{N_{n}}(a_n) = x$ and $f_{T}^{N_{n}}(a_n) \in T_-$, there is a number $L \in \N$ such that $f_{T}^{N_{n}}(a_n) \in T'_-$ for any $n \in \N_{\geq L}$. 
Since the sequence $(f_{T}^{N_{n}}(a_n))_{n \in \N_{\geq L}}$ of points $f_{T}^{N_{n}}(a_n) \in T' \subset [a_-,x)_T$ converges to $x \in [a_-,x]_T$, there is a point $a'_- \in \partial T'_-$ with either $T'_- = (a'_-,x)_T$ or $T'_- = [a'_-,x)_T = [a_-,x)_T = T_-$. 
Then $(f_{T}^{N_{L}}(a_L), x)_T \subseteq T'_-$. 
By taking the subsequence $(a_{n})_{n \in \N_{\geq L}}$ of $(a_{n})_{n \in \N}$ (and so replacing $T_+ = (x,a_1]_T$ with $(x, a_L]_T$), and by replacing $T_-$ with $[f_{T}^{N_{L}}(a_L), x)_T \subseteq T_-$, we may assume that the intersection $T \cap (U_\infty \sqcup O(a_1))$ is the disjoint union of $(x,a_1]_T$ and the open subarc $[f_T(a_1),x)_T = f_T((x,a_1]_T)$. 
Since $(U_\infty \sqcup O(a_1)) \cap T = T - \{ x \} = [f_T(a_1),x)_T \sqcup (x,a_1]_T$, we have that $\mathop{\mathrm{dom}}f_T = (x,a_1]_T$ and $\mathop{\mathrm{codom}}f_T = [f_T(a_1),x)_T$.
\end{proof}

By $f_{T}(a_1) < x < a_1$ in $T$, since $U_\infty \sqcup O(a_1)$ is a trivial flow box containing $\mathop{\mathrm{dom}}f_T = (x,a_1]_T$ and $\mathop{\mathrm{codom}}f_T = [f_T(a_1),x)_T$, the mapping $f_T \colon (x,a_1]_T \to [f_{T}(a_1),x)_T$ is an orientation-reversing homeomorphism. 
Applying Lemma~\ref{lem:bdry_circuit} to $I := (f_{T}(a_1),a_1)_T$, there is a circuit $\nu$ in $\partial U_\infty$ with wandering holonomy at $O(x) \subset \nu$.  
\end{proof}
 
Lemma~\ref{lem:wandering_h} implies the following sufficient condition for the denseness of closed orbits with respect to a flow with finitely many singular points on a surface. 

\begin{lemma}\label{lemp:sufficient_dense}
Let $v$ be a flow with finitely many singular points without non-closed
recurrent points on a compact surface $S$ without strictly limit circuits or circuits with wandering holonomy. 
Then $\overline{\mathop{\mathrm{Cl}}(v)} = \Omega(v)$. 
\end{lemma}

\begin{proof}
By the definition of a non-wandering point, we have that $\overline{\mathop{\mathrm{Cl}}(v)} \subseteq \Omega(v)$. 
Fix any point $x \in \Omega(v)$. 
If $x \in \Omega(v) - \overline{\Cv}$, then Lemma~\ref{lem:wandering_h} implies that there is either a strictly limit circuit at $O(x)$ or a circuit with wandering holonomy at $O(x)$, which contradicts the hypothesis.
Thus, we have that $x \in \overline{\Cv}$. 
\end{proof}

Proposition~\ref{prop:necessary_dense} and Lemma~\ref{lemp:sufficient_dense} imply Theorem~\ref{thm:finite_dense}. 
 
\section{Characterizations of the denseness of closed orbits of flows with finitely many singular points}

We topologically characterize orbits for flows with finitely many singular points on compact surfaces.

\subsection{Topological characterization for the case with finitely many connected components of the singular point set}

To generalize Theorem~\ref{thm:finite_dense} and other main results into a characterization of the denseness of closed orbits for a flow with finitely many connected components of the singular point set, we introduce some concepts (see \cite{yokoyama2017decompositions} for details of these constructions). 

\subsubsection{Blow-downs for singular points}
The blow-down operation for singular points was constructed in \cite{yokoyama2017decompositions}. 
For self-containedness, we describe this operation applied to the singular point set of a flow on a surface, along with a new flow induced by this modification. 
More precisely, we state a simplified construction which implies the same flow on the resulting surface: 

Let $S$ be a surface of finite genus without boundary and $v$ a flow on $S$ whose singular point set has at most finitely many connected components. 
Since the singular point set $\mathop{\mathrm{Sing}}(v)$ is closed, the complement $S - \mathop{\mathrm{Sing}}(v)$ is open and therefore is an open surface. 
From \cite[Theorem~3]{richards1963classification}, any connected component $U$ of the open surface $S - \Sv$ is homeomorphic to a surface obtained from a connected closed surface $S_U$ after removing a closed totally disconnected subset.
By fixing such a closed surface $S_U$ for any connected component $U$ of the open surface $S - \Sv$ and by identifying $U$ as a subset of $S_U$, denote by $\bm{S_{\mathrm{col}}}$ the disjoint union of such closed surfaces (i.e. $S_{\mathrm{col}} := \bigsqcup_{U} S_U$ where $U$ runs over the set of connected components of $S - \Sv$). 
Define $\partial_{\mathrm{col}} := S_{\mathrm{col}} - (S - \mathop{\mathrm{Sing}}(v))$, which is the union of closed totally disconnected subsets. 
Then the subset $S_{\mathrm{col}} - \partial_{\mathrm{col}} = S - \mathop{\mathrm{Sing}}(v)$ is an open dense subset of $S_{\mathrm{col}}$. 
Define a flow $\bm{v_{\mathrm{col}}}$, called {\bf the blow-down flow}, on $S_{\mathrm{col}}$ by $v_{\mathrm{col}}\vert_{S_{\mathrm{col}} - \partial_{\mathrm{col}}} = v\vert_{S - \mathop{\mathrm{Sing}}(v)}$ and $\mathop{\mathrm{Sing}}(v_{\mathrm{col}}) = \partial_{\mathrm{col}}$. 

Let $S$ be a surface of finite genus with boundary and $v$ a flow on $S$ whose singular point set has at most finitely many connected components.  
By taking the lift $v_{\mathrm{db}}$ of the flow $v$ to the double $S_{\mathrm{db}}$ of $S$, the subset $S \subset S_{\mathrm{db}}$ is an invariant subset of $v_{\mathrm{db}}$.
Then define {\bf the blow-down surface} $\bm{S_{\mathrm{col}}}$ by $S_{\mathrm{col}} := (S_{\mathrm{db}})_{\mathrm{col}} \cap \overline{S - \Sv}^{S_\mathrm{col}}$, where $\overline{A}^{S_\mathrm{col}}$ is the closure of a subset $A \subseteq S_\mathrm{col}$ in the surface $S_\mathrm{col}$. 
Then $S_{\mathrm{col}}$ is a topological surface, which is a surface with corners, and invariant with respect to $(v_{\mathrm{db}})_{\mathrm{col}}$. 
By Rad{\'o}'s theorem, we may assume that any connected component of $S_{\mathrm{col}}$ is a compact surface.  
Define {\bf the blow-down flow} $\bm{v_{\mathrm{col}}}$ on the surface $S_{\mathrm{col}}$ by $v_{\mathrm{col}} := (v_{\mathrm{db}})_{\mathrm{col}}\vert_{S_{\mathrm{col}}}$ up to topological equivalence. 

\begin{figure}
\begin{center}
\includegraphics[scale=0.1]{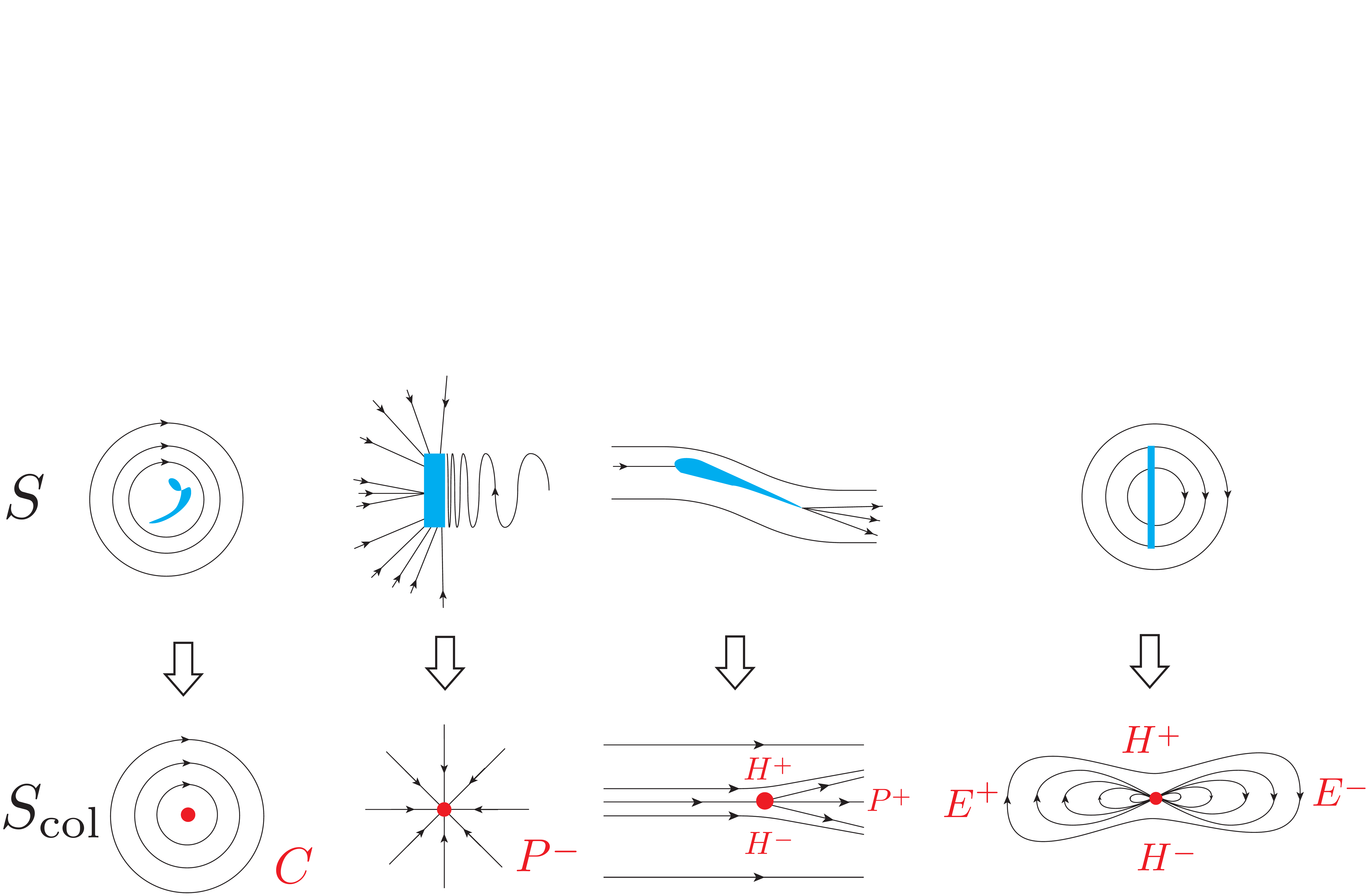}
\end{center}
\caption{Blow-downs of finitely many connected components of the singular point set. Here, the symbol $C$ represents a center and the symbols $P^{\pm}$ (resp. $H^{\pm}$, $E^{\pm}$) represent parabolic (resp. hyperbolic, elliptic) sectors (see \cite{yokoyama2017decompositions} for details).}
\label{blowup}
\end{figure}

In any case, we have that $O_v(x) = O_{v_{\mathrm{col}}}(x)$ for any point $x \in S - \mathop{\mathrm{Sing}}(v) = S_{\mathrm{col}} -  \mathop{\mathrm{Sing}}(v_{\mathrm{col}})$. 
Moreover, the restriction $v_{\mathrm{col}}\vert_{S_U}$ to any connected component $S_U \subseteq S_{\mathrm{col}}$ has at most finitely many singular points.

For any subset $A \subseteq S$, define $A_{\mathrm{col}} := \overline{A \setminus \Sv}^{S_\mathrm{col}}$. 
We have the following observation.

\begin{lemma}\label{lem:_col}
Let $S$ be a surface of finite genus and $v$ a flow on $S$ whose singular point set has at most finitely many connected components. 
Then the following statements hold for any closed subset $A \subseteq S$ and any non-singular point $x \in S$: 
\begin{enumerate}[\rm(1)]
\item $A_{\mathrm{col}} \setminus \mathop{\mathrm{Sing}}(v_{\mathrm{col}}) = A \setminus \Sv$.
\item $\omega(x) \setminus \Sv = \omega_{v_{\mathrm{col}}}(x) \setminus \mathop{\mathrm{Sing}}(v_{\mathrm{col}}) = (\omega(x))_{v_{\mathrm{col}}}  \setminus \mathop{\mathrm{Sing}}(v_{\mathrm{col}})$.
\item If $\omega(x) \setminus \Sv \neq \emptyset$, then $(\omega(x))_{\mathrm{col}} = \omega_{v_{\mathrm{col}}}(x)$.
\end{enumerate}
\end{lemma}

\begin{proof}
By the definition of $v_{\mathrm{col}}$, we have that $S - \mathop{\mathrm{Sing}}(v) = S_{\mathrm{col}} -  \mathop{\mathrm{Sing}}(v_{\mathrm{col}})$.
Fix any closed subset $A$ of $S$. 
Then $A \setminus \Sv \subseteq S - \mathop{\mathrm{Sing}}(v) = S_{\mathrm{col}} -  \mathop{\mathrm{Sing}}(v_{\mathrm{col}})$ is closed with respect to the subspace $S_{\mathrm{col}} -  \mathop{\mathrm{Sing}}(v_{\mathrm{col}}) = S - \mathop{\mathrm{Sing}}(v)$. 
Therefore, we have the following inclusion: 
\[
A_{\mathrm{col}} - (A \setminus \Sv) = \overline{A \setminus \Sv}^{S_\mathrm{col}} - (A \setminus \Sv) \subseteq \mathop{\mathrm{Sing}}(v_{\mathrm{col}})
\]
This implies the following equality: 
\[
\begin{split}
A_{\mathrm{col}} \setminus \mathop{\mathrm{Sing}}(v_{\mathrm{col}}) &= A_{\mathrm{col}} \cap (S_\mathrm{col} \setminus \mathop{\mathrm{Sing}}(v_{\mathrm{col}})) 
\\
&= \overline{A \setminus \Sv}^{S_\mathrm{col}} \cap (S_\mathrm{col} \setminus \mathop{\mathrm{Sing}}(v_{\mathrm{col}})) = A \setminus \Sv
\end{split}
\]
Fix any non-singular point $x \in S$. 
Since $\omega(x)$ is a closed subset of $S$, the following equality holds:
\[
\begin{split}
(\omega(x))_{\mathrm{col}} \setminus \mathop{\mathrm{Sing}}(v_{\mathrm{col}}) = \omega(x) \setminus \Sv &= \bigcap_{n\in \mathbb{R}}\overline{v(\R_{>n},x)}  \setminus \Sv
\\
&= \bigcap_{n\in \mathbb{R}}\overline{v_{\mathrm{col}}(\R_{>n},x)}^{S_\mathrm{col}} \setminus \mathop{\mathrm{Sing}}(v_{\mathrm{col}})
\\
&= \omega_{v_{\mathrm{col}}}(x)  \setminus \mathop{\mathrm{Sing}}(v_{\mathrm{col}})
\end{split}
\]

The closedness of $\omega_{v_{\mathrm{col}}}(x)$ implies the following inclusions: 
\[
(\omega(x))_{\mathrm{col}} = \overline{\omega(x) \setminus \Sv}^{S_{\mathrm{col}}}  = \overline{ (\omega(x))_{\mathrm{col}} \setminus \mathop{\mathrm{Sing}}(v_{\mathrm{col}})}^{S_{\mathrm{col}}} \subseteq \omega_{v_{\mathrm{col}}}(x)
\]

Suppose that $\omega(x) \setminus \Sv \neq \emptyset$.

\begin{claim}\label{claim:col_commute}
$\omega_{v_{\mathrm{col}}}(x) = (\omega(x))_{\mathrm{col}}$. 
\end{claim}

\begin{proof}[Proof of Claim~\ref{claim:col_commute}]
Assume that $\omega_{v_{\mathrm{col}}}(x) - (\omega(x))_{\mathrm{col}} \neq \emptyset$.
Then the set difference $F := \omega_{v_{\mathrm{col}}}(x) - (\omega(x))_{\mathrm{col}} \subseteq \mathop{\mathrm{Sing}}(v_{\mathrm{col}})$ is a finite subset and so is closed.
The normality of $S_{\mathrm{col}}$ implies that there are disjoint open \nbds $U_{(\omega(x))_{\mathrm{col}}}$ and $U_{F}$ of $(\omega(x))_{\mathrm{col}}$ and $F$ respectively. 
Then $\omega_{v_{\mathrm{col}}}(x) = \gamma_{\mathrm{col}} \sqcup F$ is disconnected, which contradicts the connectivity of $\omega_{v_{\mathrm{col}}}(x)$. 
\end{proof}

This completes the proof. 
\end{proof}

\subsubsection{Concepts related to blow-downs}
We recall and define some concepts to state the characterization of the denseness of closed orbits in the non-wandering set. 

\begin{definition}
An $\alpha$-limit or $\omega$-limit set of a point is a  {\bf non-trivial quasi-circuit} if it is a boundary component of an open annulus, contains a non-recurrent point, and consists of non-recurrent points and 
singular points. 
\end{definition}

We have the following correspondence between a non-trivial non-periodic circuit and its blow-up, a non-trivial quasi-circuit. 

\begin{lemma}\label{lem:blowdown_nqc}
Let $v$ be a flow on a compact surface $S$ whose singular point set has at most finitely many connected components. 
Then the $\omega$-limit set $\omega(x)$ of a point $x$ is a non-trivial quasi-circuit if and only if $\omega(x)_\mathrm{col}$ is a non-trivial non-periodic circuit. 
In any case, we have that $\omega(x)_\mathrm{col} = \omega_{v_\mathrm{col}}(x)$. 
\end{lemma}

\begin{proof}
Let $\gamma = \omega(x)$ be the $\omega$-limit set of a point $x$ of a flow $v$ on a compact surface $S$. 
By Lemma~\ref{lem:_col}, we have that $\gamma_{\mathrm{col}} = \omega_{v_{\mathrm{col}}}(x)$.

Suppose that $\gamma$ is a non-trivial quasi-circuit. 
There is an open annulus $\A \subset S_\mathrm{col} - \mathop{\mathrm{Sing}}(v_{\mathrm{col}}) = S - \Sv$ one of whose boundary components is $\gamma$ and contains a non-recurrent point, and consists of non-recurrent points and singular points. 
Since any connected component of $S_\mathrm{col}$ contains at most finitely many singular points, the $\omega$-limit set $\omega_{v_{\mathrm{col}}}(x) = \gamma_{\mathrm{col}}$ consists of non-recurrent points and finitely many singular points with respect to $v_{\mathrm{col}}$. 
Applying Lemma~\ref{lem:lc_transversal} to $v_{\mathrm{col}}$, the $\omega$-limit set $\omega_{v_{\mathrm{col}}}(x) = \gamma_{\mathrm{col}}$ is a limit circuit, which is a non-trivial non-periodic circuit. 

Conversely, suppose that the subset $\gamma_{\mathrm{col}}$ is a non-trivial non-periodic circuit. 
Then $\gamma_{\mathrm{col}}$ contains a non-recurrent orbit and consists of singular points and non-recurrent orbits, and so does the $\omega$-limit set $\gamma$. 
From Lemma~\ref{lem:lc_transversal}, the $\omega$-limit set $\omega_{v_{\mathrm{col}}}(x) = \gamma_{\mathrm{col}}$ is a non-periodic limit circuit such that there is its semi-attracting collar basin $\A_{\mathrm{col}} \subseteq \mathrm{P}(v_{\mathrm{col}}) \subset S_{\mathrm{col}} -  \mathop{\mathrm{Sing}}(v_{\mathrm{col}}) = S - \mathop{\mathrm{Sing}}(v)$.  
By replacing $x$ with a point in $O^+_{v_{\mathrm{col}}}(x) \cap \A_{\mathrm{col}} = O^+(x) \cap \A_{\mathrm{col}}$, we may assume that $x \in \A_{\mathrm{col}} \subset S_{\mathrm{col}} -  \mathop{\mathrm{Sing}}(v_{\mathrm{col}}) = S - \mathop{\mathrm{Sing}}(v)$.
Moreover, the $\omega$-limit set $\omega_{v_{\mathrm{col}}}(x) = \gamma_{\mathrm{col}}$ is a boundary component of $\A_{\mathrm{col}}$ with respect to $S_{\mathrm{col}}$. 
From $O^+_{v_{\mathrm{col}}}(x) = O^+(x) \subset \A_{\mathrm{col}}$, there is a boundary component $\mu \subset S$ of $\A_{\mathrm{col}}$ with respect to $S$ such that $\omega(x) = \gamma \subseteq \mu$. 

\begin{claim}\label{claim:00ka}
$\gamma = \omega(x) = \mu$. 
\end{claim}
\begin{proof}[Proof of Lemma~\ref{claim:00ka}]
By the existence of the semi-attracting collar basin $\A_{\mathrm{col}}$ of $\omega_{v_{\mathrm{col}}}(x)$, there are a non-singular point $y \in \omega(x) \setminus \Sv = \omega_{v_{\mathrm{col}}}(x) \setminus \mathop{\mathrm{Sing}}(v_{\mathrm{col}})$ and a transverse open arc $I$ whose boundary consists of $y$ and a point in $O^+(x)$ such that the first return map $f_{v,I}$ is attracting, the union $U := \bigcup_{z \in I} O^+(z)$ is a collar of $\gamma$ and the set difference $U - I$ is a trivial flow box. 
Since the compact surface $S$ is Fr{\'e}chet-Urysohn, the attraction of $f_{v,I}$ implies that any point in $\mu$ is contained in $\overline{O^+(x)} - (O^+(x) \sqcup \{ x \}) \subseteq \omega(x) = \gamma$. 
This means that $\omega(x) = \gamma = \mu$. 
\end{proof}
Claim~\ref {claim:00ka} implies that $\omega(x)$ is a boundary component of the open annulus $\A_{\mathrm{col}}$ in $S$. 
Since $\omega(x) = \gamma$ contains a non-recurrent orbit and consists of singular points and non-recurrent orbits, the $\omega$-limit set $\omega(x)$ is a non-trivial quasi-circuit.
\end{proof}

Though Claim~\ref{claim:00ka} can be directly deduced from \cite[Theorem~A or Lemma~3.3]{yokoyama2021poincare}, for the self-containedness, since a proof of Claim~\ref{claim:00ka} is not long, we demonstrated it. 

\begin{definition}
A non-trivial quasi-circuit $\gamma$ is a {\bf quasi-semi-attracting} limit quasi-circuit (with respect to a small collar $\A$) if there is a point $x \in \A$ with $O^+(x) \subset \A$ such that $\omega(x)  = \gamma$. 
\end{definition}

\begin{definition}
A non-trivial quasi-circuit $\gamma$ is a {\bf quasi-semi-repelling} limit quasi-circuit (with respect to a small collar $\A$) if there is a point $x \in \A$ with $O^-(x) \subset \A$ such that $\alpha(x)  = \gamma$. 
\end{definition}

\begin{definition}
A non-trivial quasi-circuit is a {\bf limit quasi-circuit} {\rm(}with respect to a small collar $\A${\rm)} if it is a quasi-semi-attracting or quasi-semi-repelling limit quasi-circuit with respect to $\A$. 
\end{definition}

Notice that the above definition of limit quasi-circuit seems to differ from the original definition in \cite{yokoyama2021poincare}; however, the two definitions coincide by \cite[Theorem~A]{yokoyama2021poincare}.
We define one-sidedness for quasi-circuits as follows. 

\begin{definition}
A non-trivial quasi-circuit $\gamma$ is {\bf one-sided} at a non-recurrent orbit $\mu \subset \gamma$ if for any small neighborhood $W$ of $\gamma$ there is a collar $V \subset W$ of $\gamma$ such that the union $V \sqcup \gamma$ is a neighborhood of $\mu$. 
\end{definition}

\begin{definition}
A limit quasi-circuit $\gamma$ is a {\bf strictly limit quasi-circuit} if either $\gamma$ is one-sided or there are a non-recurrent orbit $\mu \subseteq \gamma$ and a transverse open arc $T$ intersecting $\mu$ such that $f_v|_{T_1}$ is either attracting or repelling, and $f_v|_{T_2}$ has no periodic points, where $f_v :T \to T$ is the first return map, $T_1$ and $T_2$ are the two connected components of $T \setminus \mu$. 
\end{definition}

By Lemma~\ref{lem:blowdown_nqc}, we obtain the following equivalence between a strictly limit circuit and its blow-up, a strictly limit quasi-circuit.

\begin{lemma}\label{lem:blowdown_slqc}
Let $v$ be a flow on a compact surface $S$ whose singular point set has at most finitely many connected components. 
Then a closed connected invariant subset $\gamma$ is a strictly limit quasi-circuit if the subset $\gamma_{\mathrm{col}}$ is a strictly limit circuit. 
\end{lemma}

We define a blow-up of a circuit with wandering holonomy as follows. 

\begin{definition}
An invariant subset $\gamma$ is a {\bf blow-up of a circuit with wandering holonomy} if the subset $\gamma_{\mathrm{col}}$ is a circuit with wandering holonomy for $v_{\mathrm{col}}$. 
\end{definition}

Notice that a blow-up of a circuit with wandering holonomy need not be a (non-trivial) quasi-circuit. 
Using the proof of Lemma~\ref{lem:blowdown_nqc}, together with Lemma~\ref{lem:blowdown_slqc} and the definition of the blow-up of a circuit with wandering holonomy, by $S - \mathop{\mathrm{Sing}}(v) = S_{\mathrm{col}} -  \mathop{\mathrm{Sing}}(v_{\mathrm{col}})$, we obtain the following. 

\begin{lemma}\label{lem:blow_down}
The following statements hold for a flow $v$ with finitely many connected components of the singular point set on a compact surface $S$: 
\\
{\rm(1)} $\Omega(v) - \Sv = \Omega(v_{\mathrm{col}}) - \mathop{\mathrm{Sing}}(v_{\mathrm{col}})$, $\mathop{\mathrm{Per}}(v) = \mathop{\mathrm{Per}}(v_{\mathrm{col}})$,  and $\mathrm{R}(v) = \mathrm{R}(v_{\mathrm{col}})$. 
\\
{\rm(2)} $\overline{\mathop{\mathrm{Cl}}(v)} = \Omega(v)$ if and only if $\overline{\mathop{\mathrm{Cl}}(v_{\mathrm{col}})}^{S_\mathrm{col}} = \Omega(v_{\mathrm{col}})$.
\\
{\rm(3)} The $\omega$-limit set $\omega_v(x)$ of a non-singular point $x$ on $S$ is a limit quasi-circuit if and only if the $\omega$-limit set $\omega_{v_{\mathrm{col}}}(x)$ is a non-periodic limit circuit. 
\\
{\rm(4)} 
There are no strictly limit quasi-circuits for $v$ if and only if there are no strictly limit circuits for $v_{\mathrm{col}}$. 
\\
{\rm(5)} There are no blow-ups of circuits with wandering holonomy for $v$ if and only if there are no circuits with wandering holonomy for $v_{\mathrm{col}}$. 
\end{lemma}

\begin{proof}
By the definition of a blow-down surface, we have that $\Omega(v) - \Sv = \Omega(v_{\mathrm{col}}) - \mathop{\mathrm{Sing}}(v_{\mathrm{col}})$. 
Therefore, we obtain that $\mathop{\mathrm{Per}}(v) = \mathop{\mathrm{Per}}(v_{\mathrm{col}})$  and $\mathrm{R}(v) = \mathrm{R}(v_{\mathrm{col}})$. 
Then we have the following equality: 
\[
\begin{split}
\overline{\Cv} - \Cv &= (\Sv \cup \overline{\Pv}) - \Cv 
\\
&= \overline{\Pv} \cap (S - \Cv) 
\\
&= \overline{\mathop{\mathrm{Per}}(v_{\mathrm{col}})}^{S_\mathrm{col}} \cap (S_{\mathrm{col}} - \mathop{\mathrm{Cl}}(v_{\mathrm{col}})) 
\\
&= (\mathop{\mathrm{Sing}}(v_{\mathrm{col}}) \cup \overline{\mathop{\mathrm{Per}}(v_{\mathrm{col}})})^{S_\mathrm{col}} - \mathop{\mathrm{Cl}}(v_{\mathrm{col}}) = \overline{\mathop{\mathrm{Cl}}(v_{\mathrm{col}})}^{S_\mathrm{col}} - \mathop{\mathrm{Cl}}(v_{\mathrm{col}})
\end{split}
\]
Since $\Omega(v) - \Cv = \Omega(v_{\mathrm{col}}) - \mathop{\mathrm{Cl}}(v_{\mathrm{col}})$, the closedness of $\Sv$ and $\mathop{\mathrm{Sing}}(v_{\mathrm{col}})$ implies the following equality: 
\[
\begin{split}
\Omega(v) - \overline{\Cv} &= (\Omega(v)- \Cv) - (\overline{\Cv} - \Cv) 
\\
&= (\Omega(v_{\mathrm{col}}) - \mathop{\mathrm{Cl}}(v_{\mathrm{col}})) - (\overline{\mathop{\mathrm{Cl}}(v_{\mathrm{col}})}^{S_\mathrm{col}} - \mathop{\mathrm{Cl}}(v_{\mathrm{col}}))
\\
&= \Omega(v_{\mathrm{col}}) - \overline{\mathop{\mathrm{Cl}}(v_{\mathrm{col}})}^{S_\mathrm{col}}
\end{split}
\]
This implies assertion (2).
By Lemma~\ref{lem:blowdown_nqc}, assertion (3) holds.
Assertion (4) follows from Lemma~\ref{lem:blowdown_slqc}. 
By the definition of a blow-up of a circuit with wandering holonomy, assertion (5) holds.
\end{proof}

Recall that a topological space is $T_0$ (or Kolmogorov) if for all points $x \neq y \in X$ there is an open subset $U$ of $X$ such that $\vert \{x, y \} \cap U \vert =1$. 
Theorem~\cite[Theorem~3.3]{yokoyama2019properness}, Theorem~\ref{thm:finite_dense}, and Lemma~\ref{lem:blow_down} imply the following generalization of Theorem~\ref{thm:finite_dense}. 

\begin{corollary}\label{thm:main03}
Let $v$ be a flow with finitely many connected components of the singular point set on a compact surface $S$.
The following conditions are equivalent:
\\
$(1)$ $\overline{\mathop{\mathrm{Cl}}(v)} = \Omega(v)$.
\\
$(2)$ There are neither non-closed recurrent points, nor strictly limit quasi-circuits, nor blow-ups of circuits with wandering holonomy.
\\
$(3)$ Each orbit is proper, and there are neither strictly limit quasi-circuits nor blow-ups of circuits with wandering holonomy.
\\
$(4)$ The orbit space $S/v$ is $T_0$, and there are neither strictly limit quasi-circuits nor blow-ups of circuits with wandering holonomy.
\end{corollary}

\subsection{Denseness of closed orbits in the non-wandering set and correspondence for flows on non-compact surfaces}

Since certain physical flows, such as those arising in fluid dynamics, are naturally defined on non-compact domains, we generalize our results to non-compact surfaces.

\subsubsection{End completions of surfaces with finite genus} 

From \cite[Theorem~3]{richards1963classification}, any connected surface of finite genus is homeomorphic to a surface obtained from a closed surface after removing (a finite or infinite number of) closed totally disconnected subsets. 
Therefore, the end compactification $S_{\mathrm{end}}$ of a connected surface $S$ of finite genus is a closed surface. 

For a flow $v$ on a surface $S$ of finite genus, considering ends to be singular points, we obtain the resulting flow $\bm{v_{\mathrm{end}}}$ on the surface $S_{\mathrm{end}}$ which is a union of closed surfaces. 
We obtain the following equality. 

\begin{lemma}\label{lem:corr_end_col}
Let $v$ be a flow with finitely many connected components of the singular point set on a surface $S$ with finite genus and finite ends. 
Then the surfaces $(S_{\mathrm{end}})_{\mathrm{col}}$, $(S_{\mathrm{col}})_{\mathrm{end}}$, and $S_{\mathrm{col}}$ are homeomorphic to each another. 
Moreover, the flows $(v_{\mathrm{end}})_{\mathrm{col}}$, $(v_{\mathrm{col}})_{\mathrm{end}}$, and $v_{\mathrm{col}}$ are topologically equivalent to each another. 
\end{lemma}

\begin{proof}
Suppose that $S$ has no boundary.
Then $(S_{\mathrm{col}})_{\mathrm{end}}$ can be identified with $S_{\mathrm{col}}$. 
Therefore, the flow $(v_{\mathrm{col}})_{\mathrm{end}}$ is topologically equivalent to $v_{\mathrm{col}}$. 
Since $S_{\mathrm{end}} - S \subseteq \mathop{\mathrm{Sing}}(v_{\mathrm{end}})$, we have that $S_{\mathrm{end}} - \mathop{\mathrm{Sing}}(v_{\mathrm{end}}) = S - \Sv$ and so $(S_{\mathrm{end}})_{\mathrm{col}}= S_{\mathrm{col}} = (S_{\mathrm{col}})_{\mathrm{end}}$. 
By the constructions of the flows $(v_{\mathrm{end}})_{\mathrm{col}}$ and $(v_{\mathrm{col}})_{\mathrm{end}}$, the flow $(v_{\mathrm{end}})_{\mathrm{col}}$ is topologically equivalent to $(v_{\mathrm{col}})_{\mathrm{end}}$. 

Thus, we may assume that $\partial S \neq \emptyset$.
Then $((S_{\mathrm{db}})_{\mathrm{col}})_{\mathrm{end}}$ can be identified with $(S_{\mathrm{db}})_{\mathrm{col}}$. 
Therefore, the flow $(v_{\mathrm{col}})_{\mathrm{end}}$ is topologically equivalent to $v_{\mathrm{col}}$. 
Moreover, the flow  $(v_{\mathrm{db}})_{\mathrm{end}}$  can be identified with $(v
_{\mathrm{end}})_{\mathrm{db}}$. 
Then $(S_{\mathrm{end}})_{\mathrm{db}} - S_{\mathrm{db}} = (S_{\mathrm{db}})_{\mathrm{end}} - S_{\mathrm{db}} \subseteq \mathop{\mathrm{Sing}}((v_{\mathrm{db}})_{\mathrm{end}}) = \mathop{\mathrm{Sing}}((v_{\mathrm{end}})_{\mathrm{db}})$. 
This implies that 
\[
(S_{\mathrm{end}})_{\mathrm{db}} - \mathop{\mathrm{Sing}}((v_{\mathrm{end}})_{\mathrm{db}}) = (S_{\mathrm{db}})_{\mathrm{end}} - \mathop{\mathrm{Sing}}((v_{\mathrm{db}})_{\mathrm{end}}) = S_{\mathrm{db}} - \mathop{\mathrm{Sing}}(v_{\mathrm{db}})
\]
and so that $((S_{\mathrm{db}})_{\mathrm{end}})_{\mathrm{col}} = (S_{\mathrm{db}})_{\mathrm{col}} = ((S_{\mathrm{db}})_{\mathrm{col}})_{\mathrm{end}}$.
By the definition of a blow-down surface for a surface with boundary, we obtain $(S_{\mathrm{end}})_{\mathrm{col}} = S_{\mathrm{col}} = (S_{\mathrm{col}})_{\mathrm{end}}$. 
By the constructions of the flows $(v_{\mathrm{end}})_{\mathrm{col}}$ and $(v_{\mathrm{col}})_{\mathrm{end}}$, the flow $(v_{\mathrm{end}})_{\mathrm{col}}$ is topologically equivalent to $(v_{\mathrm{col}})_{\mathrm{end}}$. 
\end{proof}

We have the following observation, which is a correspondence between a flow and the flow obtained by its end completion. 

\begin{lemma}\label{lem:blow_down02}
The following statements hold for a flow $v$ on a surface $S$ with finite genus and finite ends: 
\\
{\rm(1)} Each orbit with respect to $v$ is proper if and only if each orbit with respect to $v_{\mathrm{end}}$ is proper. 
\\
{\rm(2)} The orbit space $S/v$ is $T_0$ if and only if the orbit space $S_{\mathrm{end}}/v_{\mathrm{end}}$ is $T_0$.
\end{lemma}

\begin{proof}
Note that the end completion for $S$ adds finitely many singular points with respect to the resulting flow $v_{\mathrm{end}}$. 
Therefore, assertion (1) is true because a non-singular orbit is proper in $S$ if and only if it is proper in $S_{\mathrm{end}}$. 
\cite[Theorem~3.3]{yokoyama2019properness} implies that each orbit with respect to $v$ (resp. $v_{\mathrm{end}}$) is proper if and only if $S/v$  (resp. $S_{\mathrm{end}}/v_{\mathrm{end}}$) is $T_0$. 
This along with the assertion {\rm(1)} implies the assertion {\rm(2)}. 
\end{proof}

\subsubsection{Virtual quasi-circuits}
Let $v$ be a flow on a surface $S$ with finite genus and finite ends.  
We define some concepts to state limit behaviors on non-compact surfaces.

\begin{definition}
An invariant subset is a {\bf virtual quasi-circuit} if it is the resulting subset from a quasi-circuit on $S_{\mathrm{end}}$ with respect to $v_{\mathrm{end}}$ by removing all the ends. 
\end{definition}
 
\begin{definition}
An invariant subset is a {\bf limit (resp. strictly limit) virtual quasi-circuit} if it is the resulting subset from a limit {\rm(resp.} strictly limit{\rm)} quasi-circuit on $S_{\mathrm{end}}$ with respect to $v_{\mathrm{end}}$ by removing all the ends. 
\end{definition}

\begin{definition}
An invariant subset $\mu$ is a {\bf virtual blow-up of a circuit with wandering holonomy} if there is a blow-up $\gamma$ of a circuit with wandering holonomy on $S_{\mathrm{end}}$ with respect to $v_{\mathrm{end}}$ such that the resulting subset from $\gamma$ by removing all the ends is $\mu$. 
\end{definition}

Corollary~\ref{thm:main03} and Lemma~\ref{lem:blow_down02} imply the following characterization of the denseness of closed orbits for a flow with finitely many connected components of the singular point set on a surface $S$ with finite genus and finite ends. 

\begin{corollary}\label{thm:main}
The following conditions are equivalent for a flow $v$ with finitely many connected components of the  singular point set on a surface $S$ with finite genus and finite ends:
\\
$(1)$ $\overline{\mathop{\mathrm{Cl}}(v)} = \Omega(v)$.
\\
$(2)$
There are neither non-closed recurrent points, nor strictly limit virtual quasi-circuits, nor virtual blow-ups of circuits with wandering holonomy.
\\
$(3)$ Every orbit is proper, and there are neither strictly limit virtual quasi-circuits nor virtual blow-ups of circuits with wandering holonomy.
\\
$(4)$ The orbit space $S/v$ is $T_0$, and there are neither strictly limit virtual quasi-circuits nor virtual blow-ups of circuits with wandering holonomy.
\end{corollary}

\section{Classification of non-closed orbit behaviors}

In this section, we show a generalization of Theorem~\ref{main:rec_wandering_pt}, which is Theorem~\ref{main:rec_wandering}.  
Moreover, Theorem~\ref{main:rec_wandering} implies a characterization of the denseness of recurrence in the non-wandering set, which refines the Poincar{\'e} recurrence theorem for flows with finitely many singular points on surfaces (Theorem~\ref{thm:rec_th}), and implies Corollary~\ref{main_cor:rec_th}.
Furthermore, Theorem~\ref{main:minimal} is demonstrated as a corollary from Theorem~\ref{main:rec_wandering}. 

\subsection{Topological trichotomy of orbits}
We have the following classification of orbits. 

\begin{theorem}\label{main:rec_wandering}
Let $v$ be a flow with finitely many connected components of the singular point set on a surface $S$ with finite genus and finite ends. 
Then one of the following three statements holds for an orbit $O$ exclusively: 
\\
{\rm(1)} The orbit $O$ is recurrent {\rm (i.e.} $O \subseteq \Cv \sqcup \mathrm{R}(v)${\rm)}. 
\\
{\rm(2)} The orbit $O$  is non-recurrent and non-wandering, and satisfies $\overline{O} - O = \alpha(O) \cup \omega(O) \subseteq \Sv$ and either {\rm(2.1)}, {\rm(2.2)}, {\rm(2.3)}, or {\rm(2.4)}: 
\\
{\rm(2.1)} There is a strictly limit virtual quasi-circuit in $\Omega(v)$ at $O \subset \mathop{\mathrm{int}} \mathrm{P}(v)$. 
\\
{\rm(2.2)} There is a virtual blow-up of a circuit with wandering holonomy at $O \subset \mathop{\mathrm{int}} \mathrm{P}(v)$. 
\\
{\rm(2.3)} There is a virtual quasi-circuit in $\overline{\Pv} - \Pv \subseteq \Omega(v)$ containing $O$. 
\\
{\rm(2.4)} The orbit $O \subset \overline{\mathrm{R}(v)} - \mathrm{R}(v)$ is contained in a Q-set. 
\\
{\rm(3)} The orbit $O$ is wandering {\rm(i.e.} $O \cap \Omega(v) = \emptyset${\rm)}, $\overline{O} - O = \alpha(O) \cup \omega(O)$, and $O \subseteq \mathop{\mathrm{int}} \mathrm{P}(v)$. 
\end{theorem}

\begin{proof}
Lemma~\ref{lem:decomp} implies that $S = \mathop{\mathrm{Cl}}(v) \sqcup \mathrm{P}(v) \sqcup \mathrm{R}(v)$. 
Taking the end completion $S_{\mathrm{end}}$ of $S$, let $v_{\mathrm{end}}$ be the resulting flow on the compact surface $S_{\mathrm{end}}$. 
Since $\mathop{\mathrm{Sing}}(v_{\mathrm{end}})$ has at most finitely many connected components, any connected component of the resulting space $S_{\mathrm{ec}} := (S_{\mathrm{end}})_{\mathrm{col}}$ is a compact surface and the restriction of the resulting flow $v_{\mathrm{ec}} := (v_{\mathrm{end}})_{\mathrm{col}}$ to any connected component has at most finitely many singular points. 
Since it suffices to show the assertion for each connected component, we may assume without loss of generality that $S_{\mathrm{ec}}$ is connected.

\begin{claim}\label{claim:closed_surface}
It suffices to prove the theorem when $S_{\mathrm{ec}}$ is a closed surface. 
\end{claim}

\begin{proof}[Proof of Claim~\ref{claim:closed_surface}]
From Lemma~\ref{lem:corr_end_col}, the surface $(S_{\mathrm{end}})_{\mathrm{col}}$ can be identified with $(S_{\mathrm{col}})_{\mathrm{end}}$. 
By taking the lift $(v_{\mathrm{ec}})_{\mathrm{db}}$ of the flow $v_{\mathrm{ec}}$ to the double $(S_{\mathrm{ec}})_{\mathrm{db}}$ of $S_{\mathrm{ec}}$, the subset $S_{\mathrm{ec}} \subset (S_{\mathrm{ec}})_{\mathrm{db}}$ is an invariant subset of $(v_{\mathrm{ec}})_{\mathrm{db}}$ and so the assertion is invariant under taking $(v_{\mathrm{ec}})_{\mathrm{db}}$. 
Thus, by taking $(v_{\mathrm{ec}})_{\mathrm{db}}$ if necessary, we may assume that the surface $S_{\mathrm{ec}}$ is a closed surface. 
\end{proof}

Fix an orbit $O$ on $S$. 

\begin{claim}\label{claim:0l}
It suffices to prove the theorem when $O$ is neither recurrent nor wandering. 
\end{claim}
\begin{proof}[Proof of Claim~\ref{claim:0l}]
If $O$ is recurrent, then assertion {\rm(1)} holds. 
Suppose that $O$ is wandering. 
Then $O \cap \overline{\mathop{\mathrm{Cl}}(v) \sqcup \mathrm{R}(v)} = \emptyset$ and so $O \subseteq S - \overline{\mathop{\mathrm{Cl}}(v) \sqcup \mathrm{R}(v)}  = \mathop{\mathrm{int}} \mathrm{P}(v)$. 
By Lemma~\ref{lem:decomp}, we obtain $(\alpha(O) \cup \omega(O)) \cap O = \emptyset$. 
Since $\overline{O} = O \sqcup (\alpha(O) \cup \omega(O))$ (cf. \cite[Lemma~2.1]{yokoyama2022refinements}), we have $\overline{O} - O = \alpha(O) \cup \omega(O)$. 
This means that assertion {\rm(3)} holds. 
\end{proof}

Since $O$ is non-recurrent and so is non-closed proper, we have $\overline{O} - O = \alpha(O) \cup \omega(O)$ (cf. \cite[Lemma~2.1 and Lemma~4.1]{yokoyama2022refinements}). 

\begin{claim}\label{claim:0m}
It suffices to prove the theorem when $O \cap \overline{\mathrm{R}(v)} = \emptyset$. 
\end{claim}
\begin{proof}[Proof of Claim~\ref{claim:0m}]
Suppose that $O \subset \overline{\mathrm{R}(v)}$. 
Then $O \subset \overline{\mathrm{R}(v)} \cap \mathrm{P}(v) \subseteq \overline{\mathrm{R}(v)} - \mathrm{R}(v)$. 
By the existence of finitely many quasi-minimal sets of flows on compact surfaces, the surface $S$ has at most finitely many quasi-minimal sets. 
Therefore, $O$ is contained in a quasi-minimal set (i.e. assertion {\rm(2.4)} holds). 
\end{proof}

Then we have the following inclusions: 
\[
\begin{split}
O \subseteq  \, & (\mathrm{P}(v) \cap \Omega(v)) \setminus \overline{\Sv \sqcup \mathrm{R}(v)} 
\\
= \, &(\mathrm{P}(v) \cap \Omega(v)) \cap (S - \overline{\Sv \sqcup \mathrm{R}(v)}) 
\\
= \,  &\mathrm{P}(v) \cap \Omega(v) \cap \mathop{\mathrm{int}}(\mathop{\mathrm{Per}}(v) \sqcup \mathrm{P}(v)) 
\end{split}
\]

Since $\mathrm{P}(v) \cap \mathop{\mathrm{int}}(\mathop{\mathrm{Per}}(v) \sqcup \mathrm{P}(v)) \subseteq (\mathrm{P}(v) \cap \partial \Pv) \sqcup \mathop{\mathrm{int}} \mathrm{P}(v)$, we have the following inclusion: 
\[
O \subseteq (\mathrm{P}(v) \cap \partial \Pv) \sqcup (\mathop{\mathrm{int}} \mathrm{P}(v) \cap \Omega(v))
\]

\begin{claim}\label{claim:06}
If $O \subseteq \mathrm{P}(v) \cap \partial \Pv$, then assertion {\rm(2.3)} holds.
\end{claim}
\begin{proof}[Proof of Claim~\ref{claim:06}]
Suppose $O \subseteq \mathrm{P}(v) \cap \partial \Pv$. 
By $\mathrm{P}(v_{\mathrm{ec}}) = \mathrm{P}(v)$ and $\mathop{\mathrm{Per}}(v_{\mathrm{ec}}) = \mathop{\mathrm{Per}}(v)$, we have $O \subset \mathrm{P}(v_{\mathrm{ec}}) \cap \overline{\mathop{\mathrm{Per}}(v_{\mathrm{ec}})}^{S_{\mathrm{ec}}}$. 
Applying Lemma~\ref{lem:circuit} to $v_{\mathrm{ec}}$, the orbit $O$ is contained in a circuit in $\overline{\mathop{\mathrm{Per}}(v_{\mathrm{ec}})}^{S_{\mathrm{ec}}} - \mathop{\mathrm{Per}}(v_{\mathrm{ec}})$ with respect to $v_{\mathrm{ec}}$. 
From $(S_{\mathrm{end}})_{\mathrm{col}} = (S_{\mathrm{col}})_{\mathrm{end}}$, Lemma~\ref{lem:blowdown_nqc} implies that the orbit $O$ is contained in a quasi-circuit in $\overline{\mathop{\mathrm{Per}}(v_{\mathrm{end}})}^{S_{\mathrm{end}}} - \mathop{\mathrm{Per}}(v_{\mathrm{end}})$ with respect to $v_{\mathrm{end}}$. 
By $\mathop{\mathrm{Per}}(v_{\mathrm{end}}) = \mathop{\mathrm{Per}}(v)$, the orbit $O$ is contained in a virtual quasi-circuit in $(\overline{\mathop{\mathrm{Per}}(v_{\mathrm{end}})}^{S_{\mathrm{end}}}  - \mathop{\mathrm{Per}}(v_{\mathrm{end}})) \cap S = \overline{\mathop{\mathrm{Per}}(v)} - \mathop{\mathrm{Per}}(v)$ with respect to $v$. 
\end{proof}

By the previous claim, we may assume $O \subseteq \mathop{\mathrm{int}} \mathrm{P}(v) \cap \Omega(v) \subset S - \overline{\Cv \sqcup \mathrm{R}(v)}$. 
By Lemma~\ref{lem:Maier}, the closure $\overline{\mathrm{R}(v)}$ is a finite union of Q-sets and so has at most finitely many connected components.
Since $S$ is a surface with finite genus and finite ends, any connected component of the complement $S- \overline{\mathrm{R}(v)}$ is a surface with finite genus, even though the complement $S- \overline{\mathrm{R}(v)}$ may have infinitely many connected components. 
Let $U_0$ be the connected component of $S_{\mathrm{ec}} - \overline{\mathrm{R}(v_{\mathrm{ec}})}$ containing $O$. 
From Lemma~\ref{lem:finite_ends}, the connected component $U_0$ is an open subset of $S$ with finitely many ends. 

By taking the end completion $U_1 := (U_0)_{\mathrm{end}}$ and regarding the ends as singular points, we obtain the resulting flow $v_{U_1}$ from $v_{\mathrm{ec}} \vert_{U_0}$ on the compact surface $U_1$ that has no non-closed recurrent orbits.
By the finiteness of $\Sv$ and by the construction of $v_{U_1}$, 
we have that the flow $v_{\mathrm{ec}}$ has at most finitely many singular points.  
Since the surface $U_0$ has at most finite ends, by construction, the flow $v_{U_1}$ has at most finitely many singular points. 
Then $U_1 = \mathop{\mathrm{Cl}}(v_{U_1}) \sqcup \mathrm{P}(v_{U_1})$. 
Applying Lemma~\ref{lem:wandering_h} to $v_{U_1}$, there is a circuit $\gamma$ which is either a strictly limit circuit at $O$ or a circuit with wandering holonomy at $O$ with respect to $v_{U_1}$. 
Moreover, although $U_1$ is not a subset of $S$, we identify $U_0$ with the intersection $U_1 \cap S$.

\begin{claim}\label{claim:07}
If $\gamma$ is a circuit with wandering holonomy at $O$ with respect to $v_{U_1}$, then assertion {\rm(2.2)} holds.
\end{claim}
\begin{proof}[Proof of Claim~\ref{claim:07}]
Suppose that $\gamma$ is a circuit with wandering holonomy at $O$ with respect to $v_{U_1}$. 
Then there are a non-singular point $x \in O \subset \gamma$ and an arbitrarily small open transverse arc $I$ containing $x$ such that the first return map $f_{v,I}$ on $I$ is orientation reversing, the domain is nonempty, and the intersection of the domain of $f_{v,I}$ and the codomain of $f_{v,I}$ is empty. 
Since the intersection of the domain of $f_{v,I}$ and the codomain of $f_{v,I}$ is empty, we obtain $I \subset \mathop{\mathrm{int}} \mathrm{P}(v) \cap U_0$ and so $v(I) = v_{U_1}(I) \subseteq \mathop{\mathrm{int}} \mathrm{P}(v) \cap U_0$. 
Therefore, applying Lemma~\ref{lem:bdry_circuit} to $v$ and $I$, there is a circuit with wandering holonomy at $O$ with respect to $v$. 
\end{proof}

From the previous claim, we may assume that $\gamma$ is a strictly limit circuit at $O$ with respect to $v_{U_1}$. 
%
%
Then it suffices to show that assertion {\rm(2.1)}  holds.
Indeed, by time reversion if necessary, we may assume that $\gamma$ is the $\omega$-limit set of a point $z \in \mathrm{P}(v_{U_1})$ such that $\gamma = \omega_{v_{U_1}}(z) \subset \mathop{\mathrm{Sing}}(v_{U_1}) \sqcup \mathrm{P}(v_{U_1})$. 
Then $z \in \mathrm{P}(v_{U_1}) = \mathrm{P}(v_{U_0}) \subseteq \mathrm{P}(v_{\mathrm{ec}})$.

\begin{claim}\label{claim:omega_per_empty}
$\omega_{v_{\mathrm{ec}}}(z) \cap \mathop{\mathrm{Per}(v_{\mathrm{ec}})} = \emptyset$.
\end{claim}

\begin{proof}[Proof of Claim~\ref{claim:omega_per_empty}]
Since the total number of Q-sets for $v$ is finite due to Lemma~\ref{lem:Maier}, \cite[Proposition~2.2]{yokoyama2016topological} implies $\overline{\mathrm{R}(v)} \cap \mathop{\mathrm{Per}(v)} = \emptyset$ and therefore $\overline{\mathrm{R}(v_{\mathrm{ec}})} \cap \mathop{\mathrm{Per}(v_{\mathrm{ec}})} = \emptyset$ due to the construction of $v_{\mathrm{ec}}$. 
By $U_1 = (U_0)_{\mathrm{end}}$, we obtain that $\mathop{\mathrm{Per}(v_{U_1})} = \mathop{\mathrm{Per}(v_{U_0})} = \mathop{\mathrm{Per}(v_{\mathrm{ec}})} \cap U_1 \subseteq \mathop{\mathrm{Per}(v)}$. 
Since $U_0$ is a connected component of $S_{\mathrm{ec}} - \overline{\mathrm{R}(v_{\mathrm{ec}})}$ and contains $z$, we obtain that $\omega_{v_{\mathrm{ec}}}(z) \setminus \overline{\mathrm{R}(v_{\mathrm{ec}})} \subseteq U_0 \subseteq \omega_{v_{U_1}}(z)$. 
Therefore, the following inequality holds: 
\[
\begin{split}
\omega_{v_{\mathrm{ec}}}(z) \cap \mathop{\mathrm{Per}(v_{\mathrm{ec}})} & =  (\omega_{v_{\mathrm{ec}}}(z) \setminus \overline{\mathrm{R}(v_{\mathrm{ec}})} ) \cap \mathop{\mathrm{Per}(v_{\mathrm{ec}})} 
\\
& \subseteq \omega_{v_{U_1}}(z) \cap \mathop{\mathrm{Per}(v_{\mathrm{ec}})} 
\\
& = \omega_{v_{U_1}}(z) \cap U_1 \cap \mathop{\mathrm{Per}(v_{\mathrm{ec}})} 
\\
& = \omega_{v_{U_1}}(z) \cap \mathop{\mathrm{Per}(v_{U_1})} 
\\
& \subseteq (\mathop{\mathrm{Sing}}(v_{U_1}) \sqcup \mathrm{P}(v_{U_1})) \cap \mathop{\mathrm{Per}(v_{U_1})} 
= \emptyset
\end{split}
\]
\end{proof}

Moreover, there is a transverse open arc $T \subset U_0$ as in the definition of strictly limit circuit with $\mu=O$. 

\begin{claim}\label{claim:07-}
The $\omega$-limit set $\omega_{v_{\mathrm{ec}}}(z)$ is a limit circuit in $\mathop{\mathrm{Sing}(v_{\mathrm{ec}})} \sqcup \mathrm{P}(v_{\mathrm{ec}})$ such that $\omega_{v_{\mathrm{ec}}}(z) = (\omega_{v_{\mathrm{end}}}(z))_{\mathrm{col}}$. 
\end{claim}

\begin{proof}[Proof of Claim~\ref{claim:07-}]
Lemma~\ref{lem:lc_transversal} implies that $\omega_{v_{U_1}}(z) = \gamma$ has its small semi-attracting collar basin $\A \subset \mathrm{P}(v_{U_1})$. 
By replacing a point in $\A \cap O^+(z)$ if necessary, we may assume that $z \in \A$. 
From constructions of $v_{U_1}$ and $v_{\mathrm{ec}}$, we can identify that $\A \subset \mathrm{P}(v_{U_1}) = \mathrm{P}(v_{U_0}) \subseteq \mathrm{P}(v_{\mathrm{ec}}) = \mathrm{P}(v) \subset S$. 
Therefore, there is an open transverse arc $I \subset \A \subset S$ one of whose boundary components is a point $z_\infty$ in $O \subset \omega_v(z) \subset \partial \A$ such that the first return map $f_{v,I}$ on 
$I$ is attracting to $z_\infty$. 
By $O \subseteq \mathop{\mathrm{int}} \mathrm{P}(v) = \mathop{\mathrm{int}} \mathrm{P}(v_{\mathrm{ec}}) \subseteq S_{\mathrm{ec}} - \overline{\mathrm{R}(v_{\mathrm{ec}})}$ and $O \subseteq \omega_{v_{\mathrm{ec}}}(z)$, the $\omega$-limit set $\omega_{v_{\mathrm{ec}}}(z)$ is not a Q-set, which is the closure of a non-closed recurrent orbit.
Since $O \subseteq \omega_{v_{\mathrm{ec}}}(z) \cap \mathrm{P}(v_{\mathrm{ec}})$, Lemma~\ref{lem:lc_transversal} implies that $\omega_{v_{\mathrm{ec}}}(z)$ is a limit circuit in $\mathop{\mathrm{Sing}(v_{\mathrm{ec}})} \sqcup \mathrm{P}(v_{\mathrm{ec}})$.
Applying Lemma~\ref{lem:blowdown_nqc} to $v_{\mathrm{end}}$, we obtain that $(\omega_{v_{\mathrm{end}}}(z))_{\mathrm{col}} = \omega_{v_{\mathrm{ec}}}(z)$. 
\end{proof}

\begin{claim}\label{claim:vqc_at}
The $\omega$-limit set $\omega_v(z)$ is a strictly limit virtual quasi-circuit at $O$.
\end{claim}

\begin{proof}[Proof of Claim~\ref{claim:vqc_at}]
From $O \cup T \subset \mathop{\mathrm{int}}\mathrm{P}(v) \cap U_0$, we obtain that the restrictions of $v$, $v_{\mathrm{ec}}$ and $v_{U_1}$ to $v(T) = v_{U_1}(T) \subseteq \mathop{\mathrm{int}}\mathrm{P}(v) \cap U_0$ correspond to each other. 
Since $(\omega_{v_{\mathrm{end}}}(z))_{\mathrm{col}} = \omega_{v_{\mathrm{ec}}}(z)$ due to Claim~\ref{claim:07-}, by Lemma~\ref{lem:blowdown_slqc}, the $\omega$-limit set $\omega_{v_{\mathrm{end}}}(z)$ is a strictly limit quasi-circuit at $O$.
Then $\omega_v(z) = \omega_{v_{\mathrm{end}}}(z) \cap S = \omega_{v_{\mathrm{end}}}(z) \setminus (S_{\mathrm{end}} - S)$ is a strictly limit virtual quasi-circuit at $O$.
\end{proof}

This completes the proof. 
\end{proof}

Theorem~\ref{main:rec_wandering_pt} follows from Theorem~\ref{main:rec_wandering} via the end completion and the blow-down operation for singular points. 
In the spherical case, we have the following classification of non-closed orbits. 
 
\begin{corollary}\label{cor:wandering_pt}
Let $v$ be a flow with finitely many connected components of the singular point set on a compact surface contained in a sphere. 
Then one of the following three statements holds for an orbit $O$ exclusively: 
\\
{\rm(1)} The orbit $O$ is closed and recurrent {\rm (i.e.} $O \subseteq \Cv${\rm)}. 
\\
{\rm(2)} The orbit $O$  is non-recurrent and non-wandering, and satisfies $\overline{O} - O = \alpha(O) \cup \omega(O) \subseteq \Sv$ and either {\rm(2.1)}, {\rm(2.2)}, or {\rm(2.3)}: 
\\
{\rm(2.1)} There is a strictly limit quasi-circuit in $\Omega(v)$ at $O \subset \mathop{\mathrm{int}} \mathrm{P}(v)$. 
\\
{\rm(2.2)}  There is a quasi-circuit in $\overline{\Pv} - \Pv \subseteq \Omega(v)$ containing $O$. 
\\
{\rm(3)} The orbit $O$ is wandering, $\overline{O} - O = \alpha(O) \cup \omega(O)$, and $O \subseteq \mathop{\mathrm{int}} \mathrm{P}(v)$. 
\end{corollary}
 
\begin{proof}
By the definition of wandering holonomy, the orientability of the sphere implies the non-existence of quasi-circuits with wandering holonomy.  
Since any flow on a compact surface contained in a sphere has no non-closed recurrent orbits, Theorem~\ref{main:rec_wandering} implies the assertion. 
\end{proof}

\subsection{Refinement of Poincar{\'e} recurrence theorem for flows with finitely many singular points}

We have the following characterization of the denseness of recurrence in the non-wandering set, which refines the Poincar{\'e} recurrence theorem.

\begin{theorem}\label{thm:rec_th}
The following are equivalent for a flow with finitely many connected components of the singular point set on a surface with finite genus and finite ends: 
\\
{\rm(1)} The set of recurrent points is dense in the non-wandering set.
\\
{\rm(2)} There are neither strictly limit virtual quasi-circuits nor virtual blow-ups of circuits with wandering holonomy. 
\\
In any of the above cases, the set of wandering points corresponds to the interior of the set of non-recurrent points. 
\end{theorem}

\begin{proof}
Let $v$ be a flow with finitely many connected components of the singular point set on a surface $S$ with finite genus and finite ends. 
Taking the end completion $S_{\mathrm{end}}$ of $S$, let $v_{\mathrm{end}}$ be the resulting flow on the compact surface $S_{\mathrm{end}}$. 
Since $\mathop{\mathrm{Sing}}(v_{\mathrm{end}})$ has at most finitely many connected components, any connected component of the resulting space $S_{\mathrm{ec}} := (S_{\mathrm{end}})_{\mathrm{col}}$ is a compact surface and, the restriction to the connected component of the resulting flow $v_{\mathrm{ec}} := (v_{\mathrm{end}})_{\mathrm{col}}$ has at most finitely many singular points. 
Since it suffices to show the assertion for each connected component, we may assume without loss of generality that $S_{\mathrm{ec}}$ is connected.
Moreover, by Lemma~\ref{lem:corr_end_col}, the surface $(S_{\mathrm{end}})_{\mathrm{col}}$ can be identified with $(S_{\mathrm{col}})_{\mathrm{end}}$. 
By Theorem~\ref{main:rec_wandering}, assertion (2) implies assertion (1).

Therefore, it suffices to show that  assertion (1) implies assertion (2). 
Indeed, suppose assertion (1) holds. 
Assume that there is a strictly limit virtual quasi-circuit $\gamma_v$ of $v$. 
By the definition of a strictly limit virtual quasi-circuit, there is a strictly limit quasi-circuit $\gamma_{\mathrm{end}}$ of $v_{\mathrm{end}}$ on the compact surface $S_{\mathrm{end}}$ with $\gamma_{\mathrm{end}} \cap S = \gamma$. 
Applying Lemma~\ref{lem:blowdown_slqc} to the flow $v_{\mathrm{ec}}$ on the compact surface $S_{\mathrm{ec}}$, the subset $(\gamma_{\mathrm{end}})_{\mathrm{col}}$ is a strictly limit circuit. 
From Lemma~\ref{lem:slc_01}, we have that $\Omega(v_{\mathrm{ec}}) - \overline{\mathop{\mathrm{Cl}}(v_{\mathrm{ec}}) \sqcup \mathrm{R}(v_{\mathrm{ec}})} \neq \emptyset$.
Since $S_{\mathrm{end}} - \mathop{\mathrm{Sing}}(v_{\mathrm{end}}) = S_{\mathrm{ec}} - \mathop{\mathrm{Sing}}(v_{\mathrm{ec}})$, applying Lemma~\ref{lem:blow_down} to $v_{\mathrm{end}}$, we have that $\Omega(v_{\mathrm{end}}) - \overline{\mathop{\mathrm{Cl}}(v_{\mathrm{end}}) \sqcup \mathrm{R}(v_{\mathrm{end}})} \neq \emptyset$.
Because $S_{\mathrm{end}} - S$ consists of finitely many singular points of $v_{\mathrm{end}}$, we obtain $\Omega(v) - \overline{\mathop{\mathrm{Cl}}(v) \sqcup \mathrm{R}(v)} \neq \emptyset$, which contradicts the condition (1).
Thus, there are no strictly limit virtual quasi-circuits of $v$. 
Assume that there is a virtual blow-up $\gamma$ of a circuit with wandering holonomy with respect to $v$. 
By the definition of a virtual blow-up of a circuit with wandering holonomy, there is a blow-up $\gamma_{\mathrm{end}}$ of a circuit with wandering holonomy on the compact surface $S_{\mathrm{end}}$ with $\gamma_{\mathrm{end}} \cap S = \gamma$. 
From the definition of a blow-up of a circuit with wandering holonomy, the subset $(\gamma_{\mathrm{end}})_{\mathrm{col}}$ is a circuit with wandering holonomy of $S_{\mathrm{ec}}$. 
Lemma~\ref{lem:wandering_hol} implies that $\Omega(v_{\mathrm{ec}}) - \overline{\mathop{\mathrm{Cl}}(v_{\mathrm{ec}}) \sqcup \mathrm{R}(v_{\mathrm{ec}})} \neq \emptyset$.
By the same argument as above, we have $\Omega(v) - \overline{\mathop{\mathrm{Cl}}(v) \sqcup \mathrm{R}(v)} \neq \emptyset$, which contradicts the condition (1).
This means that the assertion holds. 
%
\end{proof}

The previous theorem implies Corollary~\ref{main_cor:rec_th} and is utilized later in the proof of Corollary~\ref{main_cor:a} below. 
Note that the existence of finitely many connected components of the singular point set is necessary.
In fact, there are toral flows $v$ without strictly limit circuits nor circuits with wandering holonomy such that $\mathrm{R}(v) = \emptyset$ and $\overline{\mathop{\mathrm{Cl}}(v)} \subsetneq \Omega(v)$ (see examples in Lemma~\ref{lem:counter_exam01} and Lemma~\ref{lem:counter_exam00}).

\subsection{Characterization of minimality}

The characterization of minimality for flows on connected surfaces follows from  Theorem~\ref{main:rec_wandering_pt} as follows. 

\begin{proof}[Proof of Theorem~\ref{main:minimal}]
Let $v$ be a non-wandering flow on a connected surface with finite genus and finite ends. 
Suppose that $v$ is minimal. 
Then $v$ is non-wandering and $S = \mathrm{LD}(v)$. 
Therefore $\Sv = \emptyset$. 
Fix any orbit $O$. 
From \cite[Theorem~VI]{cherry1937topological}, the closure $\overline{O}$ contains uncountably many non-closed recurrent orbits and so $\emptyset \neq \overline{O} - O = S -  O = \mathrm{LD}(v) - O \not\subseteq \Sv$. 

Conversely, suppose that $v$ is non-wandering and the set difference $\overline{O} - O$ for any orbit $O$ is not contained in the singular point set.
Since the set difference $\overline{O} - O$ for any closed orbit $O$ is empty, there are no closed orbits (i.e. $S = \mathrm{P}(v) \sqcup \mathrm{R}(v)$). 
From Theorem~\ref{main:rec_wandering_pt}, the non-existence of singular points implies that any orbits are recurrent and so $S = \mathrm{R}(v)$. 
Since $S_{\mathrm{end}}$ consists of non-closed recurrent points in $S = \mathrm{R}(v)$ and finitely many singular points $S_{\mathrm{end}} - S$ of $v_{\mathrm{end}}$, applying Lemma~\ref{lem:top23} to $v_{\mathrm{end}}$ on the compact surface $S_{\mathrm{end}}$, we get that $\mathrm{E}(v) = \mathrm{E}(v_{\mathrm{end}}) = \emptyset$ and so that $S = \mathrm{LD}(v)$. 
The local denseness implies that the closure $\overline{O(x)}$ for any point $x \in \mathrm{LD}(v) = S$ is a neighborhood of $x$. 
Fix an orbit $O$. 
By \cite[Proposition~2.2]{yokoyama2016topological}, the union $\{ x \in S \mid \overline{O(x)} = \overline{O} \} = \overline{O}$ is a neighborhood of every point of $\overline{O}$ and so is open. 
Since any nonempty open and closed subset of the connected surface $S$ is the whole surface $S$, we obtain $\overline{O} = S$. 
This means that $v$ is minimal. 
\end{proof}

\section{Characterization of the existence of finitely many non-wandering orbits and the correspondence between non-wandering properties and closedness}

This section focuses on the finiteness of the non-wandering set and clarifies its relationship with the closed point set. 
In particular, we establish Corollary~\ref{main_cor:a} regarding the finiteness of non-wandering orbits, which essentially follows from Corollary~\ref{main_cor:rec_th}. 
We also prove Corollary~\ref{main:c} concerning the correspondence between non-wandering properties and closedness, which follows from Theorem~\ref{main:rec_wandering_pt}. 
To begin with, to prove Corollary~\ref{main:c}, we state the following.

\begin{lemma}\label{lem:circuit_w_wh_f}
Let $v$ be a flow with finitely many singular points on a compact surface. 
If each recurrent orbit is a closed orbit, then there are at most finitely many non-wandering orbits with wandering holonomy. 
\end{lemma}

\begin{proof}
Let $v$ be a flow with finitely many singular points on a compact surface with $\mathrm{R}(v) = \emptyset$. 
Then $S = \Cv \sqcup \mathrm{P}(v)$. 
Since any point whose $\omega$-limit or $\alpha$-limit set is a limit circuit is wandering, Lemma~\ref{lem:lc_transversal} implies that the $\omega$-limit and $\alpha$-limit sets of any non-recurrent non-wandering points are singular points. 

Assume that there are infinitely many non-wandering orbits $O_n$ ($n \in \mathbb{N}$) with wandering holonomy. 
Then orbits $O_n$ are non-recurrent and non-wandering. 
By Theorem~\ref{main:rec_wandering_pt}, the finiteness of singular points implies that $O_n$ are separatrices. 
By the finiteness of singular points, there is a singular point $\omega$ (resp. $\alpha$) which is the $\omega$-limit (resp. $\alpha$-limit) set of $O_n$ for infinitely many $n \in \mathbb{N}$. 
By taking a subsequence of $(O_n)_{n \in \mathbb{N}}$, we may assume that $\omega(O_n) = \omega$ and $ \alpha(O_n) = \alpha$ for any $n \in \mathbb{N}$. 
For any $n \in \mathbb{N}$, fix a connected component $D_n$ of $S - (\{ \alpha, \omega \} \sqcup \bigsqcup_{n \in \mathbb{N}} O_n)$ with $\partial D_n = \{ \alpha, \omega \} \sqcup O_n \sqcup O_{n+1}$. 
Since $S$ is of finite genus, by finitely many singular points, by taking a subsequence of $(O_n)_{n \in \mathbb{N}}$, we may assume that any $D_n$ is an invariant open disk with $D_n \subset S - \Sv$ such that any two $D_n$, $D_m$ with $n \neq m$ are disjoint. 
For any $n \in \mathbb{N}$, by $\{ \alpha, \omega \} \sqcup O_{n+1} = \partial D_n \cap \partial D_{n+1}$, the union $D_n \sqcup O_{n+1} \sqcup D_{n+1}$ is an invariant open disk which is a \nbd of $O_{n+1}$. 
Since $O_2$ is a non-wandering orbit with wandering holonomy, there is a non-singular point $x \in O_2$ and there is an arbitrarily small open transverse arc $I \subset D_n \sqcup O_{n+1} \sqcup D_{n+1}$ containing $x$ such that the first return map $f_{v,I}$ on $I$ is orientation reversing, the domain of $f_{v,I}$ is nonempty, the intersection of the domain of $f_{v,I}$ and the codomain of $f_{v,I}$ is empty, and the set difference $I - \{x \}$ consists of open transverse arcs $I_1 \subset D_1$ and $I_2 \subset D_2$. 
By $\{ \alpha, \omega \} \sqcup O_{2} = \overline{O_2}$, we have $x \notin \mathrm{dom}(f_{v,I})$. 
Since $D_1$ (resp. $D_2$) is invariant, the restriction of $f_I$ to $\mathrm{dom}(f_{v,I}) \cap I_1$ (resp. $\mathrm{dom}(f_{v,I}) \cap I_2$) is orientation-preserving.
Therefore, $f_{v,I}$ is orientation-preserving, which contradicts the orientation-reversing property of $f_{v,I}$. 
Thus, the assertion holds. 
\end{proof}

We have the following sufficient condition for a circuit to be a strictly limit circuit.

\begin{lemma}\label{lem:induced_str_lim_circuit}
Let $v$ be a flow with finitely many singular points on a compact surface $S$. 
Suppose that the closed point set $\mathop{\mathrm{Cl}}(v)$ is closed and that there are no non-closed recurrent points.
For any non-periodic limit circuit $\gamma$ and any separatrix $\mu \subset \gamma$, the circuit $\gamma$ is a strictly limit circuit at $\mu$. 
\end{lemma}

\begin{proof}
By the non-existence of non-closed recurrent points, we have $S = \Cv \sqcup \mathrm{P}(v)$. 
For any non-recurrent point $x$ in $\Omega(v) - \Cv$, by Theorem~\ref{main:rec_wandering_pt}, there is a circuit that is either a strictly limit circuit at $O(x)$ or a circuit with wandering holonomy at $O(x)$. 
This means that $\Omega(v)$ consists of closed orbits, strictly limit circuits, and non-wandering orbits with wandering holonomy. 
For any non-periodic limit circuit $\gamma$ and any separatrix $\mu \subset \gamma$, since $\mu \subset \Omega(v) \cap \mathrm{P}(v) \subseteq \Omega(v) - \Cv = \Omega(v) - \overline{\Cv}$, Lemma~\ref{lem:lim_cycle_strict} implies that $\gamma$ is a strictly limit circuit at $\mu$. 
\end{proof}

We demonstrate Corollary~\ref{main:c} as follows.

\begin{proof}[Proof of Corollary~\ref{main:c}]
By Lemma~\ref{lem:induced_str_lim_circuit}, assertions (2) and (3) are equivalent.
Theorem~\ref{main:ch_finite_limit_sets} implies that assertions (3), (4), and (5) are equivalent. 

\begin{claim}\label{claim:00r}
Assertion {\rm(3)} implies assertion {\rm(1)}. 
Moreover, the non-wandering set $\Omega(v)$ consists of closed orbits, strictly limit circuits, and non-wandering orbits with wandering holonomy. 
\end{claim}
\begin{proof}[Proof of Claim~\ref{claim:00r}]
Suppose that $\mathrm{R}(v) = \emptyset$, any limit circuit consists of at most finitely many orbits, and there are at most finitely many closed orbits (i.e. assertion (3) holds). 
Since any limit circuit is contained in $\Omega(v)$, by $\overline{\Pv} = \Pv$ and $\Cv \subseteq \Omega(v)$, Theorem~\ref{main:rec_wandering_pt} implies that $\Omega(v)$ consists of closed orbits, strictly limit circuits, and non-wandering orbits with wandering holonomy. 
From Theorem~\ref{main:fininite_existence_limit_circuits}, there are at most finitely many strictly limit circuits. 
By $\mathrm{R}(v) = \emptyset$, the finiteness of singular points and Lemma~\ref{lem:circuit_w_wh_f} imply that there are at most finitely many non-wandering orbits with wandering holonomy. 
From the hypothesis that any limit circuit consists of at most finitely many orbits, by the existence of finitely many closed orbits and strictly limit circuits, the non-wandering set $\Omega(v)$ consists of finitely many orbits. 
Therefore assertion {\rm(1)} holds
\end{proof}

\begin{claim}\label{claim:00p}
Assertion {\rm(1)} implies assertion {\rm(3)}. 
\end{claim}
\begin{proof}[Proof of Claim~\ref{claim:00p}]
Suppose that $\Omega(v)$ consists of finitely many orbits (i.e. assertion (1) holds). 
Since any limit circuits are contained in $\Omega(v)$, the finiteness of $\Omega(v)$ implies that any limit circuits consist of at most finitely many orbits. 
By $\Cv \subseteq \Omega(v)$, the finiteness of $\Omega(v)$ implies that $\Cv$ consists of finitely many orbits with $\overline{\Pv} = \Pv$. 
By $\Cv \sqcup \mathrm{R}(v) \subseteq \Omega(v)$,  \cite[Lemma 3.1]{yokoyama2019properness} implies $\mathrm{R}(v) = \emptyset$. 
Therefore assertion {\rm(3)} holds. 
\end{proof}
This completes the proof. 
\end{proof}

\subsection{Correspondence between non-wandering properties and closedness}

We first provide a sufficient condition for the non-wandering set to coincide with the set of closed orbits.


\begin{lemma}\label{lem:correspond_pn}
Let $v$ be a flow with finitely many connected components of the singular point set on a compact surface $S$.
Suppose that $\mathop{\mathrm{Cl}}(v)$ is closed and there are neither non-closed recurrent points, nor strictly limit quasi-circuits, nor blow-ups of circuits with wandering holonomy. 
Then $\mathop{\mathrm{Cl}}(v) = \Omega(v)$. 
\end{lemma}

\begin{proof}
The non-existence of non-closed recurrent points implies that $S = \Cv \sqcup \mathrm{P}(v)$. 
Since $\Cv$ is closed, by Theorem~\ref{thm:rec_th}, we have $\Cv = \overline{\Cv \sqcup \mathrm{R}(v)} = \Omega(v)$. 
\end{proof}

We establish the following characterization of the equality $\mathop{\mathrm{Cl}}(v) = \Omega(v)$.

\begin{corollary}\label{cor_cor:a}
The following conditions are equivalent for a flow $v$ with finitely many connected components of the singular point set on a compact surface $S$:
\\
$(1)$ $\mathop{\mathrm{Cl}}(v) = \Omega(v)$.
\\
$(2)$ The closed point set $\mathop{\mathrm{Cl}}(v)$ is closed, and there are neither non-closed recurrent points, nor strictly limit quasi-circuits, nor blow-ups of circuits with wandering holonomy.
\\
$(3)$ The set $\mathrm{P}(v)$ of non-recurrent points is open, and there are neither non-closed recurrent points, nor strictly limit quasi-circuits, nor blow-ups of circuits with wandering holonomy.

In any of the above cases, we have $S = \mathop{\mathrm{Cl}}(v) \sqcup \mathrm{P}(v) = \Omega(v) \sqcup \mathrm{P}(v)$.
\end{corollary}

\begin{proof}
By Lemma~\ref{lem:correspond_pn}, assertion (2) implies assertion (1). 
Suppose that $\mathop{\mathrm{Cl}}(v) = \Omega(v)$.
Applying Corollary~\ref{prop:necessary_dense02} to $v_{\mathrm{col}}$, the closed point set $\mathop{\mathrm{Cl}}(v_{\mathrm{col}})$ is closed and there are neither non-closed recurrent points, nor non-periodic limit circuits, nor circuits with wandering holonomy with respect to $v_{\mathrm{col}}$.
From $S - \mathop{\mathrm{Sing}}(v) = S_{\mathrm{col}} -  \mathop{\mathrm{Sing}}(v_{\mathrm{col}})$, the closed point set $\mathop{\mathrm{Cl}}(v)$ is closed. 
Lemma~\ref{lem:blowdown_nqc} implies that there are no strictly limit quasi-circuits. 
By the definition of a blow-up of a circuit with wandering holonomy, there are no blow-ups of circuits with wandering holonomy with respect $v$. 
This means that assertion (1) implies assertion (2). 
If $v$ has no non-recurrent orbits, then $S = \Cv \sqcup \mathrm{P}(v)$, which implies the equivalence of assertions (2) and (3).
\end{proof}

Corollary~\ref{main_cor:a} follows directly from the previous corollary. 
Moreover, Lemma~\ref{lem:blow_down02} and Corollary~\ref{cor_cor:a} imply the following correspondence.

\begin{corollary}\label{thm:main02}
The following conditions are equivalent for a flow $v$ with finitely many connected components of the  singular point set on a surface $S$ with finite genus and finite ends:
\\
$(1)$ $\mathop{\mathrm{Cl}}(v) = \Omega(v)$.
\\
$(2)$ The closed point set $\mathop{\mathrm{Cl}}(v)$ is closed, and there are neither non-closed recurrent points, nor strictly limit virtual quasi-circuits, nor virtual blow-ups of circuits with wandering holonomy.
\\
$(3)$ The closed point set $\mathop{\mathrm{Cl}}(v)$ is closed, every orbit is proper, and there are neither strictly limit virtual quasi-circuits nor virtual blow-ups of circuits with wandering holonomy.
\\
$(4)$ The closed point set $\mathop{\mathrm{Cl}}(v)$ is closed, the orbit space $S/v$ is $T_0$, and there are neither limit virtual quasi-circuits nor virtual blow-ups of circuits with wandering holonomy.

In any of the above cases, the set $\mathrm{P}(v)$ of non-recurrent points is open. 
\end{corollary}

\section{Examples} 

First, we show that the finiteness of non-degenerate $\omega$-limit sets and that of  $\alpha$-limit sets are independent as follows.

\subsection{Independency of finiteness for non-degenerate $\omega$-limit sets and $\alpha$-limit sets}\label{example:001}

We have the following example, which implies the necessity of the finiteness of singular points in Theorem~\ref{main:inifinie_limit_sets}. 

\begin{lemma}\label{lem:exam01}
There is a flow on a sphere with infinitely many non-degenerate $\omega$-limit sets of points but no non-degenerate $\alpha$-limit sets of points. 
\end{lemma}

\begin{proof}
Consider a basin of a sink of a flow. 
By replacing a sink with a limit circuit, as in Figure~\ref{fig:ex_infinite_sing}, one can obtain a flow with infinitely many $\omega$-limit sets that are limit circuits and so are not singletons. 
\begin{figure}
\begin{center}
\includegraphics[scale=0.4]{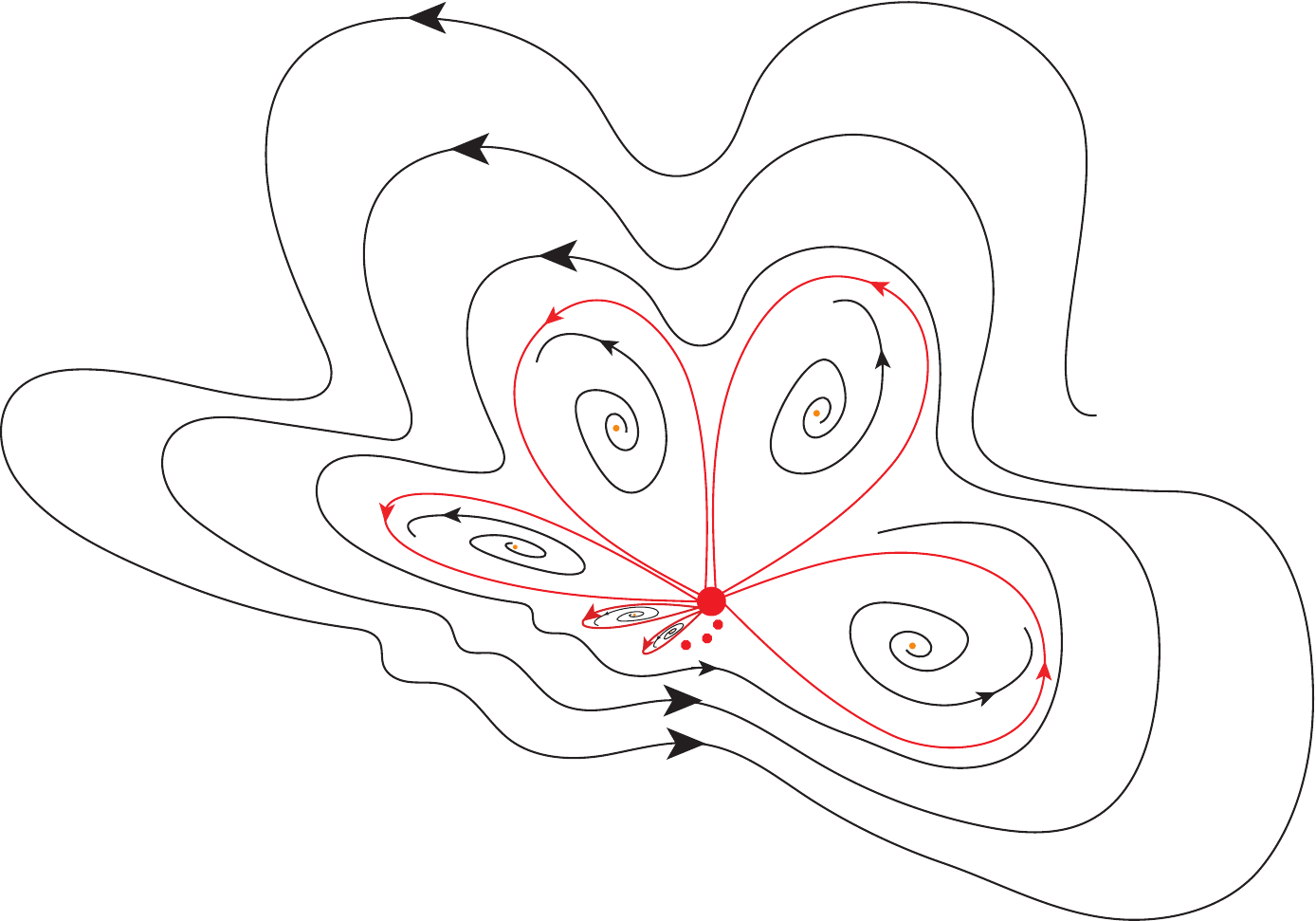}
\end{center}
\caption{Infinitely many non-degenerate $\omega$-limit set of points, which are limit circuits, and no non-degenerate $\alpha$-limit sets of points.}
\label{fig:ex_infinite_sing}
\end{figure}

Moreover, consider a flow $w$ consisting of two singular points, which are a sink and a source, and connecting orbits between them on a sphere. 
By replacing the sink with a limit circuit as in Figure~\ref{fig:ex_infinite_sing} as above, the resulting flow is a flow with infinitely many non-degenerate $\omega$-limit sets but no non-degenerate $\alpha$-limit sets of points. 
\end{proof}

The remainder of this section provides flow examples illustrating the necessity of the conditions in Corollary~\ref{main_cor:rec_th}, Corollary~\ref{main_cor:a}, and Lemma~\ref{lem:top23}. 

\subsection{Necessity of the finiteness of singular points}
The finiteness of singular points is necessary for the equivalence in Corollary~\ref{main_cor:rec_th} and Corollary~\ref{main_cor:a} as follows. 

\begin{lemma}\label{lem:counter_exam01}
There is a flow $v$ on a torus $\T^2$ with countably many singular points satisfying the following three conditions: 
\\
$(1)$ $\mathrm{R}(v) = \emptyset$ and there are neither limit circuits nor circuits with wandering holonomy.
\\
$(2)$ $\overline{\mathop{\mathrm{Cl}}(v)} = \Cv = \Sv \subsetneq \Omega(v)$.
\\
$(3)$ The flow $v$ has a non-periodic non-limit circuit in $\Omega(v)$. 
\end{lemma}

\begin{proof}
Consider a toral flow $w$ which consists of one homotopically non-trivial limit cycle $C$ and non-closed proper orbits whose $\omega$-limit sets and $\alpha$-limit sets are $C$. 
Then $\T^2 = C \sqcup \mathrm{P}(w)$. 
Fix a point $z \in C$ and a non-closed proper orbit $O$.
Write an open trivial flow box $D := \T^2 -(C \sqcup O) \subset \mathrm{P}(w)$.
Choose a closed transversal $T$ intersecting $C$ exactly once at $z$,  a point $x \in O \cap T$, and a monotonic sequence $(t_n)_{n \in \Z}$ on $O$ with $\lim_{n \to \infty} t_n = \infty$ and $\lim_{n \to - \infty} t_n = - \infty$ such that $x_n = w_{t_n}(x)$ is the $n_{\mathrm{th}}$ return of $x$ to $T$ with $\lim_{n \to \infty} x_n = \lim_{n \to - \infty} x_n = z$, and that any connected component of the set difference $D \setminus T$ is an open trivial flow box.
Denote by $D_n$ the open trivial flow box with four corners $x_n, x_{n+1}, x_{n+2}, x_{n+1}$ such that $D \setminus T = \bigcup_{n \in \Z} D_n$.
Write an open trivial flow box $D'_{2n} = D_{2n} \sqcup D_{2n + 1} \sqcup ((T \setminus O) \cap \overline{D_{2n}} \cap \overline{D_{2n+1}}) \subset D$ with four corners $x_{2n}, x_{2n+2}, x_{2n+3}, x_{2n+1}$.
By replacing the closure of each $D'_{2n}$ by an invariant square with a flow as on the left in Figure~\ref{Spiral}, and replacing $C \subset \mathop{\mathrm{Per}}(w)$ with a circuit consisting of a singular point $z$ and a non-recurrent orbit $C - \{z \}$, we obtain the resulting flow $v$ on the torus $\T^2$ such that the singular point set $\mathop{\mathrm{Sing}}(v) = \{ z \} \sqcup \{ x_n, y_n \}_{n \in \Z} \subset T$ is countable. 
\begin{figure}
\begin{center}
\includegraphics[scale=0.4]{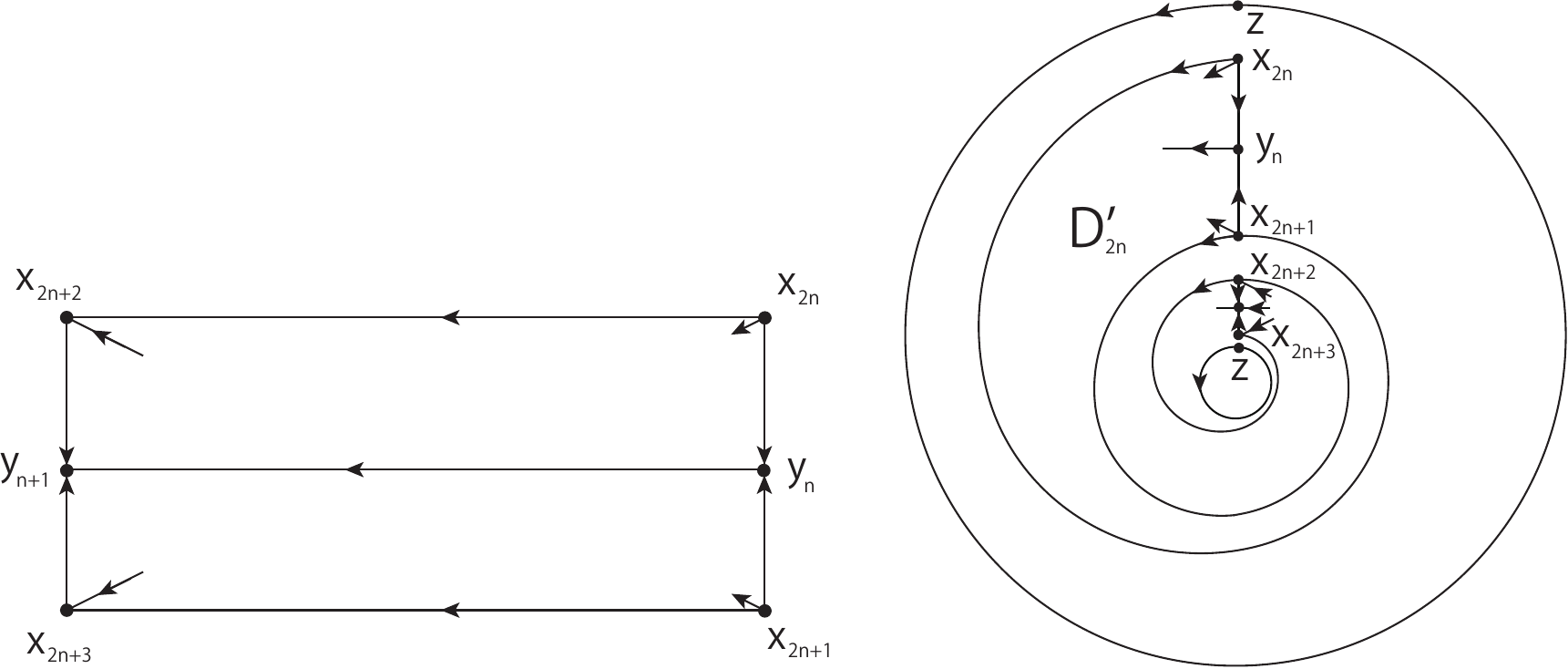}
\end{center}
\caption{An invariant square $D'_{2n}$}
\label{Spiral}
\end{figure}
Since the closure of any invariant square $D'_{2n}$ consists of singular points and non-recurrent points, we have that $\T^2 = \mathop{\mathrm{Sing}}(v) \sqcup \mathrm{P}(v)$ and $\Omega(v) = C \sqcup \{ x_n, y_n \}_{n \in \Z} \supsetneq  \mathop{\mathrm{Sing}}(v) = \overline{\mathop{\mathrm{Cl}}}(v)$.
The non-periodic circuit $C$ of $v$ is neither the $\omega$-limit set of a point nor the $\alpha$-limit set of a point and so $C$ is not a limit circuit. 
This means that there are no limit circuits of $v$. 
Moreover, since the first return maps of any closed transverse arc are orientation-preserving, there are no circuits with wandering holonomy. 
\end{proof}

\begin{lemma}\label{lem:counter_exam00}
There is a toral flow $v$ with uncountably many singular points satisfying the following two conditions: 
\\
$(1)$ $\mathrm{R}(v) = \emptyset$ and there are no circuits. 
\\
$(2)$ $\overline{\mathop{\mathrm{Cl}}(v)} = \Sv \subsetneq \Omega(v)$.
\end{lemma}

Although a flow as in the previous lemma appears in \cite[Lemma~7.1]{yokoyama2021poincare}, we provide the following construction to keep this paper self-contained.

\begin{proof}
Consider the suspension flow $w$ of a Denjoy $C^1$-diffeomorphsim $f \colon \mathbb{S}^1 \to \mathbb{S}^1$ with a minimal set $\mathcal{C} \subset \mathbb{S}^1$, which is a Cantor set, on the mapping torus $\mathbb{T}^2 := (\mathbb{S}^1 \times [0,1])/ (x,0) \sim (f(x),1)$. 
Then the minimal set $\mathcal{M}$ satisfies that $\mathcal{M} \cap (\mathbb{S}^1 \times {1/2}) = \mathcal{C} \times \{ 1/2 \} \subset S$. 
Moreover, we obtain $\Omega(w) = \mathrm{R}(w) = \mathcal{M}$. 
Since the complement $\mathbb{S}^1 - \mathcal{C}$ of the Cantor set $\mathcal{C}$ in the circle $\mathbb{S}^1$ consists of wandering domains and so consists of non-recurrent points, we have $\mathbb{T}^2 - \mathcal{M} = \mathrm{P}(w)$. 
Fix an open small \nbd $U$ of $ \mathcal{C} \times \{ 1/2 \}$. 
Let $X$ be the $C^1$ vector field that generates the flow $w$. 
Since $w$ is the suspension flow, we have that $X\vert_{\mathbb{S}^1 \times [1/3, 2/3]} = (0,1)$. 
Using a bump function $\varphi \colon \T^2 \to \R$ with $\mathbb{T}^2 -U \subseteq \varphi^{-1}(1)$ and $\varphi^{-1}(0) = \mathcal{C} \times \{ 1/2 \}$, we obtain an integrable continuous vector field $\varphi X$ which generates a flow $v_X$ such that $\mathcal{M} = \mathop{\mathrm{Sing}(v_X)} \sqcup \{ \text{separatrix of } v_X \}$ and  $O_{v_X}(y) = O_{v_X}(y)$ for any point $y \in \mathbb{T}^2 - \mathcal{M}$.
Then $\mathbb{T}^2 =  \mathop{\mathrm{Sing}}(v_X) \sqcup \mathrm{P}(v_X)$ and so $\overline{\mathop{\mathrm{Cl}}(v_X)} = \mathop{\mathrm{Sing}}(v_X)  = \mathcal{C} \times \{ 1/2 \} \subsetneq \mathcal{M} = \Omega(v_X)$. 

By Rad{\'o}'s theorem~\cite{Rado1925smooth}, there is a homeomorphism $h \colon \T^2 \to S$ to a surface $S$. 
Then the induced flow $v$ of $v_X$ by $v(t,x) := h(v_X(t,h^{-1}(x))$ is a flow on the surface $S$ which is topologically equivalent to $v_X$.
This means that $v$ has the desired properties. 
\end{proof}

\subsection{Necessity of the compactness of surfaces}

By removing singular points $\mathcal{C} \times \{ 1/2 \}$ in the flow constructed in the previous proof, we have the following statement, which implies that compactness is necessary for the equivalence in Corollary~\ref{main_cor:rec_th}. 

\begin{corollary}\label{cor:counter_exam00}
There is a non-singular flow $v$ on an infinitely many punctured torus $S$ with $\emptyset = \overline{\Cv \sqcup \mathrm{R}(v)} \neq \Omega(v) \subsetneq \mathrm{P}(v) = S$ such that there are neither virtual quasi-circuits nor virtual blow-ups of circuits with wandering holonomy.
\end{corollary}

To state the necessity of compactness for Lemma~\ref{lem:top23} and Corollary~\ref{main_cor:rec_th}, we state the following statement. 

\begin{lemma}\label{ex:counter_exam01+}
There is a flow $v$ on a non-compact surface $S$ with two ends and infinite genus without non-periodic circuits such that $\emptyset \neq \mathrm{R}(v) \subsetneq \overline{\mathop{\mathrm{Cl}}(v)} = \Omega(v)$. 
\end{lemma}

\begin{proof}
To construct such an example, recall the Denjoy construction. 
Fix an irrational number $r \in \R - \Q$, a translation $\widetilde{f_0} \colon \R \to \R$ by $\widetilde{f_0} (x) := x +r$, and an irrational rotation $f_0 \colon \mathbb{S}^1 \to \mathbb{S}^1$ on a circle $\mathbb{S}^1 := \R/\Z$ by $f_0([x]) := [x+r]$. 
Fix a sequence $(n_l)_{l \in \Z}$ with $n_0 = 0$ and $n_l < n_{l+1}$ such that $f^{n_l}([0])$ tends to $[0]$ from the positive (resp. negative) side as $l \to \infty$ (resp. $l \to -\infty$).
Consider a sequence $\{d_n := 1/(1+n^2) \mid n \in \Z \}$ and a sequence $\{I_n \mid n \in \Z \}$ of half-open intervals $I_n := [-1/(2+2n^2), 1/(2+2n^2))$ of the length $d_n$.  
Put $L := \sum_{n \in \Z} d_n = \sum_{n \in \Z} 1/(1+n^2) < \infty$ and define $y_n := \widetilde{f_0}^n(0) - \lfloor \widetilde{f_0}^n(0) \rfloor \in [0,1)$. 
Insert the interval $I_n$ at the point $y_n \in [0,1)$ and denote by $S^1 := \R/(1+L)\Z$ the resulting circle and let $f \colon S^1 \to S^1$ be a homeomorphism resulting from the insertion with a unique minimal set $\mathcal{C}$, which is a Cantor set. 
Note that the homeomorphism $f$ is not uniquely determined, as there is a choice in how each interval $I_n$ is defined.
In this case, any such choice is acceptable provided that $f$ is a homeomorphism.
Consider the suspension flow $v_f$ of the homeomorphism $f \colon \mathbb{S}^1 \to \mathbb{S}^1$ with a minimal set, which is a Cantor set.
This flow $v_f$ is defined on the mapping torus $\mathbb{T}^2 := (\mathbb{S}^1 \times [0,1])/ (x,0) \sim (f(x),1)$.
From Rad{\'o}'s theorem~\cite{Rado1925smooth} and Gutierrez's smoothing theorem~\cite[Smoothing theorem]{gutierrez1986smoothing}, we may assume that $v_f$ is $C^1$. 
Denote by $X$ the $C^1$ vector field on $\mathbb{T}^2$ which generates the $C^1$ flow $v_f$.
Let $\mathcal{M}$ be the unique minimal set of $X$ whose intersection $\mathcal{M} \cap (\mathbb{S}^1 \times \{ 1/2 \})$ is a Cantor set. 
Since $\mathbb{S}^1 \times \{1/2 \}$ is a closed transversal, by deforming $X$ into a vector field whose generating flow is topologically equivalent to $v_f$ if necessary, we may assume that the $C^1$ vector field $X$ satisfies that $X\vert_{\mathbb{S}^1 \times [1/3,2/3]} = (0,1)$.
Fix any Riemannian metric on the torus $\mathbb{T}^2$.

Define the points $0_{n,-} < 0_{n,+}$ as the points of the boundary of $I_n$. 
By construction, the sequences $(0_{n_l,-})_{l \in \Z}$ and $(0_{n_l,+})_{l \in \Z}$ of points in $S^1$ satisfy that $\lim_{ l \to -\infty} 0_{n_l,-} = \lim_{l \to -\infty} 0_{n_l,+} = 0_{0,-}$ and $\lim_{l \to \infty} 0_{n_l,-} = \lim_{l \to \infty} 0_{n_l,+}  = 0_{0,+}$. 
Fix a monotone decreasing sequence $(r_l)_{l \in \mathbb{N}}$ of points in $\mathop{\mathrm{int}} I_0 = I_0 - \partial I_0  = I_0 - \{ 0_{0,-}, 0_{0,+} \}$ with $\lim_{l \to \infty} r_l = 0_{0,-}$. 
Put $p_l := f^{n_l}(r_{\vert l \vert}) \in I_{n_l} \subset S^1$ for any $l \in \Z_{\neq0}$.
Then $\lim_{l \to -\infty} p_l = 0_{0,-}$ and $\lim_{l \to \infty} p_l  = 0_{0,+}$. 
Take pairwise disjoint closed intervals $B_l \subset I_0$ which are neighborhoods of points $r_l$ for any $l \in \mathbb{N}$. 
Then the images $D_l :=  f^{n_l}(B_{\vert l \vert}) \subset I_{n_l}$ for any $l \in \Z_{\neq 0}$ are neighborhoods of $p_l$ respectively. 
%
By deforming the first component of $X\vert_{\mathbb{S}^1 \times [1/3,2/3]} = (0,1)$ within $\mathbb{S}^1 \times (1/3, 3/8)$ and applying the inverse operation on $\mathbb{S}^1 \times (5/8, 2/3)$, we can obtain the flow by the resulting vector field that is topologically equivalent to $v_f$ generated by $X$, and so we may assume that the $C^1$ vector field $X$ satisfies that $X\vert_{\mathbb{S}^1 \times [3/8,5/8]} = (0,1)$ and that, for any $l \in \mathbb{N}$, the length of $D_l$ equals the length of $D_{-l}$ and the orientation-preserving isometry from $D_l$ to $D_{-l}$ maps $p_l$ to $p_{-l}$.
By removing a closed subset $( \{ 0_{0,-}, 0_{0,+} \} \sqcup \bigsqcup_{l \in \Z_{\neq 0}} D_ {l}) \times \{ 1/2 \} \subset \mathbb{S}^1 \times \{ 1/2 \} \subset \T^2$, we obtain the following resulting surface $T_1 := \mathbb{T}^2 - (( \{ 0_{0,-}, 0_{0,+} \} \sqcup \bigsqcup_{l \in \Z_{\neq 0}} D_ {l}) \times \{ 1/2 \})$.  

Considering the metric completion of $T_1$ with respect to the metric on $T_1 \subset \mathbb{T}^2$ induced by the Riemannian metric on $\mathbb{T}^2$, let $(T_1)_{\mathrm{mc}}$ be the metric completion  of \( T_1 \). 
By removing two points in $(T_1)_{\mathrm{mc}}$ corresponding to $\{ 0_{0,-}, 0_{0,+} \} \times \{ 1/2 \}$, the resulting space $(T_1)_{\mathrm{mc}} - \{ 0_{0,-}, 0_{0,+} \} \times \{ 1/2 \}$ can be identified with \( T_1 \) together with a countable set of arcs $D^-_l$ and $D^+_l$, where $D^-_l := \lim_{\varepsilon \to 0+} D_l \times \{ 1/2 - \varepsilon \} \subset (T_1)_{\mathrm{mc}}$ and $D^+_l := \lim_{\varepsilon \to 0+} D_l \times \{ 1/2 + \varepsilon \} \subset (T_1)_{\mathrm{mc}}$.
Then $D^+_l$ and $D^-_l$ share two endpoints.
Identifying $D^-_l$ (resp. $D^+_l$) with $D^+_{-l}$ (resp. $D^-_{-l}$) by the orientation-preserving isometry, the resulting space $S'$ is an open topological surface with two ends and infinite genus. 
For any $l \in \Z_{\neq 0}$, let $[D^+_ {l}] \subset S'$ 
be the subset obtained by identifying $D^-_l$ 
 with $D^+_{-l}$ 
Then $S' - \bigsqcup_{l \in \Z_{\neq 0}} [D^+_ {l}] = T_1 \subset \mathbb{T}^2$ has the induced smooth structure from $\mathbb{T}^2$. 

Denote by $\pi \colon T_1 \sqcup \bigsqcup_{l \in \Z_{\neq 0}} (D_ {l} \times \{ 1/2 \}) \to S'$ a non-continuous map defined by $\pi(D_l) = [D^+_ {l}]$ and $\pi \vert_{T_1} = 1_{T_1}$, where $1_{T_1}$ is the identity map on $T_1$. 
By $X\vert_{\mathbb{S}^1 \times [3/8,5/8]} = (0,1)$, we obtain the induced integrable non-singular continuous vector field $Y$ on $S'$ which is the push forward of $X$ by $\pi$.  
Using a bump function $\psi$ on $\T^2$ with $\psi^{-1}(0) = \bigsqcup_{l \in \Z_{\neq 0}} \partial D_l \times \{ 1/2 \}$ as in Figure~\ref{fig:ex_inf_genus},
\begin{figure}
\begin{center}
\includegraphics[scale=.98]{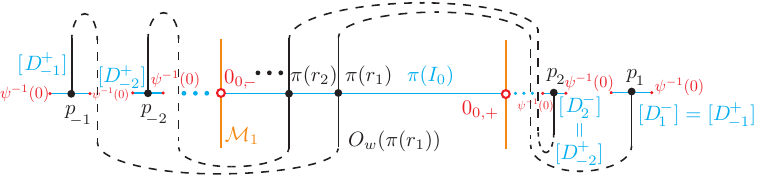}
\end{center}
\caption{Local structure near the infinite genus, including periodic orbits accumulating on $\mathcal{M}_1$.}
\label{fig:ex_inf_genus}
\end{figure} 
since $\psi \circ \pi^{-1} \colon S' \to \R$ is well-defined and piecewise smooth near $\pi(\mathbb{S}^1 \times \{ 1/2 \})$, the vector field $\psi \circ \pi^{-1} Y$ is integrable and so generates a flow $w$ on the topological surface $S'$ such that every minimal set is either $\mathcal{M}_1 := \mathcal{M} - (\{ 0_{0,-}, 0_{0,+} \} \times \{ 1/2 \}) = \mathrm{E}(w)$, a point in $\pi(\psi^{-1}(0)) = [\partial D^+_ {l}]$, or a periodic orbit in $\bigsqcup_{l \in \mathbb{N}} O_w(\pi(\mathop{\mathrm{int}}B_{l}))$. 
Then $\mathop{\mathrm{Sing}}(w) = \pi(\psi^{-1}(0)) = [\partial D^+_ {l}]$ and $\mathop{\mathrm{Per}}(w) \supseteq \bigsqcup_{l \in \mathbb{N}} O_w(\pi(r_{l}))$. 
Since $\lim_{l \to -\infty} f^{n_l}(r_{\vert l \vert}) = \lim_{l \to -\infty} p_l  = 0_{0,-}$ and $\lim_{l \to \infty} f^{n_l}(r_{\vert l \vert}) = \lim_{l \to \infty} p_l = 0_{0,+}$, 
the sequence of the length of the periodic orbits $O_w(\pi(r_{l}))$ goes to infinity. 
Therefore, by $\lim_{l \to \infty} r_l = 0_{0,-}$, the non-existence of locally dense orbits implies that $\Omega(w) = \overline{\mathop{\mathrm{Cl}}(w)} \supseteq \overline{\mathop{\mathrm{Per}}(w)} \supseteq \overline{\bigcup_{l \in \mathbb{N}} O_w(\pi(r_{l}))} \supsetneq \mathcal{M}_1 = \mathrm{E}(w) = \mathrm{R}(w) \neq \emptyset$. 

By Rad{\'o}'s theorem~\cite{Rado1925smooth}, there is a homeomorphism $h \colon S' \to S$ to a surface $S$. 
Then the induced flow $v$ of $w$ by $v(t,x) := h(w(t,h^{-1}(x))$ is a flow on the surface $S$ which is topologically equivalent to $w$.
This means that $v$ has the desired properties. 
\end{proof}

Notice that, by removing the singular point set from the surface in the flow constructed in the previous proof, we obtain a non-singular flow $v'$ on a surface with infinite ends such that $\emptyset \neq \mathrm{R}(v') \subsetneq \overline{\mathop{\mathrm{Cl}}(v')}$. This also implies the necessity of the compactness assumption for the surface in Lemma~\ref{lem:top23} and Corollary~\ref{main_cor:rec_th}. 

\subsection{Necessity of the non-existence of strictly limit circuits and circuits with wandering holonomy}

We construct a flow that demonstrates the necessity of the non-existence of strictly limit circuits and circuits with wandering holonomy in Corollary~\ref{main_cor:rec_th} and Corollary~\ref{main_cor:a}.

\subsubsection{Necessity of the non-existence of strictly limit circuits}

Now, we have the following example to state the necessity of the non-existence of strictly limit circuits in Corollary~\ref{main_cor:rec_th} and Corollary~\ref{main_cor:a}.

\begin{lemma}\label{ex:counter_exam_strict}
There is a flow $v$ with a strictly limit circuit which is generated by a vector field with four critical points on a sphere with $\overline{\mathop{\mathrm{Cl}}(v)} \subsetneq \Omega(v)$ as in Figure~\ref{fig:omega}.
\begin{figure}
\begin{center}
\includegraphics[scale=0.3]{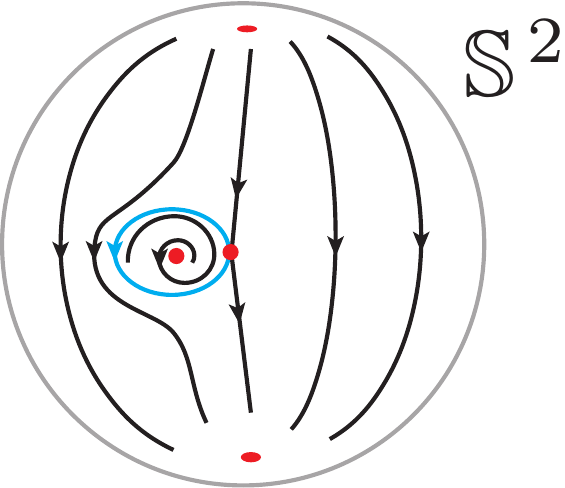}
\end{center}
\caption{A spherical flow $v$ with $\overline{\mathop{\mathrm{Cl}}(v)} \subsetneq \Omega(v)$ such that the singular point set consists of one sink, two sources, and one saddle.}
\label{fig:omega}
\end{figure} 
\end{lemma}

\begin{proof}
The singular point set $\Sv$ consists of one sink, two sources, and one saddle as in Figure~\ref{fig:omega}. 
Then, the unique homoclinic non-singular orbit of the saddle is the difference $\Omega(v) - \overline{\mathop{\mathrm{Cl}}(v)}$. 
\end{proof}

\subsubsection{Necessity of the non-existence of circuits with wandering holonomy}

Recall that a singular point $p$ of a flow $v$ on a surface is non-degenerate if the vector field $\partial v/\partial t$ is defined in a neighborhood of $p$, the Hessian matrix at $p$ is defined, and the determinant of the Hessian matrix is non-zero.
Moreover, a singular point of a flow on a surface is degenerate if it is not non-degenerate. 
The following example implies the necessity of the non-existence of circuits with wandering holonomy in Corollary~\ref{main_cor:rec_th} and Corollary~\ref{main_cor:a}.

\begin{lemma}\label{lem:counter_exam02}
There is a flow $w$ without degenerate singular points on a non-orientable closed surface $S$ with non-orientable genus four satisfying the following three conditions: 
\\
{\rm(1)} $\overline{\mathop{\mathrm{Cl}}(w)} = \mathop{\mathrm{Cl}}(w) = \mathop{\mathrm{Sing}}(w) = \bigcup_{x \in S} \omega(x) \cup \alpha(x)\subsetneq \Omega(w)$. 
\\
{\rm(2)} $\mathop{\mathrm{Per}}(w) \sqcup \mathrm{R}(w) = \emptyset$. 
\\
{\rm(3)} There are no strictly limit quasi-circuits, but there is a circuit with wandering holonomy. 
\end{lemma}

\begin{proof}
Consider a flow $v$ on a non-orientable compact surface $S_0$ as in Figure~\ref{g2-ex} with $S_0 =  \mathop{\mathrm{Sing}}(v) \sqcup \mathrm{P}(v)$ and $\bigcup_{x \in S_0} \omega(x) \cup \alpha(x) = \mathop{\mathrm{Sing}}(v) \subsetneq \mathop{\mathrm{Sing}}(v) \sqcup O' = \Omega(v)$ such that $\Omega(v)$ does not contain limit circuits and $O'$ is a proper orbit. 
Then the lift $w$ of the flow $v$ to the double $S$ of $S_0$ is a flow without degenerate singular points on the non-orientable closed surface $S$, and $w$ satisfies (1)--(3), which completes the proof. 
\begin{figure}
\begin{center}
\includegraphics[scale=0.3]{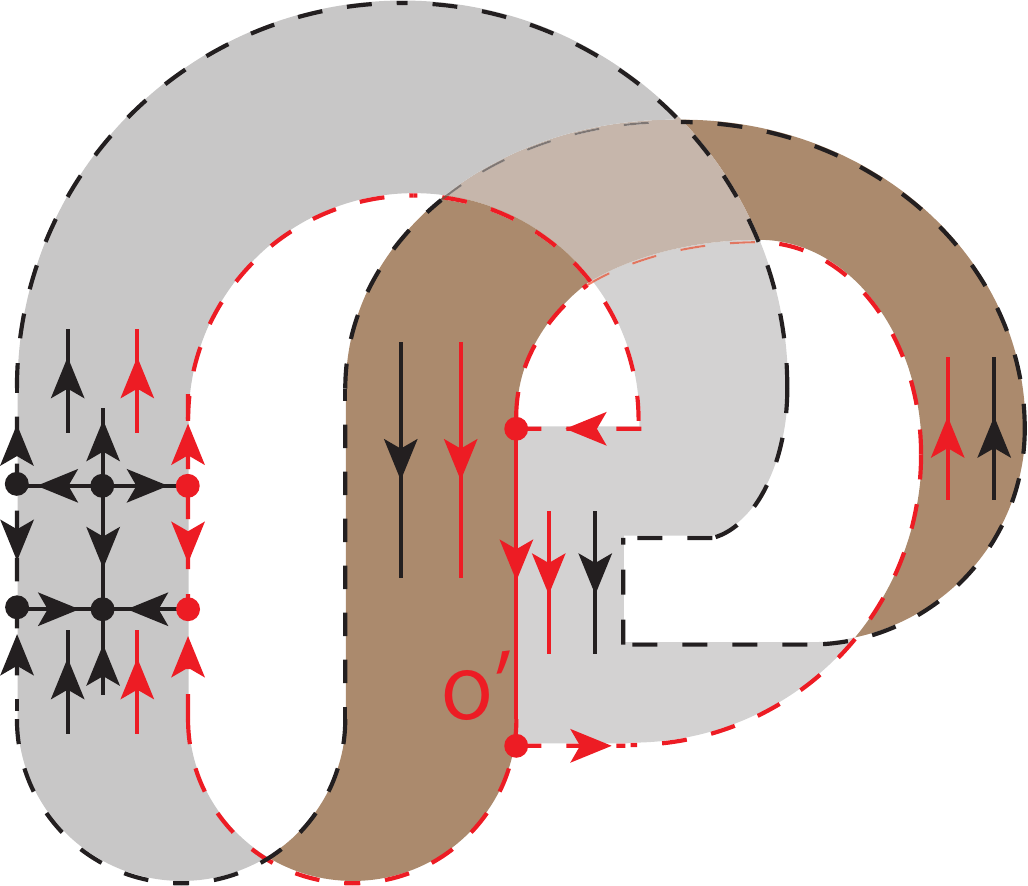}
\end{center}
\caption{A flow on a compact surface with two non-orientable genus and two boundary components which consists of one sink, one source, six $\partial$-saddles, and non-closed proper orbits.}
\label{g2-ex}
\end{figure}
\end{proof}

Recall a circuit $\gamma$ has non-orientable holonomy if there are a non-singular point $x \in \gamma$ and an arbitrarily small open transverse arc $I$ containing $x$ such that the domain of the first return map on $I$ is not empty and the first return map on $I$ is non-orientable. 

The non-existence of circuits with wandering holonomy in Corollary~\ref{main_cor:rec_th} cannot be weakened to the non-existence of circuits with non-orientable holonomy.
In fact, the denseness condition does not inhibit the existence of non-orientable holonomy. 

\begin{lemma}\label{lem:counter_exam03}
There is a flow $v$ on a Klein bottle satisfying $\overline{\mathop{\mathrm{Cl}}(v)} = \Omega(v)$ with a non-periodic circuit with fixed-point-free non-orientable holonomy. 
\end{lemma}

\begin{proof}
Consider a flow on a Klein bottle consisting of periodic orbits. 
By replacing a periodic orbit with non-orientable holonomy with a $0$-saddle with a homoclinic separatrix, we have the resulting flow $v$ which has a non-periodic circuit with fixed-point-free non-orientable holonomy and consists of one non-closed proper orbit and closed orbits such that $\mathbb{K} = \Omega(v) = \overline{\mathop{\mathrm{Cl}}(v)}$. 
\end{proof}

\section{Applications to non-wandering flows on compact surfaces}

Finally, we observe that the denseness of closed orbits for non-wandering flows on compact surfaces corresponds to the non-existence of locally dense orbits, and that any Hamiltonian flow on compact surfaces satisfies the denseness of closed orbits. 
To demonstrate this observation, we present the following characterization of non-wandering flows that captures the denseness condition.

\begin{proposition}\label{interior}
The following are equivalent for a flow $v$ on a compact surface $S$: 
\\
{\rm(1)} The flow $v$ is  non-wandering. 
\\
{\rm(2)} $\mathop{\mathrm{int}}\mathrm{P}(v) = \emptyset$. 
\\
{\rm(3)}  $\mathop{\mathrm{int}}\mathrm{P}(v) \sqcup \mathrm{E}(v) = \emptyset$ $(\mathrm{i.e.}$ $S = \mathrm{Cl}(v) \sqcup (\mathrm{P}(v) - \mathop{\mathrm{int}}\mathrm{P}(v)) \sqcup \mathrm{LD}(v) )$. 
\\
{\rm(4)} $S = \overline{\mathop{\mathrm{Cl}}(v) \sqcup \mathrm{LD}(v)}$. 
\end{proposition}

\begin{proof} 
Assertion (3) implies assertion (2). 
Recall that $S = \mathop{\mathrm{Cl}}(v) \sqcup \mathrm{P}(v) \sqcup \mathrm{R}(v)$ and that  $\mathrm{P}(v)$ is the set of non-recurrent points.  
By Lemma~\ref{lem:top23}, assertion (4) implies assertion (3).

We claim that assertion (1) implies assertion (4).
Indeed, suppose that $v$ is non-wandering. 
By \cite[Theorem III.2.12 and Theorem III.2.15]{BS}, the set of recurrent points is dense in $S$. 
The denseness of recurrent points implies that $\mathop{\mathrm{int}}\mathrm{P}(v) = \emptyset$ and so $\mathrm{E}(v) = \emptyset$.  
This implies that $S = \overline{\mathop{\mathrm{Cl}}(v) \sqcup \mathrm{LD}(v)}$. 

We claim that assertion (2) implies assertion (1).
Indeed, suppose that $\mathop{\mathrm{int}}\mathrm{P}(v) = \emptyset$. 
By Lemma~\ref{lem:top23}, we obtain that $\overline{\mathop{\mathrm{Cl}}(v) \sqcup \mathrm{LD}(v)} \cap \mathrm{E}(v) = \emptyset$. 
From $S = \mathop{\mathrm{Cl}}(v) \sqcup \mathrm{P}(v) \sqcup \mathrm{LD}(v) \sqcup \mathrm{E}(v)$, the complement $U := S - \overline{\mathop{\mathrm{Cl}}(v) \sqcup \mathrm{LD}(v)} = \mathop{\mathrm{int}}(\mathrm{P}(v) \sqcup \mathrm{E}(v))$ is an open subset of $S$. 
By Lemma~\ref{lem:top23}, this open set $U$ contains $\mathrm{E}(v)$ and thus is an open neighborhood of $\mathrm{E}(v)$.
Assume that $\mathrm{E}(v) \neq \emptyset$. 
Then $U = \mathop{\mathrm{int}}(\mathrm{P}(v) \sqcup \mathrm{E}(v))$ is nonempty. 
Since the interior of the Q-set of any exceptional orbit is empty, Lemma~\ref{lem:Maier} implies that the closure $\overline{\mathrm{E}(v)}$ is a finite union of Q-sets of exceptional orbits and therefore has the empty interior. 
Then the set difference $\mathop{\mathrm{int}}(\mathrm{P}(v) \sqcup \mathrm{E}(v)) \setminus \overline{\mathrm{E}(v)} \subseteq \mathop{\mathrm{int}}\mathrm{P}(v)$ is a nonempty open subset, which contradicts $\mathop{\mathrm{int}}\mathrm{P}(v) = \emptyset$. 
Thus, we have that $\mathrm{E}(v) = \emptyset$. 
Then $S = \mathop{\mathrm{Cl}}(v) \sqcup (\mathrm{P}(v) - \mathop{\mathrm{int}}\mathrm{P}(v)) \sqcup \mathrm{LD}(v)$. 
Therefore, the closure of the set of recurrent points is the whole surface (i.e. $S = \overline{\mathop{\mathrm{Cl}}(v) \sqcup \mathrm{LD}(v)}$) and so $v$ is non-wandering. 
\end{proof}

By the previous proposition, notice that the existence of exceptional quasi-minimal sets implies the existence of wandering domains. 
We have the following equivalence. 

\begin{lemma}\label{lem-01}
The following are equivalent for a flow $v$ on a compact surface $S$:
\\
$(1)$ The flow $v$ is non-wandering and $\overline{\mathop{\mathrm{Cl}}(v)} = \Omega(v)$.
\\
$(2)$ $S = \overline{\mathop{\mathrm{Cl}}(v)}$.
\\
$(3)$ $\mathop{\mathrm{int}} \mathrm{P}(v)  = \emptyset$ and $\overline{\mathop{\mathrm{Cl}}(v)} = \Omega(v)$.
\\
$(4)$ 
$S = \mathop{\mathrm{Cl}}(v) \sqcup (\mathrm{P}(v) - \mathop{\mathrm{int}}\mathrm{P}(v))$. 
\\
$(5)$ The flow $v$ is non-wandering and $\mathrm{LD}(v) = \emptyset$.
\end{lemma}

\begin{proof}
Recall that $S = \mathop{\mathrm{Cl}}(v) \sqcup \mathrm{P}(v) \sqcup \mathrm{R}(v)$.
Obviously, assertions $(1)$ and $(2)$ are equivalent, and assertion $(4)$ implies assertion $(2)$.
Lemma~\ref{thm41} implies that assertion $(3)$ implies assertion $(4)$.

We claim that assertion $(2)$ implies assertion $(3)$.
Indeed, if $\overline{\mathop{\mathrm{Cl}}(v)} = S$, then $\mathop{\mathrm{int}} \mathrm{P}(v)  = \emptyset$ and $\Omega(v) \subseteq S =  \overline{\mathop{\mathrm{Cl}}(v)} \subseteq \Omega(v)$ due to the definition of $\Omega(v)$. 

Notice that assertions $(1)$--$(4)$ imply assertion $(5)$. 
Proposition~\ref{interior} implies assertion (5) implies assertion $(2)$. 
\end{proof}

In the non-wandering case, the denseness of closed orbits is characterized as follows. 

\begin{corollary}\label{prop:nw}
The following are equivalent for a non-wandering flow $v$ on a compact surface:
\\
$(1)$ $\overline{\mathop{\mathrm{Cl}}(v)} = \Omega(v)$.
\\
$(2)$ $\mathrm{LD}(v) = \emptyset$.
\end{corollary}

The non-wandering property in the previous corollary is necessary. 
In fact, the Denjoy flow has no locally dense orbits, but the non-wandering set consists of non-closed recurrent orbits. 
This implies the following observation. 

\begin{corollary}
Each Hamiltonian flow on a compact surface satisfies the denseness of closed {\rm(}in particular, recurrent and non-wandering{\rm)} orbits. 
\end{corollary}

\begin{proof}
Let $v$ be a Hamiltonian flow on a compact surface. 
Since any Hamiltonian flow on a compact surface is area-preserving and so has no wandering domain, the flow $v$ is non-wandering. 
The existence of the Hamiltonian of $v$ on the surface implies the non-existence of non-closed recurrent orbits. 
Corollary~\ref{prop:nw} implies the assertion. 
\end{proof}

The compactness in the previous corollary is necessary. 
In fact, the flow generated by a vector field $X = (1,0)$ on the plane $\R^2$ is Hamiltonian but consists of non-recurrent wandering points. 

\section*{Acknowledgments}
The author acknowledges the hospitality of the Fields Institute for Research in Mathematical Sciences. 

\bibliographystyle{abbrv}
\bibliography{yt20260505}

\end{document}